\documentclass[11pt]{article}

\usepackage[utf8]{inputenc}
\usepackage{lmodern}
\usepackage{subfiles}
\usepackage{enumitem}
\setenumerate{topsep=6pt,ref={\normalfont(\roman*)},label={\normalfont(\roman*)}, itemsep=0pt}
\usepackage{pgfplots}
\pgfplotsset{compat=1.18}
\usepgfplotslibrary{groupplots}
\usepackage{amsfonts}
\usepackage{amsthm}
\usepackage{amsmath}
\usepackage{amssymb}
\usepackage{amscd}
\usepackage{mathrsfs}
\usepackage{mathtools}
\usepackage{bbm}
\usepackage{esint}

\usepackage[margin=3cm]{geometry}
\usepackage{setspace}
\usepackage{indentfirst}
\usepackage{graphicx}
\usepackage{graphics}
\usepackage{lscape}
\usepackage{pgf,tikz}
\usepackage{tikz-cd}
\usepackage{color}
\usepackage{pict2e}
\usepackage{epic}
\usepackage{epstopdf}
\usepackage{titlesec, titlefoot}
\titleformat{\section}[block]{\Large\bfseries\filcenter}{\thesection}{1em}{}
\usepackage{commath}
\usepackage{float}
\usepackage{caption}
\usepackage{etoolbox}
\usepackage[affil-it]{authblk}
\usepackage{combelow}

\usepackage[hidelinks, bookmarksdepth=3]{hyperref}
\hypersetup{bookmarksopen=true} 
\usepackage{hypcap}

\graphicspath{{./Pictures/}}
\allowdisplaybreaks

\expandafter\def\expandafter\normalsize\expandafter{%
\normalsize
\setlength\abovedisplayskip{6pt}
\setlength\belowdisplayskip{6pt}
\setlength\abovedisplayshortskip{6pt}
\setlength\belowdisplayshortskip{6pt}
}

\usetikzlibrary{arrows.meta}

\definecolor{profileblue}{RGB}{36,104,154}
\definecolor{errororange}{RGB}{181,100,42}
\definecolor{figuregray}{RGB}{151,157,162}

\theoremstyle{plain}

\renewcommand*\thesection{\arabic{section}}
\numberwithin{equation}{section} 

\newtheorem{theorem}{Theorem}[section]
\newtheorem{lemma}[theorem]{Lemma}
\newtheorem*{lemma*}{Lemma}
\newtheorem{proposition}[theorem]{Proposition}
\newtheorem{corollary}[theorem]{Corollary}
\newtheorem{conjecture}[theorem]{Conjecture}

\theoremstyle{definition}
\newtheorem{definition}[theorem]{Definition}
\newtheorem{remark}[theorem]{Remark}

\expandafter\let\expandafter\oldproof\csname\string\proof\endcsname
\let\oldendproof\endproof
\renewenvironment{proof}[1][\proofname]{%
\oldproof[\upshape \bfseries #1]%
}{\oldendproof}

\makeatletter
\def\@makechapterhead#1{%
\vspace*{50\p@}%
{\parindent \z@ \raggedright \normalfont
\interlinepenalty\@M
\Huge\bfseries  \thechapter.\quad #1\par\nobreak
\vskip 40\p@
}}
\makeatother

\newcommand{\eps}{\varepsilon}

\DeclareMathOperator{\ddiv}{div}
\DeclareMathOperator{\supp}{supp}

\DeclareMathOperator \dist{dist}

\DeclareMathOperator{\ind}{ind}

\DeclareMathOperator{\curl}{curl}

\def \a{\alpha}

\def \p{\partial}
\def \mc{\mathcal}
\def \mb{\mathbb}

\begin{document}

\title{\textbf{Magnetic vector fields with isolated singularities}}

	\author[1]{{\Large L\'aszl\'o Sz\'ekelyhidi}}
    \author[1]{{\Large David Ziener}}

    \affil[1]{\small Max Planck Institute for Mathematics in the Sciences,
    Inselstrasse 22, 04103 Leipzig, Germany
    \protect\\
    \texttt{laszlo.szekelyhidi@mis.mpg.de}
    \protect\\
    \texttt{david.ziener@mis.mpg.de}}

	\date{}

    \maketitle

	\begin{abstract}
    We construct unit-length, divergence-free vector fields tangent to the boundary on arbitrary smoothly bounded domains in $\mathbb R^3$. These vector fields are continuously differentiable except for a finite number of points on the boundary. Our approach identifies a form of h-principle and a close analogy with Nash's construction of smooth isometric embeddings.
    \end{abstract}

    \medskip

    \noindent\textbf{Mathematics Subject Classification (2020).} Primary 35F60; Secondary  58C15.

    \smallskip

    \noindent\textbf{Keywords and phrases.}
    Unit-length divergence-free vector fields, stray-field-free magnetizations, micromagnetics, isolated singularities, Nash--Moser iteration.

	\vspace{1cm}

\section{Introduction}
Suppose $U \subset \mb R^{3}$ is a bounded body, which models some ferromagnetic material with associated magnetization $\boldsymbol{m}:U \to \mb R^3$. To such a magnetization which we suppose also satisfies the saturation condition
\[
|\boldsymbol{m}|=1,\quad \text{in }U,
\]
we can associate the micromagnetic energy \cite{Brown} which in the setting without an external field and in non-dimensional variables takes the form (see \cite{HubertSchafer1998,DeSimoneKohnMullerOtto2000})
\begin{equation}
    \label{eq:MagneticEnergy}
    E(\boldsymbol{m})=w^2 \int_U |\nabla \boldsymbol{m}|^2dx+Q\int_{U}\varphi(\boldsymbol{m})dx+\int_{\mb R^3}|H_d|^2\,dx.
\end{equation}
The three different energy contributions are called exchange, anisotropy and magnetostatic (or stray-field) energy. Local minima of \eqref{eq:MagneticEnergy} represent the stable magnetizations of a magnetic body.

The different terms of the energy have the following roles: The exchange energy favors a slowly varying magnetization and sets a length scale of order $w$ which is called the \textit{magnetostatic exchange length}. The anisotropy energy favors certain directions of $\boldsymbol{m}$ according to the zero set of the smooth function $\varphi:\mb S^2 \to \mb R_+$. The \textit{quality factor} $Q$ is a material parameter which determines the weight of the anisotropy energy with respect to the magnetostatic energy. A material is called hard if $Q>1$, soft if $Q<1$ and ideally soft if $Q=0$. Finally, the magnetostatic energy tries to eliminate the distributional divergence. Indeed, the demagnetizing field $H_d$ is related to $\boldsymbol{m}$ through Maxwell's equations
\begin{equation}
    \label{eq:DemagField}
       \curl(H_d) = 0, \quad \ddiv(H_d + \boldsymbol{m} \mathbf{1}_U) = 0 \quad \text{in } \mathcal{D}'(\mathbb{R}^3), \qquad H_d(x) \to 0\quad \text{as } |x| \to \infty.
\end{equation}
Due to the competition of the different energy terms in the micromagnetic energy, one can observe the separation of a magnetic body into regions of slowly varying magnetizations, so-called \textit{domains} separated by \textit{domain walls} in which the magnetization rapidly changes. One major goal in micromagnetics is understanding the properties of these domains and their domain walls, see \cite{HubertSchafer1998} and its references.

In this paper we are interested in ideally soft ferromagnetic materials \cite[Section 3.3.3]{HubertSchafer1998}, that is, $Q=0$. Moreover, we suppose that the size of our magnetic body is much larger than the magnetostatic exchange length such that the exchange energy is negligible. Hence, directly from \eqref{eq:MagneticEnergy}--\eqref{eq:DemagField}, in this idealized setting a stable magnetization $\boldsymbol{m}$ on $U$ should satisfy
\begin{equation}
    \label{eq:StrayFieldFree}
    |\boldsymbol{m}|=1,\quad \ddiv \boldsymbol{m}=0,\quad \text{in }U,\qquad \boldsymbol{m} \cdot \nu =0,\quad \text{on }\p U,
\end{equation}
where $\nu$ is the outer unit normal to $U$. In the following we will call such a vector field a \textit{magnetic vector field}. In fact due to the absence of any anisotropy it is not clear if one should expect the formation of domains in this case. It is instructive to consider first the two-dimensional situation, which corresponds to a thin-film element. Suppose $\omega \subset \mb R^2$ is some bounded simply connected region. In the idealized setting of large sample size, zero anisotropy and without demagnetizing field Van den Berg \cite{VanDenBerg1986} argued that a magnetization $\boldsymbol{m}=\boldsymbol{m}(x,y)$ should satisfy
\begin{equation}
    \label{eq:StrayFieldFree2D}
    \boldsymbol{m}_z=0,\quad |\boldsymbol{m}|=1,\quad \ddiv \boldsymbol{m}=0,\quad \text{in }\omega,\qquad \boldsymbol{m} \cdot \nu =0,\quad \text{on }\p \omega.
\end{equation}
In fact to solve \eqref{eq:StrayFieldFree2D}, one can equivalently write $\boldsymbol{m}=\nabla^\perp v$ and solve the eikonal equation
\begin{equation}
    \label{eikonal}
    |\nabla v|=1,\quad \text{in }\omega,\qquad v =\mathrm{const},\quad \text{on }\p\omega.
\end{equation}
Van den Berg uses the method of characteristics to obtain solutions to \eqref{eikonal}. Since the characteristics may intersect, $\nabla v$ will in general exhibit line discontinuities, which one can interpret as domain walls with thickness shrunken to $0$. It is simple to see \cite[Section~6]{DeSimoneKohnMullerOtto2006} that these solutions are nonunique. In particular it is an interesting question to select physically relevant solutions by comparing their micromagnetic energies. Since the singular domain walls in van den Berg's construction should be replaced with continuously varying magnetization on a microscopic level, this also includes studying the energy of different wall patterns. In soft ferromagnetic films the two simplest wall patterns are the following. First, the Bloch wall in which the magnetization varies only in a direction $n$, performing a $180^\circ$ rotation in a plane perpendicular to $n$ \cite[Section~3.6.1(A)]{HubertSchafer1998}. Second, the Néel wall where the magnetization rotates in the film-plane \cite[Section~3.6.4(C)]{HubertSchafer1998}. Another more complicated wall-type is the cross-tie wall, which consists of a periodic arrangement of vortices. Here one should mention the remarkable paper by Alouges-Rivière-Serfaty \cite{AlougesRiviereSerfaty2002} which explains analytically the formation of these structures in thin films and derives regimes in which a cross-tie wall is energetically preferable to a Néel wall. There is a vast amount of mathematical literature studying different wall structures in micromagnetics and it is impossible to give an exhaustive list. We again refer to the survey article \cite{DeSimoneKohnMullerOtto2006} and its references. See also Rivière's survey \cite{Riviere2002}.

In the setting of a three-dimensional bulk region when solving \eqref{eq:StrayFieldFree} one should not distinguish between domains and domain walls. This can be illustrated by the following surprising construction. Arrott-Heinrich-Aharoni \cite{ArrottHeinrichAharoni1979} (see also \cite{HubertRave1998}) construct a solution to \eqref{eq:StrayFieldFree} on a cylinder which possesses only two point singularities at the top and bottom of the cylinder. Let us briefly explain their Ansatz: In cylindrical coordinates $(\rho,\varphi,z)$ denote by $e_\varphi$ the unit vector in the $\varphi$-direction. Consider a potential of the form $A=A_\varphi(\rho,z)e_\varphi$ and let
\begin{equation}
    \label{eq:Arrott}
     \boldsymbol{m}=\operatorname{curl}(A_\varphi e_\varphi)
    +\sqrt{1-\bigl|\operatorname{curl}(A_\varphi e_\varphi)\bigr|^2}
    \,e_\varphi.
\end{equation}
If we choose $A_\varphi$ so that it vanishes on the boundary of the cylinder and $|\operatorname{curl}(A_\varphi e_\varphi)|\leq 1$, one easily sees that $\boldsymbol{m}$ satisfies \eqref{eq:StrayFieldFree}. To prevent singularities in $\boldsymbol{m}$ forming on the axis, $A_\varphi$ has to decay suitably towards the axis. For example, an explicit choice within the Ansatz of \cite[p. 162, equation (3.78a)]{HubertSchafer1998} is
\[
A_\varphi(\rho,z)=\frac{\rho(1-\rho^2/R^2)\,z(L-z)}{2\sqrt{z^2(L-z)^2+L^2\rho^2}},
\]
where we denote by $L$ the length of the cylinder, by $R$ its radius, and let $0\leq z\leq L$.\\
Through this observation Hubert and Schäfer \cite[p.~163]{HubertSchafer1998} posed the following conjectures, which ask if the construction in \cite{ArrottHeinrichAharoni1979} can be generalized to arbitrary three-dimensional bodies.
\begin{conjecture}
    \label{conj:Con1}
    ``For all non-pathological bodies that we define as bodies with only a finite number of corners, there exist magnetization distributions $\boldsymbol{m}(\boldsymbol{r})$ with $|\boldsymbol{m}|=1$, $\operatorname{div}\boldsymbol{m}=0$ and $\boldsymbol{m}\cdot\boldsymbol{n}=0$ ($\boldsymbol{n}=\text{surface normal}$), which are continuous and differentiable everywhere except in a finite number of singular points.''
\end{conjecture}
A stronger version of the conjecture states the following \cite[p.~163]{HubertSchafer1998}.
\begin{conjecture}
    \label{conj:Con2}
    ``All these singular points can be placed at the surface of the body.''
\end{conjecture}
Note that since a magnetic vector field should be unit length and tangent to the boundary by the Hopf Index Theorem we cannot expect a magnetic vector field to be free of singularities in general.

To the authors' knowledge, the only construction of magnetic vector fields in three dimensions which possess a finite number of point singularities is the construction in \cite{ArrottHeinrichAharoni1979}, which through the Ansatz \eqref{eq:Arrott} can easily be generalized to axially symmetric bodies. See also the paper of Vlaming and van den Berg \cite{VlamingVanDenBerg1988} in which a possible strategy by foliations for more general three-dimensional bodies is proposed, but not carried out in a mathematically rigorous way.

In Theorem \ref{thm:main} below, we provide in the setting of bodies without corners a very general answer to both conjectures, by constructing magnetic vector fields on arbitrary smoothly bounded bodies in $\mathbb R^3$, with the topologically minimal number of singularities.
\subsection{Mathematical preliminaries}
We fix once and for all a smoothly bounded domain\footnote{Following the mathematical convention a domain is an open connected set. This is not to be confused with the notion of magnetic domains we discussed previously.} $U\subset\mathbb R^3$ and denote by $\nu$ its outer unit normal. For an open subset $E\subset\overline U$ and $k\in\mathbb N_0$, we write $C^k(E;\mathbb R^3)$ for the space of continuous vector fields on $E$, whose derivatives up to order $k$ on $E \cap U$ extend continuously to $E$. Each boundary component $S$ of $U$ is a closed oriented surface and we denote its Euler characteristic by $\chi(S)$. Moreover, we will write from now on $X$ for (magnetic) vector fields instead of $\boldsymbol{m}$.

To talk more precisely about singularities, we introduce the index. Let $p\in S$ and suppose $W$ is a neighborhood of $p$. Let  $T \in C^1(W\setminus\{p\};\mathbb R^3)$ be tangent to $S$ and nonvanishing on $W \setminus \{p\}$. Choose a closed disk $D\subset W$ with $p$ in its interior and a positively oriented orthonormal frame $(e_1,e_2)$ defined on $D$. Write
\[
T=T^1\,e_1+T^2\,e_2
\qquad\text{on }D\setminus\{p\}.
\]
The \textit{index} of $T$ at $p$ is
\[
\operatorname{ind}_p(T)
:=
\deg\left(\frac{(T^1,T^2)}{\sqrt{(T^1)^2+(T^2)^2}}: \partial D\longrightarrow\mathbb S^1\right)
\]
where $\partial D$ carries the boundary orientation. In other words the index is the winding number of $T$ around $p$. It is easy to see that this definition is independent of the choice of $D$ and the frame $(e_1,e_2)$.\\
Denote by $S_1,\ldots,S_\ell$ the connected components of $\partial U$ with
$\chi(S_i)\neq 0$. Choose $p_i\in S_i$ for $i=1,\ldots,\ell$, and set
\begin{equation}
    \label{def:Sigma}
    \Sigma=\{p_1,\ldots,p_\ell\}.
\end{equation}
From now on the set $\Sigma$ is fixed throughout the paper. Moreover, we fix once and for all $0<\delta<1/2$. Our main result is the following theorem.
\begin{theorem}
    \label{thm:main}
    There exists $X\in H^1(U;\mb R^3)\cap C^1(\overline U\setminus\Sigma;\mathbb R^3)$ such that
    \[
    \ddiv X=0\quad\text{in }U,\qquad |X|=1\quad\text{on }\overline U\setminus\Sigma,\qquad
    X\cdot\nu=0\quad\text{on }\partial U\setminus\Sigma.
    \]
    Moreover, for every $i=1,\ldots,\ell$, we have
    \[
        \ind_{p_i}(X|_{S_i})=\chi(S_i)
    \]
    and
    \[
        |DX(x)|=O(|x-p_i|^{-1-\delta})
        \quad\text{as }x\to p_i,\quad
        x\in\overline U\setminus\Sigma.
    \]
\end{theorem}
In fact by looking at the proof of Theorem \ref{thm:main}, we have that $X=\curl B$, for some $B \in C^2(\overline{U}\setminus\Sigma;\mb R^3)$.\\

\begin{remark}
Our proof gives much more than just the existence of a single vector field $X$ as above. First, observe that, if one replaces the constraint ``$|X|=1$'' with ``$X$ is nowhere vanishing'', the system becomes extremely flexible. To be more precise, in Lemma \ref{lem:RelativeExactExtension} we prove a \textit{relative h-principle} for nowhere vanishing divergence free vector fields. See \cite{EliashbergMishachev2002} for the general framework.

If one aims to construct a unit-length divergence free $X$, flexibility of the system is known in the low-regularity regime of $L^\infty$ vector fields see \cite[Section~4.5.2]{Szekelyhidi2013}. Here one can make a direct connection to the constructions of $C^{1,\alpha}$ isometric embeddings \cite{ContiDeLellisSzekelyhidi2012,CaoHirschInauen2025}. In particular one could prove a version of Lemma \ref{lem:RelativeExactExtension} and construct for any $\alpha<1/3$ a unit-length divergence free $X \in C^{0,\alpha}$.

To build a unit-length divergence free field of higher regularity we follow the idea for constructing smooth isometric embeddings \cite{Nash1956}: We first build a nonvanishing divergence free vector field which is almost unit-length and then remove the length error with a Nash--Moser correction. In particular, due to the flexibility for nonvanishing divergence free vector fields described above, one can view our result as an h-principle for $C^1$ divergence free vector fields.
\end{remark}

\begin{remark}
Hubert and Schäfer \cite[p.~163]{HubertSchafer1998} suggest that, if solutions to \eqref{eq:StrayFieldFree} with finitely many point singularities exist, they are ``probably energetically favoured compared to hypothetical solutions
with line singularities or with an infinite number of point singularities''. At the moment we only know that our solutions have finite exchange energy. Understanding whether or not they are energetically optimal in some sense will be part of future work.
\end{remark}

\subsection{Overview of the proof}
\label{sec:overview}
The proof of Theorem \ref{thm:main} relies fundamentally on the ideas of Nash in his construction of smooth isometric embeddings \cite{Nash1956}. To make this connection clear, we start by briefly sketching the main strategy of Nash's proof. For an excellent exposition we refer to \cite{Andrews2002}.
\subsubsection*{Nash's construction}
 The isometric embedding problem is to find, for a given metric $g=(g_{ij})$ and some sufficiently large  $q\in\mathbb N$, a smooth embedding  $u$ into $\mathbb R^q$ whose derivatives satisfy, in local coordinates,
\begin{equation}
    \label{eq:overview-isometry}
    (g_u)_{ij}:=\p_i u\cdot\p_j u=g_{ij}.
\end{equation}
Here $g_u=(Du)^T Du$ is the induced metric. Nash looks for an embedding of the form $u=(y,z)$ which means that the induced metric is of the form
\begin{equation}
    \label{eq:overview-metric-splitting}
    g_{(y,z)}=g_y+g_z.
\end{equation}
The purpose of $y$ will be to approximate the given metric $g$ by oscillations, while $z$ allows for a correction of the remaining error.\\
First, one constructs a fixed smooth embedding $z$ which is strictly short (meaning that $g-g_z$ is positive definite), and which is \textit{free}. Freeness is the nondegeneracy condition that will ensure the solvability of the linearized equation in Step 2 below.\\
Moreover, Nash obtains a finite decomposition
 \begin{equation}
    \label{eq:overview-primitive-metrics}
    g-g_z=\sum_\nu \rho_\nu^2
    \nabla\varphi_\nu\otimes\nabla\varphi_\nu,
\end{equation}
with smooth $\rho_\nu$ and $\varphi_\nu$.\\
\textbf{Step 1: Nash Spirals.}
Nash uses an Ansatz of the form
\begin{equation}
    \label{eq:overview-twist}
    y_\lambda=
    \left(
    \frac{\rho_\nu}{\lambda}\cos(\lambda\varphi_\nu),
    \frac{\rho_\nu}{\lambda}\sin(\lambda\varphi_\nu)
    \right)_\nu.
\end{equation}
Each component pair contributes
 $\rho_\nu^2\nabla\varphi_\nu\otimes\nabla\varphi_\nu
+\lambda^{-2}\nabla\rho_\nu\otimes\nabla\rho_\nu$
to the induced metric. Thus \eqref{eq:overview-metric-splitting}--\eqref{eq:overview-primitive-metrics} and a simple computation give
 \begin{equation}
    \label{eq:overview-metric-error}
    g_{(y_\lambda,z)}
    =g+\frac{1}{\lambda^2}
    \sum_\nu \nabla\rho_\nu\otimes\nabla\rho_\nu.
\end{equation}
The oscillatory terms cancel in this computation. Consequently, the metric error tends to zero in $C^k$ as $\lambda\to\infty$, even though the higher derivatives of $y_\lambda$ may become large. The important point is that $z$ remains fixed throughout this step. Hence the error in \eqref{eq:overview-metric-error} can be made arbitrarily small without changing the free embedding about which the final correction will be performed.\\
\textbf{Step 2: Nash--Moser correction.}
Set $h=g-g_{(y_\lambda,z)}$. Keeping $y_\lambda$ fixed, we seek a small perturbation $v$ of $z$ such that
\begin{equation}
    \label{eq:overview-metric-correction}
    g_{z+v}=g_z+h.
\end{equation}
Then $u=(y_\lambda,z+v)$ satisfies \eqref{eq:overview-isometry}. To solve \eqref{eq:overview-metric-correction} the first idea is to use the implicit function theorem. The linearization of \eqref{eq:overview-metric-correction} at $v=0$ is
 \begin{equation}
    \label{eq:LinearizationIsometricEmbeddings}
    \p_i z\cdot\p_j v
    +\p_j z\cdot\p_i v=h_{ij}.
\end{equation}
For this linearized problem, a computation shows that, under the assumption  $\p_i z\cdot v=0$, the equation reduces to a linear system, which is solvable under the freeness condition. Hence we obtain a right inverse for the linearized metric equation. However, here Nash faces a problem: The map
\[
v\mapsto g_{z+v}-g_z
\]
maps $C^{k+1} \rightarrow C^k$. The right inverse obtained from freeness on the other hand gives only
\[
||v||_{C^k}\lesssim ||h||_{C^k}.
\]
This \textit{loss of derivatives} prevents an application of the standard implicit function theorem. Nash overcomes this difficulty by carefully using smoothing operators in the correction procedure, which allow him to control higher derivatives while the metric error tends to zero. For sufficiently small $h$, this gives the desired solution of \eqref{eq:overview-metric-correction}. Since $z$ is fixed, the required smallness can be achieved by increasing $\lambda$.

This combination of linearized correction and smoothing was subsequently developed into a general method, notably in the work of Moser \cite{Moser1961,Moser1966}, Zehnder \cite{Zehnder1975}, and Hamilton \cite{Hamilton1982}. The resulting Nash--Moser theorems are applicable far beyond the isometric embedding problem. These include for example applications in KAM theory \cite{Zehnder1975}, the construction of solutions of nonlinear wave equations \cite{CraigWayne1993}, and free-boundary problems in physical geodesy \cite{Hormander1976} and water-wave theory \cite{Lannes2005}. More geometric applications are deformations of complex structures on manifolds with boundary \cite{Hamilton1977} and local normal forms for smooth Poisson structures \cite{Conn1985}.
\subsubsection*{Our construction}
We can now describe the strategy behind the proof of Theorem \ref{thm:main}. Our goal is to construct  $X\in C^1(\overline U\setminus \Sigma;\mathbb R^3)$, such that
\begin{equation}
    \label{eq:DivFreeUnit}
    \ddiv X=0\quad\text{in }U,\qquad |X|=1\quad\text{on }\overline U\setminus\Sigma,\qquad X\cdot\nu=0\quad\text{on }\partial U\setminus\Sigma.
\end{equation}
Our idea is to use the unit length constraint as the nonlinear equation, while ensuring that the divergence remains zero through adding only perturbations of the form $\curl B$, where $B \in C^{2}(\overline U\setminus \Sigma;\mathbb R^3)$.

As in Nash's construction, we carry out the proof by using his two fundamental ideas. First, through a finite number of convex integration perturbations we can make a vector field arbitrarily close to unit-length while preserving its divergence. Second, a Nash--Moser theorem removes the final length error. We outline our strategy.

Before the two main ideas can be used we have to do some preliminary work. We have to find a vector field $X_0 \in C^\infty(\overline U\setminus \Sigma;\mathbb R^3)$, which allows us to add the oscillatory perturbations and, specific to our problem, already has the prescribed indices at the points of $\Sigma$. To do so, we introduce in Section \ref{sec:Prelim} the necessary analytic tools which will be used throughout the paper. Then in Section \ref{sec:Prep} we carry out the construction of $X_0=\curl B_0$. The main property of this  vector field is that it already solves the full problem in suitably scaled neighborhoods of the boundary while being strictly short in the remaining part of  $\overline U \setminus \Sigma$. A technical detail in the construction is that the error  $\Delta:=1-|X_0|^2$ has to satisfy suitable derivative bounds. This ensures that the perturbations we pick in the proof of the main theorem extend smoothly, so $X_0$ remains unchanged close to the boundary. See Sections \ref{sec:LocPert} and \ref{sec:ProofMain}.

\textbf{Building block.} 
We introduce the basic building block. Suppose $\Omega \subset \mb R^3$ is open and $X \in C^\infty(\Omega;\mb R^3)$ is nowhere vanishing. As an initial step we have to understand how to realize changes in a vector field through adding a curl. The usual Ansatz is to consider $B(x,\vartheta) \in \mb R^3$, smooth in $x\in \Omega$ and $2\pi$-periodic in $\vartheta$. For a unit vector $\xi \in \mb R^3$ and $\lambda>0$ we compute from the chain rule
\begin{align}
    X+\curl\left(\frac1\lambda B(x,\lambda \xi \cdot x)\right)=X+\left(\xi \times \p_\vartheta B+\frac1\lambda\curl_x(B)\right)(x,\lambda \xi \cdot x).
\end{align}
Hence, if we want to realize a change $w(x,\vartheta)$, we have to solve
\begin{equation}
    \label{eq:PotMot}
    w(x,\vartheta)=(\xi \times \p_\vartheta B)(x,\vartheta).
\end{equation}
Suppose that
\begin{equation}
    \label{eq:PotentialRequirementw}
    \overline{w}=0,\quad w\cdot \xi=0,
\end{equation}
where an overline denotes the $\vartheta$-average. If we denote by $\p_\vartheta^{-1}$ the unique average $0$ primitive due to \eqref{eq:PotentialRequirementw} we can define
\begin{equation}
    \label{eq:PotentialDefi}
    \mc P_\xi w=-\xi \times \p^{-1}_\vartheta w,
\end{equation}
and an elementary computation shows that $B=\mc P_\xi w$ solves \eqref{eq:PotMot}.

In view of the above discussion, we want to find $Y(x,\vartheta)$ such that the increment $w=Y-X$ satisfies \eqref{eq:PotentialRequirementw}. As our goal is to build a unit-length field, for $\rho \in C_c^\infty(\Omega)$ we want
\begin{equation*}
    |Y|^2=|X|^2+\rho^2,
\end{equation*}
so $Y$ increases the length of $X$ by a prescribed amount. Suppose
\[
    b:=|X-(X\cdot\xi)\xi|>0 \quad\text{near }\supp\rho,
\]
and denote
\[
    e_1=\frac{X-(X\cdot\xi)\xi}{b}.
\]
Under this assumption we can define
\begin{equation}
    Y=(X\cdot \xi)\xi+\sqrt{b^2+\rho^2} \bigl(\cos(a\sin\vartheta)e_1 +\sin(a\sin\vartheta)(\xi\times e_1)\bigr),
\end{equation}
where $a\in C^\infty(\Omega)$ is chosen such that
\[
\overline{Y-X}=0,\quad (Y-X)\cdot \xi=0,\quad |Y|^2=|X|^2+\rho^2.
\]
Such a choice of $a$ is possible if $\rho/b$ is sufficiently small. See Lemma \ref{lem:LocalPert}, in particular Step 1, for the precise construction. We perturb through
\[
X_\lambda=X+\curl\left(\frac{1}{\lambda}\mc P_\xi(Y-X)\right),
\]
and obtain by the discussion above
\begin{equation}
    \label{eq:ErrorLambda1}
    |X_\lambda|^2=|X|^2+\rho^2+O(\lambda^{-1}).
\end{equation}
This construction is well-known. In fact, the increment $Y-X$ is a scaled version of the corrugations used in Lemma 4.4 in \cite{ContiDeLellisSzekelyhidi2012}, which are based on Kuiper's corrugations \cite{Kuiper1955}. See also \cite{DeLellisInauenSzekelyhidi2018} and \cite{CaoSzekelyhidi2022Global}.

\begin{remark}
    There is also a direct connection to the Bloch wall discussed earlier. Suppose for simplicity $\boldsymbol{m}=|\boldsymbol{m}|e_1$ is constant with $|\boldsymbol{m}|<1$. If we suppose moreover that $\rho=\sqrt{1-|\boldsymbol{m}|^2}$, and take $\xi=e_3$ one can compute that the perturbation is of the form
    \[
    \boldsymbol{m}_\lambda(z)= \cos(a\sin(\lambda z))e_1+ \sin(a\sin(\lambda z))e_2+ O(\lambda^{-1}).
    \]
    The difference compared to a Bloch wall is that the vector field does not perform a $180^\circ$ rotation but instead varies only slightly according to the angle $a$. Compare Figure \ref{fig:bloch-building-block} below with \cite[Fig.~3.53(a), p.~201]{HubertSchafer1998}.
    \par\medskip
    \noindent\begin{minipage}{\linewidth}
    \centering

    \begin{tikzpicture}[x=1.38cm,y=1.38cm,
  frame/.style={black,line width=.52pt,line join=round},
  section/.style={draw=black,fill=white,line width=.43pt},
  moment/.style={black,line width=.85pt,-{Stealth[length=2.05mm,width=1.25mm]}},
  axis/.style={black,line width=.43pt,-{Latex[length=1.7mm,width=1.12mm]}},
  every node/.style={inner sep=1.4pt,font=\small}]
  \def\amplitude{0.40}
  \def\tubelength{4.8}
  \def\depthx{-0.57}
  \def\depthy{-0.135}
  \def\top{1.17}
  \def\side{1.02}

  \foreach \s in {0,\tubelength}{
    \draw[frame]
      ({\s+\side*\depthx},{\top+\side*\depthy}) --
      ({\s-\side*\depthx},{\top-\side*\depthy}) --
      ({\s-\side*\depthx},{-\top-\side*\depthy}) --
      ({\s+\side*\depthx},{-\top+\side*\depthy}) -- cycle;
  }
  \foreach \v in {-\top,\top}{
    \foreach \u in {-\side,\side}{
      \draw[frame] ({\u*\depthx},{\v+\u*\depthy}) --
                   ({\tubelength+\u*\depthx},{\v+\u*\depthy});
    }
  }
  \foreach \i in {0,...,8}{
    \pgfmathsetmacro{\s}{\tubelength*\i/8}
    \pgfmathsetmacro{\angledeg}{deg(\amplitude*sin(45*\i))}
    \begin{scope}[shift={(\s,0)}]
      \begin{scope}[x={(\depthx*1.38cm,\depthy*1.38cm)},y={(0cm,1.38cm)}]
        \path[section] (0,0) circle[radius=1];
      \end{scope}
      \pgfmathsetmacro{\vx}{\depthx*sin(\angledeg)}
      \pgfmathsetmacro{\vy}{cos(\angledeg)+\depthy*sin(\angledeg)}
      \draw[moment] ({-\vx},{-\vy}) -- (\vx,\vy);
    \end{scope}
  }

  \draw[frame] ({\side*\depthx},{-\top+\side*\depthy}) --
               ({\tubelength+\side*\depthx},{-\top+\side*\depthy}) --
               ({\tubelength+\side*\depthx},{\top+\side*\depthy}) --
               ({\side*\depthx},{\top+\side*\depthy}) -- cycle;
  \draw[frame] ({\side*\depthx},{\top+\side*\depthy}) --
               ({-\side*\depthx},{\top-\side*\depthy}) --
               ({\tubelength-\side*\depthx},{\top-\side*\depthy}) --
               ({\tubelength+\side*\depthx},{\top+\side*\depthy});
  \draw[frame] ({\tubelength-\side*\depthx},{\top-\side*\depthy}) --
               ({\tubelength-\side*\depthx},{-\top-\side*\depthy}) --
               ({\tubelength+\side*\depthx},{-\top+\side*\depthy});

  \draw[axis] (-1.01,-.33) -- (-1.01,.43) node[above] {$e_1$};
  \draw[axis] (-.09,-1.57) -- ({\tubelength+.50},-1.57)
    node[right] {$\lambda z$};
  \foreach \position/\ticktext in {0/0,.25/{\pi/2},.5/\pi,.75/{3\pi/2},1/2\pi}{
    \pgfmathsetmacro{\s}{\tubelength*\position}
    \draw[line width=.4pt] (\s,-1.54) -- (\s,-1.61);
    \node[below] at (\s,-1.62) {$\ticktext$};
  }
\end{tikzpicture}
    \captionsetup{justification=centering}
    
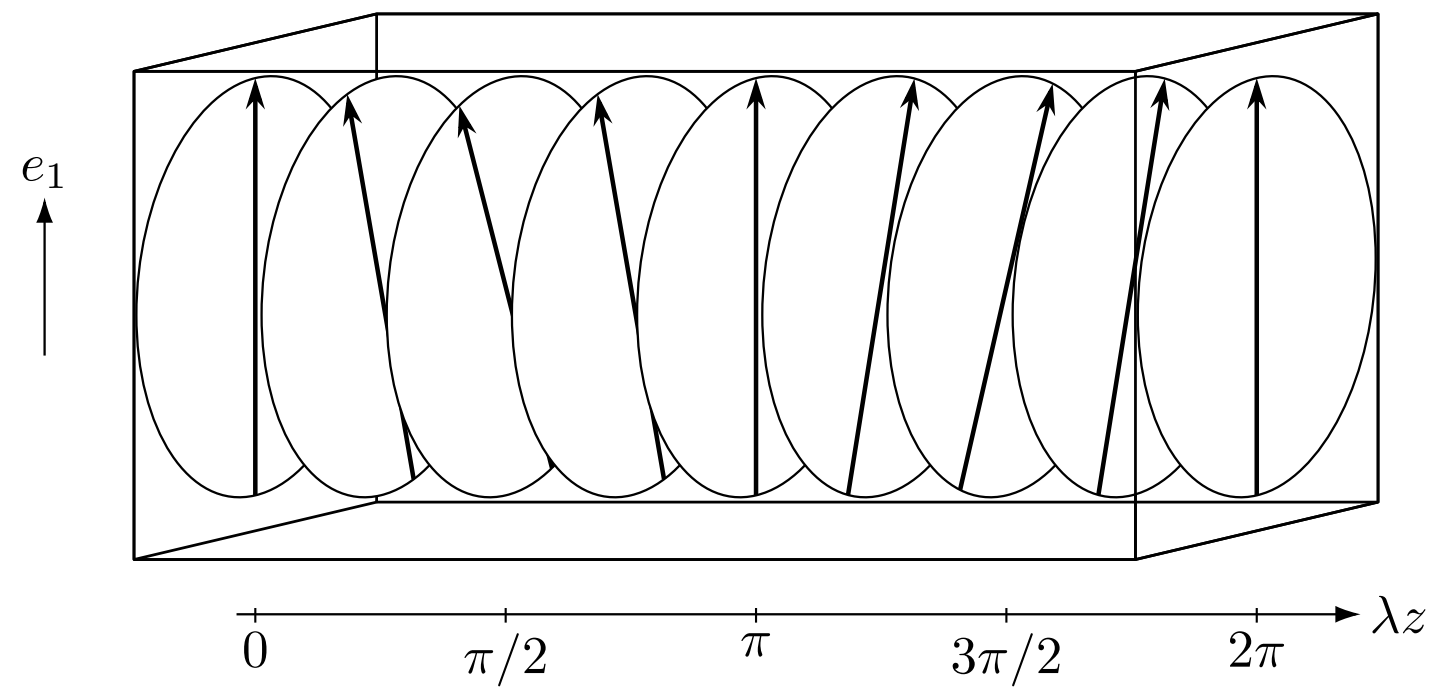
\captionof{figure}{Perturbation.}
    \label{fig:bloch-building-block}
\end{minipage}
\par\medskip
    \end{remark}
We use these perturbations at two points in the paper:
\begin{enumerate}
    \item As discussed before we use the corrugations in the proof of Lemma \ref{lem:RelativeExactExtension} to give a constructive proof of a \textit{relative h-principle} for nowhere vanishing divergence free vector fields. In particular, from this we can easily construct smooth nowhere vanishing divergence free vector fields which are almost unit-length.
    \item To remove the small length error and obtain a unit-length divergence free $X$ of higher $C^k$ regularity we want to use a Newton iteration. Linearizing the unit-length equation in \eqref{eq:DivFreeUnit}, for given smooth $f$ we have to find a potential $B$ such that
    \begin{equation}
        \label{eq:LinearizedUnitLength}
         X\cdot\curl B=f.
    \end{equation}
    A first idea is to construct an order zero inverse for \eqref{eq:LinearizedUnitLength}. For a smooth function $u$, we have from the product rule
    \begin{equation}
        \label{eq:ProdRuleHelc}
        X\cdot\curl(uX)
        =X\cdot(\nabla u\times X+u\curl X)
        =u(X\cdot\curl X).
    \end{equation}
    Hence we obtain a solution to \eqref{eq:LinearizedUnitLength} via
    \begin{equation}
        \label{eq:Inversetauzero1}
        B=\frac{f}{X\cdot\curl X}X
    \end{equation}
    under the nondegeneracy assumption of \textit{nowhere vanishing helicity} that is
    \[
    X\cdot \curl X \neq 0.
    \]
    This formula for the inverse is a special instance of \cite[p. 150]{Gromov1986}. Note that in order for the Newton iteration to succeed with this Ansatz, our perturbations must take care of two things simultaneously. First, they have to make the length error sufficiently small; second, they have to preserve nowhere vanishing helicity.
\end{enumerate}
This latter point is a central difficulty, especially when compared with Nash's proof. We recall that in the setting of isometric embeddings, Nash can separate these two roles of the perturbations due to his $(y,z)$--embedding Ansatz: the $z$--embedding takes care of the nondegeneracy condition for solving the linearized equation, while the metric induced by $y$ is made close to $g-g_z$ through adding perturbations.
    
However, in our case first of all we do not have the freedom of extra codimension to separate the two roles of the perturbations. Secondly, the perturbations may produce zeros in the helicity and hence the Ansatz \eqref{eq:Inversetauzero1} can not be used in our scheme. This forces us to invert \eqref{eq:LinearizedUnitLength} through a differential operator of higher order, which requires a more complicated nondegeneracy condition, see Lemma \ref{lem:InverseIntro} and Lemma \ref{lem:FiniteTypeAfterPer}. In particular, this means that the closeness condition in the Nash--Moser theorem has to be of higher order, see Proposition \ref{prop:NashMoser}. Therefore, when using the perturbations in Section \ref{sec:LocPert} the second time, we are forced to iteratively build versions of the perturbations, which correct a length error to higher order in $\lambda^{-1}$.

\textbf{Nash--Moser correction.} The higher order right inverse we construct in Lemma \ref{lem:InverseIntro} has a loss of derivative. Hence, as in the setting of isometric embeddings, a standard implicit function theorem cannot be applied. Therefore, following \cite{Moser1961,Moser1966}, in Section \ref{sec:NM} we prove a Nash--Moser correction theorem which requires the length error to be sufficiently small in a (weighted) $C^m$--norm, with a quantitative dependence on higher derivatives of the initial vector field and its nondegeneracy condition. Note that in order to preserve the prescribed boundary indices, we have to make sure that the Nash--Moser theorem produces a correction which vanishes near the boundary.

Finally, in Section \ref{sec:ProofMain} we complete the proof of the main theorem as follows. We use the vector field $X_0$, add the local perturbations in a careful way and then remove the final length error by an application of the Nash--Moser theorem from Section \ref{sec:NM}.
\subsubsection*{AI Disclosure}
This project is part of the second author's PhD program, which started in January 2025. By early 2026, the authors had a complete draft proving the main theorem for domains which are diffeomorphic to a ball and had outlined a strategy for arbitrary domains without the use of AI tools. As the usefulness of AI for mathematical research increased, AI tools, in particular OpenAI's ChatGPT 5.6 Sol and ChatGPT 6 Astra, were incorporated in the research and used to assist with computations and technical estimates. In particular, ChatGPT was used to choose the cutoff function in the proof of Proposition \ref{prop:Prep} and helped significantly with the computations in Section \ref{sec:LocPert}. For non-mathematical applications, AI was used to generate the figures and improve the wording of some paragraphs. The paper was written and checked by the authors, who take full responsibility for its correctness.

\subsubsection*{Acknowledgements}
The authors thank Matteo Giardi for fruitful discussions in the preparation of this work. LS would also like to thank Tristan Rivi\`ere for originally suggesting this problem to him and for many very useful discussions during his time at ETH Z\"urich.

\section{Preliminaries}
\label{sec:Prelim}
With $U$ and $\Sigma$ as in Theorem \ref{thm:main}, we write $N=\overline U\setminus\Sigma$ and $\partial N=\partial U\setminus\Sigma$. The symbols $|\cdot|$ and $\cdot$ denote the Euclidean norm and inner product. We use $D^k$ for derivatives of order $k$ and $\nabla f$ for the gradient of a scalar function. For fields defined only on the boundary, derivative estimates refer to their constant extensions along normal lines.\\
In view of the asymptotics described in Theorem \ref{thm:main} we introduce the following weight. Let $\sigma\in C^1(\mb R^3) \cap C^\infty(\mathbb R^3\setminus\Sigma)$ be nonnegative and positive on $\mathbb R^3\setminus\Sigma$, such that
\begin{equation}
    \label{eq:ChoiceSigma}
    \sigma(x)=\begin{cases}
    |x-p|^{1+\delta},\quad \text{near each }p\in\Sigma,\\
    1,\quad \text{outside a fixed compact set}.
    \end{cases}
\end{equation}
By this choice of $\sigma$, we have
\begin{equation}
    \label{eq:DerBoundSig}
    ||\nabla\sigma||_{L^\infty(\mathbb R^3)}<\infty,
    \qquad
    |D^k\sigma|\le C_k\sigma^{1-k}
    \quad\text{on }\mathbb R^3\setminus\Sigma,\quad k\ge1.
\end{equation}
Let $E\subset\mathbb R^3\setminus\Sigma$. For a vector field or function $u$ we use the norms
 \[
||u||_{k;E}=\max_{0\leq j \leq k} \sup_{E} \sigma^j|D^j u|.
\]
For a  vector potential $B$ we use
\[
||B||_{k,1;E}=\max_{0\leq j \leq k} \sup_{E} \sigma^{j-1}|D^j B|.
\]
Moreover, we write
\[
    ||u||_{k}:=||u||_{k;N},\qquad ||B||_{k,1}:=||B||_{k,1;N}.
\]
By contrast, $||\cdot||_{C^k(E)}$ denotes the usual norm of coordinate components and their partial derivatives up to order $k$ on $E\subset\mathbb R^3\setminus\Sigma$. The reason we use the weight $\sigma^{j-1}$ for vector potentials is that with this normalization we have
\[
    ||\curl B||_k \leq C ||B||_{k+1,1}.
\]
We also need a weighted distance which we define as
\[
    \dist_\sigma(x,y)=\inf_{\gamma}\int_0^1\frac{|\dot \gamma(t)|}{\sigma(\gamma(t))}dt.
\]
Here the infimum is taken over piecewise smooth curves $\gamma:[0,1]\to N$ with $\gamma(0)=x$ and $\gamma(1)=y$.\\
The next lemma gives a Whitney-type covering that we will use for the local perturbations.
\begin{lemma}
\label{lem:BoundCharts}
There exists a countable, locally finite family of balls
\[
    U_i(t)=B_{t\ell_i}(x_i),
    \qquad U_i(2)\subset\mathbb R^3\setminus\Sigma,
\]
such that
\[
    \mathbb R^3\setminus\Sigma=\bigcup_i U_i(1),
    \qquad
    \sigma(x)\simeq\ell_i
    \quad\text{for }x\in U_i(2).
\]
Moreover, there is a partition
\begin{equation}\label{eq:ChartPartition}
    I=\bigsqcup_{q=1}^Q I_q,
    \qquad
    i,j\in I_q,\quad i\ne j
    \quad\Longrightarrow\quad
    U_i(2)\cap U_j(2)=\emptyset.
\end{equation}
If $E\subset U$ satisfies $d=\dist_\sigma(E,\partial N)>0$, an analogous family can be chosen with
\[
    E\subset\bigcup_i U_i(1),
    \qquad U_i(2)\subset U.
\]
For this family, the constants depend on $d$.
\end{lemma}
\begin{figure}[htbp]
\centering
\resizebox{0.46\textwidth}{!}{%
\begin{tikzpicture}[
    font=\small,
    line cap=round,
    line join=round
]
    \path[use as bounding box]
        (-2.05,-1.1) rectangle (3.65,1.1);

    \foreach \cx/\rad in {
        -1.0/0.88,
         0.5/0.43,
         1.5/0.19,
         2.1/0.075,
         2.5/0.028,
         2.75/0.012
    } {
        \path[
            fill=profileblue!12,
            draw=profileblue,
            line width=0.65pt
        ]
            (\cx,0) circle[radius=\rad];
    }

    \draw[fill=white,draw=black!75,line width=0.95pt]
        (3,0) circle[radius=2.8pt];

    \node at (-1.0,0) {$U_i(2)$};
    \node[right,inner sep=6pt] at (3,0) {$p$};
\end{tikzpicture}%
}
\caption{Balls shrinking towards $p$.}
\label{fig:BoundCharts}
\end{figure}
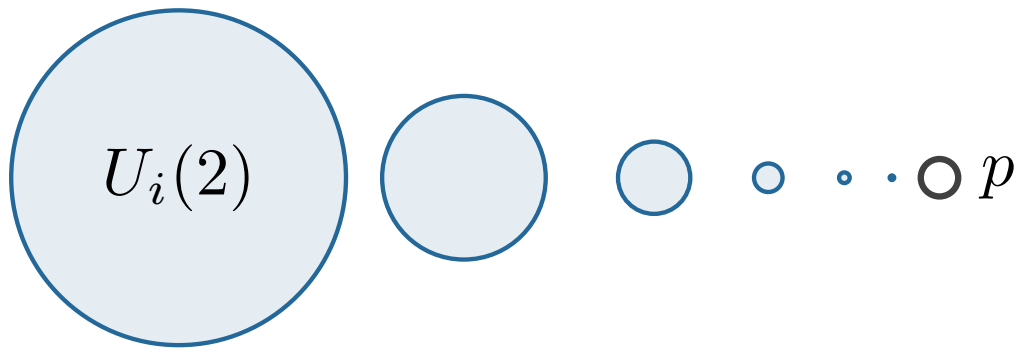

\begin{proof}
Recall that $L:=||\nabla \sigma||_{L^\infty(\mb R^3)}<\infty$. Choose $a>0$ sufficiently small that $aL\le1/20$. Hence
\[
    \sigma(y)\simeq\sigma(x) \quad\text{whenever }|y-x|\le2a\sigma(x).
\]
In the language of H\"ormander this means that the norms $|v|/(a\sigma(x))$ are slowly varying. Thus, H\"ormander's covering theorem \cite[Theorem~1.4.10]{Hormander2003} gives countably many $x_i\in\mathbb R^3\setminus\Sigma$ such that, with $\ell_i=a\sigma(x_i)$, the balls $U_i(1)$ cover
$\mathbb R^3\setminus\Sigma$ with uniformly bounded multiplicity. That is, there exists $M\in\mathbb N$ such that
\begin{equation}
    \label{eq:MultiplicityBound}
    \#\{j\in I:x\in U_j(1)\}\le M
    \qquad\text{for every }x\in\mathbb R^3\setminus\Sigma.
\end{equation}
The same Lipschitz estimate gives $U_i(2)\subset\mathbb R^3\setminus\Sigma$ and
$\sigma\simeq\ell_i$ on $U_i(2)$.\\
Fix an index $i$. If $U_i(2)\cap U_j(2)\ne\emptyset$, then
\[
    |\ell_i-\ell_j|\le aL|x_i-x_j|<2aL(\ell_i+\ell_j),
\]
so $\ell_i\simeq\ell_j$. Thus, for all indices $j$ as above the sets $U_j(1)$ have volume comparable to $\ell_i^3$ and are contained in $B_{C\ell_i}(x_i)$. By \eqref{eq:MultiplicityBound},
\begin{align*}
    \ell_i^3\,\#\{j\in I:U_i(2)\cap U_j(2)\ne\emptyset\}
    \simeq
    \sum_{\substack{j\in I\\U_i(2)\cap U_j(2)\ne\emptyset}}
    |U_j(1)|
    \le M\,|B_{C\ell_i}(x_i)|
    \simeq M\ell_i^3.
\end{align*}
Thus the number of intersecting sets is bounded independently of $i$. Pick an integer $Q$ with
\[
    \#\{j\in I:U_i(2)\cap U_j(2)\ne\emptyset\}\le Q-1 \qquad\text{for every }i\in I.
\]
Enumerate the indices and assign successively each $i$ to one of $I_1,\ldots,I_Q$, while avoiding any set already containing an index $j$ with $U_i(2)\cap U_j(2)\ne\emptyset$. Note that for each $i$ at most $Q-1$ sets have to be avoided, so this is always possible. In particular
\[
    U_i(2)\cap U_j(2)=\emptyset \qquad\text{for }i\ne j\in I_q.
\]
This gives \eqref{eq:ChartPartition}. The family is locally finite: every point belongs to
some $U_i(1)$, which intersects only finitely many of the sets $U_j(2)$.\\
For the interior cover, repeat the construction with the choice $a<d/4$, and keep only the indices $i$ for which $U_i(1)\cap E\ne\emptyset$. If $e\in E\cap U_i(1)$ and $z\in U_i(2)$, the segment $[e,z]$ lies in $U_i(2)$ and satisfies
\begin{equation}
    \label{eq:weightLengthlessd}
    \int_{[e,z]}\frac{ds}{\sigma}\le \frac{|e-z|}{(1-2aL)\sigma(x_i)} \le\frac{3\ell_i}{(1-2aL)\sigma(x_i)} <4a<d.
\end{equation}
If $z\notin U$, a subsegment of $[e,z]$ would lie in $N$ and join $e$ to a point of $\partial U\setminus\Sigma=\partial N$, so \eqref{eq:weightLengthlessd} would contradict the definition of $d$.
Thus $U_i(2)\subset U$. The rest of the properties are directly inherited from the ambient cover, with constants now of course allowed to depend on $d$.
\end{proof}
For either cover, define
\[
    \Phi_i(y)=x_i+\ell_i y,\qquad y\in B_2(0),
\]
and set
\begin{equation}
    \label{eq:ScaledMetricWeight}
    \widehat\sigma_i=\frac{\sigma\circ\Phi_i}{\ell_i}.
\end{equation}
The comparison $\sigma\simeq\ell_i$ on $U_i(2)$ and \eqref{eq:DerBoundSig} easily imply
\begin{equation}\label{eq:ScaledGeometryBounds}
    ||\widehat\sigma_i||_{C^k(B_2(0))}+||\widehat\sigma_i^{-1}||_{C^k(B_2(0))} \le C_k,\qquad k\ge0.
\end{equation}
For a function or vector field $u$ and a vector potential $B$, use the normalizations
\begin{equation}\label{eq:ChartNormalizations}
    \widehat u_i=u\circ\Phi_i,\qquad
    \widehat B_i=\ell_i^{-1}B\circ\Phi_i.
\end{equation}
Since
\[
    D_y^j\widehat u_i
    =\ell_i^j(D^j u)\circ\Phi_i,\qquad
    D_y^j\widehat B_i
    =\ell_i^{j-1}(D^jB)\circ\Phi_i,
\]
we have
\begin{equation}
    \label{eq:Normequivalence}
    ||\widehat u_i||_{C^k(B_2(0))}\simeq ||u||_{k;U_i(2)},\qquad ||\widehat B_i||_{C^k(B_2(0))} \simeq||B||_{k,1;U_i(2)}.
\end{equation}
All comparison constants are independent of $i$. For the interior cover, they may also depend on $d$. The main reason for introducing the cover and rescalings, is that it allows us to perform the local construction in Section \ref{sec:LocPert} in a fixed ball and the usual $C^k$ norms. Through \eqref{eq:Normequivalence} we can then easily transfer statements to the weighted norms on the covering.\\
We give next the smoothing operators that will be used in the proof of Proposition \ref{prop:NashMoser}. This construction is standard (see \cite[Section 2.3.4]{Gromov1986}), but since we have to make sure the smoothing operators respect the support of the vector fields in a suitable way, we give an outline of the proof.
\begin{lemma}
    \label{lem:Smoothing}
    Fix $\bar k\in\mathbb N$ and $d>0$. There is a constant $C\geq 1$, depending only on $\bar k$, such that for $R\geq 2C(1+d^{-1})$, there are linear smoothing operators $S_R$ on vector potentials $B$ supported in $\{x:\dist_\sigma(x,\partial N)\ge d\}$ which satisfy
     \begin{align}
        \label{eq:Smoothing1}
    ||S_RB||_{r+s,1}
        &\le C_{r,s}R^s||B||_{r,1},
        &&r,s\ge0,\quad r+s\le\bar k,\\
        \label{eq:Smoothing2}
    ||(I-S_R)B||_{r,1}
        &\le C_{r,s}R^{-s}||B||_{r+s,1},
        &&r,s\ge0,\quad r+s\le\bar k.
    \end{align}
    Moreover,
     \[
    \supp S_RB \subset\{x:\dist_\sigma(x,\supp B)\le C/R\} \subset\{x:\dist_\sigma(x,\partial N)\ge d/2\}.
    \]
\end{lemma}
\begin{proof}
By the support assumption on $B$, we can extend it smoothly by zero outside $U$ to $\mb R^3 \setminus \Sigma$. Use the standard convolution construction of \cite[Section 2.3.4]{Gromov1986} in the ball cover of $\mathbb R^3\setminus\Sigma$ from Lemma~\ref{lem:BoundCharts}. Choose a partition of unity $(\chi_i)$ with $\supp\chi_i\subset U_i(3/2)$ and cutoffs $\widetilde\chi_i$ supported in $U_i(2)$ and equal to $1$ near $\supp\chi_i$, such that
\[
    ||\chi_i\circ\Phi_i||_{C^k(B_2(0))}+||\widetilde\chi_i\circ\Phi_i||_{C^k(B_2(0))}\le C_k \qquad\text{for every }k\ge0,
\]
where $C_k$ is independent of $i$. Set $\widehat B_i=\ell_i^{-1}B\circ\Phi_i$ and extend $(\widetilde\chi_i\circ\Phi_i)\widehat B_i$ smoothly by zero to $\mathbb R^3$.\\
Choose $\psi\in C_c^\infty(B_1(0))$ with integral $1$ and vanishing moments of orders $1,\ldots,\bar k-1$. Write $\psi_R(x)=R^3\psi(Rx)$ and define on $\mathbb R^3\setminus\Sigma$
\begin{equation}
    \label{eq:SmoothingOperator}
    S_RB=\sum_i\chi_i\ell_i
    \left(\psi_R*\bigl[(\widetilde\chi_i\circ\Phi_i)\widehat B_i\bigr]\right)\circ\Phi_i^{-1}.
\end{equation}
The standard Euclidean convolution estimates, the norm equivalence \eqref{eq:Normequivalence} and Lemma \ref{lem:BoundCharts} give \eqref{eq:Smoothing1}--\eqref{eq:Smoothing2}.\\
Finally, note that convolution enlarges each coordinate support by at most $R^{-1}$. Let $x$ lie in the support of \eqref{eq:SmoothingOperator} on $\mathbb R^3\setminus\Sigma$. There exist an index $i$ and $y\in\supp B$ such that
\[
    x,y\in U_i(2),\qquad |x-y|\le\ell_i/R.
\]
Since $\sigma\simeq\ell_i$ on $U_i(2)$, we derive
\[
    \int_{[y,x]}\frac{ds}{\sigma}
    \le C\frac{|x-y|}{\ell_i}
    \le C/R.
\]
For $R\geq2C(1+d^{-1})$, we have $C/R\leq d/2<d$. If $[y,x]\not\subset U$, a subsegment would be contained in $N$ and join $y$ to a point of $\partial N$, with weighted length at most $C/R<d$. This contradicts $\dist_\sigma(y,\partial N)\geq d$. Thus $[y,x]\subset U$ and
\[
    \dist_\sigma(x,\supp B)\leq\dist_\sigma(x,y)\leq C/R.
\]
It easily follows
\[
    \supp S_RB\subset\{x:\dist_\sigma(x,\supp B)\leq C/R\}\subset\{x:\dist_\sigma(x,\partial N)\geq d/2\},
\]
which completes the proof.
\end{proof}
Next we introduce the nondegeneracy condition which will allow us to construct an inverse for the linearized equation. See Lemma \ref{lem:InverseIntro} below. As discussed in the introduction, in Nash's isometric embedding theorem the condition which plays the same role is called \textit{freeness}.\\
The reason we have to work with an order-dependent condition is that the local perturbations we use may destroy the nondegeneracy condition for some fixed $m$, while nondegeneracy is retained at higher orders, see Lemma \ref{lem:FiniteTypeAfterPer}.\\
To motivate the nondegeneracy condition, fix a smooth nowhere vanishing
vector field $X$. Linearizing the unit-length equation in
Theorem \ref{thm:main}, for given smooth $f$ we want to solve
\begin{equation}
    \label{eq:LinearizedEquation}
    X\cdot\curl B=f.
\end{equation}
Recall from the introduction, if $X\cdot\curl X$ never vanishes, we can solve
\eqref{eq:LinearizedEquation} through
\begin{equation}
    \label{eq:Inversetauzero}
    B=\frac{f}{X\cdot\curl X}X.
\end{equation}
A similar idea works if a derivative of $X\cdot\curl X$ in a
direction perpendicular to $X$ is nowhere zero. To see this,
let $H$ be a smooth vector field with $H\cdot X=0$, and write
for smooth $v$, $Hv=H\cdot\nabla v$. Pick $Y$ such that $Y\times X=H$. Then
\begin{equation}
    \label{eq:ProductRuleuY}
    X\cdot\curl(vY)
    =X\cdot(\nabla v\times Y+v\curl Y)
    =Hv+vX\cdot\curl Y.
\end{equation}
Applying \eqref{eq:ProductRuleuY} with $v=u(X\cdot\curl X)$ gives
\[
\begin{aligned}
    &X\cdot\curl\bigl(u(X\cdot\curl X)Y\bigr)\\
    &\qquad=H\bigl(u(X\cdot\curl X)\bigr)
      +u(X\cdot\curl X)(X\cdot\curl Y)\\
    &\qquad=uH(X\cdot\curl X)
      +(X\cdot\curl X)\bigl(Hu+uX\cdot\curl Y\bigr)\\
    &\qquad=uH(X\cdot\curl X)
      +X\cdot\curl\bigl((Hu+uX\cdot\curl Y)X\bigr).
\end{aligned}
\]
The final equality follows from \eqref{eq:ProdRuleHelc}.
Rearranging gives the choice
\begin{equation}
    \label{eq:ChoiceB}
    B=u(X\cdot\curl X)Y-\bigl(Hu+uX\cdot\curl Y\bigr)X,
\end{equation}
which satisfies
\[
    X\cdot\curl B=uH(X\cdot\curl X).
\]
Hence, wherever $H(X\cdot\curl X)\neq0$ choosing
\begin{equation}
    \label{eq:Invtau1}
    u=\frac{f}{H(X\cdot\curl X)}
\end{equation}
in \eqref{eq:ChoiceB} means that $B$ solves \eqref{eq:LinearizedEquation}. The following definition gives a quantitative and higher-order version of this condition.\\
We write
\[
    \Pi_{X^\perp}Y=Y-\frac{Y\cdot X}{|X|^2}X \qquad(X\ne0)
\]
for orthogonal projection onto $X^\perp$.
\begin{definition}
    \label{def:Nondegen}
    Let $\mathbf e_1,\mathbf e_2,\mathbf e_3$ be the standard basis of $\mb R^3$. For a nowhere vanishing  $X\in C^\infty(N;\mathbb R^3)$, set
    \[
        q_X=\sigma X\cdot\curl X,\qquad
        H_\alpha^X=\sigma \Pi_{X^\perp}\mathbf e_\alpha,
        \qquad \alpha=1,2,3.
    \]
    We also use $H_\alpha^X$ for the directional derivative $H_\alpha^X\cdot\nabla$. For $I=(\alpha_1,\ldots,\alpha_r)$, write $H_I^X=H_{\alpha_1}^X\cdots H_{\alpha_r}^X$, $|I|=r$, and $H_\emptyset^X=\mathrm{Id}$. Define
    \[
        \tau_m(X)=\left(\sum_{|I|\leq m}|H_I^Xq_X|^2\right)^{1/2}.
    \]
     A vector field $X$ is uniformly nondegenerate of order $m$ in $E\subset N$ if  $\tau_m(X)\geq\kappa>0$ on $E$.
\end{definition}
Using the nondegeneracy condition above we can prove the following lemma. We refer to Section 5 for its proof. Write for $X \in C^\infty(N;\mathbb R^3)$
\[
P_X B =2X\cdot\curl B,\quad B \in C^\infty(N;\mathbb R^3).
\]
\begin{lemma}
\label{lem:InverseIntro}
Suppose that $X \in C^\infty(N;\mathbb R^3)$ is nowhere vanishing and $\tau_m(X)>0$ in $N$. Then there exists a differential right inverse $L_X$ of order $m$ for $P_X$. In particular,
\[
    \supp L_X f\subset\supp f.
\]
The coefficients of $L_X$ depend smoothly on the derivatives of $X$ through order $2m+1$. For every $k\ge0$ and every $f\in C^\infty(N)$,
\begin{equation*}
    ||L_X f||_{k,1}\le C_k\left(1+||X||_{2m+1+k}+|||X|^{-1}||_0+||\tau_m(X)^{-1}||_0\right)^{P_k}||f||_{m+k}.
\end{equation*}
The constants $C_k,P_k>0$ depend only on $m$ and $k$.
\end{lemma}
\begin{remark}
    The solution \eqref{eq:Inversetauzero} corresponds to $\tau_0(X)>0$, while \eqref{eq:Invtau1} corresponds to $\tau_1(X)>0$. 
\end{remark}

\section{Global Short Field}
\label{sec:Prep}
The central result of this section is the construction of a vector field which will serve as the preparatory field that we modify later using local perturbations in order to satisfy the conditions of Proposition \ref{prop:NashMoser}. The idea is that with this preparation we take care of the topological properties of the problem, which in our case are the singularities on the boundary. At the same time by making the field strictly short away from the boundary we have enough room to add the local perturbations of Section \ref{sec:LocPert}.

\begin{proposition}
    \label{prop:Prep}
    There exist a smooth potential $B_0$ with $X_0=\curl B_0\in C^\infty(N;\mathbb R^3)$ and constants $c_0,\kappa_0,d_*>0$ such that with  $\Delta:=1-|X_0|^2$ the following holds.
    \begin{enumerate}
        \item \label{prop:it1} We have  $c_0\le |X_0|\le 1$ on $N$ and $\dist_\sigma(\supp\Delta,\partial N)\geq d_*$.
        \item \label{prop:it2} The  field satisfies $X_0\cdot\nu=0$ on $\partial N$ and $X_0$ has the indices as in Theorem \ref{thm:main}.
        \item \label{prop:it3} We have $\Delta+\tau_0(X_0)\geq\kappa_0$ on $N$.
        \item \label{prop:it4} We have $\sqrt{\Delta}\in C^\infty(N)$ and, for every $k\geq0$ and $0<\eps<1$,
         \[
        \sigma^k|D^k\sqrt{\Delta}|\leq C_{k,\eps}(\sqrt{\Delta})^{1-\eps}
        \quad\text{on }N.
        \]
        \item \label{prop:it5} Finally, $||X_0||_k<\infty$ for all $k \ge 0$.
    \end{enumerate}
\end{proposition}
First, we construct in Lemma \ref{lem:SmoothFrameInd} a frame which yields in Corollary \ref{coro:TangentFieldWithCorrectIndex} suitable boundary data for $X_0$. In Lemma \ref{lem:FormOnEnd}, the frame will also be used to write down an explicit potential in a neighborhood of $p \in \Sigma$ with nowhere vanishing curl.
\begin{lemma}
\label{lem:SmoothFrameInd}
Let $S\subset\partial U$ be a boundary component of genus $\gamma$, and $p\in S$. Choose positively oriented boundary-normal coordinates $x=(s,t)$ centered at $p$. Write
\[
r=|x|, \qquad U_\varepsilon(p)=\{x:0<r<\varepsilon,\ t\geq 0\}.
\]
For $\varepsilon>0$ small enough, there exist a smooth $\varphi:U_\varepsilon(p)\to\mathbb R$ and a positively oriented orthonormal frame $(n,e_1,e_2)$ on $U_\varepsilon(p)$, with the following properties.
\begin{enumerate}
\item \label{lemF:it1}
We have
\[
cr^{-1-\delta}\leq |\nabla\varphi|\leq Cr^{-1-\delta},\qquad n=\frac{\nabla\varphi}{|\nabla\varphi|}.
\]
\item \label{lemF:it2}
On  $U_\varepsilon(p)\cap\partial U$,
\[
\varphi=0, \qquad n=-\nu,
\]
Consequently, $e_1$ and $e_2$ are tangent to  $\partial U$ on $\p U\cap U_\eps(p)$.
\item \label{lemF:it3}
The field $e_2$ has
\[
\ind_p(e_2)=\chi(S)=2-2\gamma.
\]
\item \label{lemF:it4} For every $k\geq 0$, we have
\[
|D^k n|+|D^k e_1|+|D^k e_2| \leq C_k r^{-k} \qquad\text{on }U_\varepsilon(p).
\]
\end{enumerate}
\end{lemma}

\begin{proof}
    We first give the proof for $\gamma=0$ in the Euclidean half-space $\{t\geq0\}$, with $p=0$. Write $x=(s,t)=r\omega$, where $\omega\in\mb S^2_+:=\{\omega\in\mb S^2:\omega_3\geq0\}$. Define
    \begin{equation}
        \label{eq:AnsatzVarphi}
        \varphi(x)=\frac{1}{r^\delta}\omega_3=\frac{t}{r^{1+\delta}},\qquad
        \nabla\varphi(x)=r^{-1-\delta}(\mathbf e_3-(1+\delta)\omega_3\omega).
    \end{equation}
    Since $|\mathbf e_3-(1+\delta)\omega_3\omega|^2=1-(1-\delta^2)\omega_3^2$, we have
    \[
        \delta r^{-1-\delta}\le|\nabla\varphi|\le r^{-1-\delta},\qquad
        n=\frac{\nabla\varphi}{|\nabla\varphi|}
        =\frac{\mathbf e_3-(1+\delta)\omega_3\omega}{\sqrt{1-(1-\delta^2)\omega_3^2}},
    \]
    which proves \ref{lemF:it1}. Set
    \[
        e_1=\frac{\Pi_{n^\perp}(2\omega_1 \omega-\mathbf{e_1})}{|\Pi_{n^\perp}(2\omega_1 \omega-\mathbf{e_1})|},\qquad
        e_2=n\times e_1.
    \]
    It is simple to see that this is well-defined and $(n,e_1,e_2)$ is a positively oriented orthonormal frame. On $t=0$, we have $\varphi=0$ and $n=\mathbf e_3=-\nu$. Consequently, $e_1,e_2$ are tangent to the boundary, proving \ref{lemF:it2}.\\
    All three fields are smooth functions of $\omega=x/r$. Differentiating $x/r$ produces a factor $r^{-1}$, hence
    \[
        |D^k n|+|D^k e_1|+|D^k e_2|
        \leq C_kr^{-k},\quad \forall k\geq0,
    \]
    which proves \ref{lemF:it4}.\\
    It remains to prove \ref{lemF:it3}. On a circle around the origin in $\{t=0\}$, write $\omega=(\cos\alpha,\sin\alpha,0)$. Then
    \[
        e_2=(-\sin(2\alpha),\cos(2\alpha),0).
    \]
    Thus $e_2$ winds twice around the origin, so $\ind_0(e_2)=2$, proving \ref{lemF:it3}.\\
    We briefly sketch the changes for the general case. For genus $\gamma$, choose a Morse function $f\in C^\infty(\mb S^2_+)$, with $f=\omega_3$ near the equator, $\nabla_{\mb S^2}f\neq0$ on $\{f=0\}$, and $\{f<0\}$ consisting of $\gamma$ disjoint discs whose closures lie in $\operatorname{int}(\mb S^2_+)$. Following the Ansatz \eqref{eq:AnsatzVarphi} set
    \[
        \varphi(r\omega)=r^{-\delta}f(\omega),\qquad
        \nabla\varphi(r\omega)=r^{-1-\delta}
        (\nabla_{\mb S^2}f-\delta f\omega).
    \]
    The defining properties of $f$ give $|\nabla\varphi|\simeq r^{-1-\delta}$. Since $\mb S^2_+$ is contractible, $n=\nabla \varphi/|\nabla \varphi|$ can be completed to a smooth positively oriented orthonormal frame, with the boundary properties and derivative bounds as above. To compute the index, write $w=n\times\omega=w_1e_1+w_2e_2$. Its zeros are the critical points of $f$, where $n=-\mathrm{sgn}(f)\omega$. The Poincaré-Hopf Theorem gives total gradient index $1$ on $\mb S^2_+$ and $\gamma$ on $\{f<0\}$. Comparing orientations therefore gives
    \[
        \sum_{q:\nabla_{\mb S^2}f(q)=0}\ind_q(w_1,w_2)
        =\gamma-(1-\gamma)=2\gamma-1.
    \]
    Since $w$ winds once along the equator, this implies $\ind_0(e_2)=1-(2\gamma-1)=2-2\gamma$.\\
    The extension to a curved boundary is straightforward.
\end{proof}
The following corollary is derived by a standard cutoff argument from Lemma \ref{lem:SmoothFrameInd}. We omit the proof.
\begin{corollary}
    \label{coro:TangentFieldWithCorrectIndex}
    For each $p \in \Sigma$ construct the fields $e_2$ as in Lemma \ref{lem:SmoothFrameInd}. These local fields extend to a smooth unit-tangent field $T$ on $\p U\setminus \Sigma$. In particular near each $p \in \Sigma$, we have
    \begin{equation}
        \label{eq:BoundUnitLength}
        \ind_p(T)=\chi(S_p),\quad |D^k T(x)| \lesssim_k \frac{1}{|x-p|^{k}}
    \end{equation}
    where $S_p$ is the boundary component containing $p$.
\end{corollary}

In the following we assume that $\eps$ is picked small enough such that the statement asserted in Lemma \ref{lem:SmoothFrameInd} holds for every $p \in \Sigma$.

\begin{lemma}
    \label{lem:FormOnEnd}
    Let $p \in \Sigma$. In the notation of Lemma \ref{lem:SmoothFrameInd}, define
     \[
    B_e=\frac{\sin(\varphi)e_1+(\cos(\varphi)-1)e_2}{|\nabla\varphi|},\quad \text{on }U_\eps(p).
    \]
    For $X_e=\curl B_e$, we have $X_e=T=e_2$ on $U_\eps(p)\cap \p U$ and
    \begin{equation}
        \label{eq:NormOne}
        \bigl||X_e|-1\bigr|\leq Cr^\delta.
    \end{equation}
    Moreover, for every $k \geq 0$ we have
     \begin{equation}
        \label{eq:OmeDerBounds}
        ||X_e||_{k;U_\eps(p)}\leq C_{k}.
    \end{equation}
\end{lemma}

\begin{proof}
From Lemma \ref{lem:SmoothFrameInd}, we have
 \begin{equation}
    \label{eq:FrameBound}
    |D^j \varphi|\leq C_j r^{-j-\delta},\quad \left|D^j \left(\frac{e_i}{|\nabla\varphi|}\right)\right|\leq C_j r^{1+\delta-j}.
\end{equation}
Using the definition of  $B_e$, we compute
 \begin{align}
    \label{eq:FormBe}
    X_e&=n \times \left(\cos(\varphi)e_1-\sin(\varphi)e_2\right) \nonumber \\
    &+\sin(\varphi)\curl\left(\frac{e_1}{|\nabla\varphi|}\right)+(\cos(\varphi)-1)\curl\left(\frac{e_2}{|\nabla\varphi|}\right).
\end{align}
The first term in \eqref{eq:FormBe} has length $1$, hence by using \eqref{eq:FrameBound} it follows easily that
 \[
\bigl|\,|X_e|-1\,\bigr|\leq Cr^\delta,
\]
which proves \eqref{eq:NormOne}.\\
On the boundary  $U_\eps(p)\cap \p U$, we have $\varphi=0$ and  $n=-\nu$, so $X_e=n\times e_1=e_2=T$.\\
For the derivative bounds, suppose $r\leq 1$. The chain rule together with \eqref{eq:FrameBound} gives
 \[
|D^k \sin(\varphi)|+|D^k\cos(\varphi)|\leq C_k r^{-(1+\delta)k}.
\]
From this formula we have easily that
 \[
|D^k X_e|\leq C_k|D^{k+1}B_e|\leq C_k r^{-(1+\delta)k},
\]
so multiplying by $\sigma^k\simeq r^{(1+\delta)k}$ proves \eqref{eq:OmeDerBounds} which completes the proof of the lemma.
\end{proof}
We next construct a unit-length, divergence-free field near $\partial N$. Since our boundary data $T$ is singular we cannot expect to construct such a field on a neighborhood of $\p U$. Hence we introduce a boundary weight and collars with thickness comparable to $\sigma$.\\
Choose a positive $\sigma_\partial\in C^\infty(\partial N)$ with $\sigma_\partial(x)=|x-p|^{1+\delta}$ near each $p\in\Sigma$, such that $0<\sigma_\partial\leq1$ and $\sigma_\partial|D\nu|\leq1/4$ on $\partial N$. Near $\partial N$, we use boundary-normal coordinates $(s,t)$, where $s=(s^1,s^2)$ are coordinates on $\partial N$ and $t\geq0$ is the inward normal distance. Extend $\sigma_\partial$ and $\nu$ constantly along $t$ in a normal neighborhood of $\partial N$. Set
\[
\theta(s,t)=\frac{t}{\sigma_\partial(s)},\qquad
\mathcal C_{\theta_0}=\{(s,t):0\leq\theta(s,t)\leq\theta_0\},
\]
for $\theta_0>0$ small enough such that the $(s,t)$ coordinates are defined. By a \textit{collar} we mean a set $\mathcal C_{\theta_0}$ of this form. The bound $\sigma_\partial|D\nu|\leq1/4$ is used to get the lower bound on the helicity below. Since we want to solve a transport equation in the $t$ direction it is convenient for the weight to be independent of $t$.
\begin{lemma}
    \label{lem:UnitCollar}
    There exist constants $\theta_0,\kappa_c>0$ with the following properties. On the collar $\mathcal C_{\theta_0} \subset N$ there exists  $B_c \in C^\infty(\mc C_{\theta_0};\mathbb R^3)$ such that $X_c=\curl B_c$ satisfies
    \[
    |X_c|=1,\quad \text{on }\mathcal C_{\theta_0},\qquad X_c=T,\quad \text{on }\p N.
    \]
    Moreover, we have the estimates
     \begin{equation}
        \label{eq:BoundsXc}
        \begin{aligned}
            q_{X_c}\geq\kappa_c\quad\text{on }\mathcal C_{\theta_0}, \\
            ||X_c||_{k;\mathcal C_{\theta_0}}+||B_c||_{k,1;\mathcal C_{\theta_0}}\leq C_k,\quad \forall k \geq 0.
        \end{aligned}
    \end{equation}
\end{lemma}

\begin{proof}
     Write $F(s,t)$ for the boundary-normal coordinate map, where $t$ is the inward normal distance. Hence $\partial_t=-\nu\cdot\nabla$.\\
    \textbf{Step 1: Transport Equation} We  construct a divergence-free vector field $X_c$ in a neighborhood of $\p N$. For this, extend $T$, so that $T(F(s,t))=T(F(s,0))$. This preserves its length and tangency to the slices $\{t=\tilde t\}$.  We denote
    \[
    V=\cos\theta\,T+\sin\theta\,(\nu\times T)
    \]
    and make the Ansatz
     \[
    X_c=\sqrt{1-b^2}V-b\nu,
    \]
    where $b$ satisfies $b(s,0)=0$ and $|b|<1$. Since  $T$ and hence $V$ have unit length, it follows that $X_c$ is unit-length. The condition that it is divergence-free is equivalent to solving the quasi-linear initial value problem 
     \begin{equation}
        \label{PDEforb}
        \begin{cases}
        \p_t b -\frac{b}{\sqrt{1-b^2}}V\cdot\nabla b=-\sqrt{1-b^2}\operatorname{div}V+b\operatorname{div}\nu,\\
        b(s,0)=0.
    \end{cases}
    \end{equation}
    To obtain existence near $\Sigma$ we have to work in rescaled coordinates. Let $s_0\in\p N$ and rescale the coordinates through
    \[
    z=\frac{s-s_0}{\sigma_\partial(s_0)},\quad \hat t=\frac{t}{\sigma_\partial(s_0)}.
    \]
    It is not difficult to see that we can choose $c>0$ independent of $s_0$ such that
    \begin{equation}
        \label{eq:Unifcs0}
        |z|<4c\Rightarrow \frac{1}{2}\leq
        \frac{\sigma_\partial(s_0+\sigma_\partial(s_0)z)}{\sigma_\partial(s_0)}\leq2.
    \end{equation}
    Fix a sufficiently small $0<\theta_0<1/4$, to be decreased below. Set
    \[
    Q_\rho=\{(z,\hat t):|z|<\rho,\ 0\leq\widehat t\leq2\theta_0\},\quad 0<\rho\leq4c.
    \]
    In norms on $Q_\rho$, functions are composed with $(z,\hat t)\mapsto(s_0+\sigma_\partial(s_0)z,\sigma_\partial(s_0)\hat t)$, and $C^k(Q_\rho)$ denotes the ordinary $C^k$ norm.
    Using the same notation for the composed functions and since $V\cdot\nabla t=0$, equation \eqref{PDEforb} becomes
    \[
    \begin{cases}
        \partial_{\hat t}b-\dfrac{b}{\sqrt{1-b^2}}
        \displaystyle\sum_{i=1}^2(V\cdot\nabla s^i)\partial_{z^i}b
        =\sigma_\partial(s_0)\bigl(-\sqrt{1-b^2}\operatorname{div}V+b\operatorname{div}\nu\bigr),\\
        b(z,0)=0.
    \end{cases}
    \]
    Near a point $p \in \Sigma$ we have by definition
    \begin{equation}
        \label{eq:RatioBound}
        \frac{\sigma_\p(s_0+\sigma_\p(s_0)z)}{\sigma_\p(s_0)}=\left|\frac{s_0}{|s_0|}+|s_0|^\delta z\right|^{1+\delta}.
    \end{equation}
    Together with \eqref{eq:Unifcs0}, this means that all derivatives of \eqref{eq:RatioBound} and of its inverse are uniformly bounded. Away from $\Sigma$ this is clear by compactness.\\
    By definition of $\theta$ and the rescaled coordinates, we have
    \[
    \theta(s,t)=\frac{t}{\sigma_\p(s)}=\hat t\,\frac{\sigma_\partial(s_0)}{\sigma_\partial(s_0+\sigma_\partial(s_0)z)},
    \]
    so the preceding bounds on the ratio \eqref{eq:RatioBound} give $||\theta||_{C^k(Q_{4c})}\leq C_k$.
    Since $T$ is extended constantly in $t$, its components satisfy $\partial_{\hat t}T=0$. Moreover, near $\Sigma$, \eqref{eq:BoundUnitLength} and \eqref{eq:Unifcs0} give
    \[
    |\partial_z^\alpha T|
    \leq C_\alpha\sigma_\partial(s_0)^{|\alpha|}|s|^{-|\alpha|}
    \leq C_\alpha |s_0|^{\delta |\alpha|}
    \leq C_\alpha
    \]
    For $\nu$ it is clear that $||\nu||_{C^k(Q_{4c})}\leq C_k$. Applying the chain and product rules to
    $V$, we therefore derive
    \[
    ||V||_{C^k(Q_{4c})} \leq C_k||T||_{C^k(Q_{4c})} \bigl(1+||\nu||_{C^k(Q_{4c})}\bigr)\leq C_k.
    \] 
    Since the geometry is fixed we hence have
     \begin{equation}
        \label{eq:BoundRHSChar}
        ||V||_{C^k(Q_{4c})}+||\sigma_\p(s_0)\operatorname{div}V||_{C^k(Q_{4c})}+||\sigma_\p(s_0)\operatorname{div}\nu||_{C^k(Q_{4c})}\leq C_k.
    \end{equation}
    The coefficient of $\partial_{\hat t}b$ in the rescaled equation is $1$, so the initial surface is noncharacteristic. The standard local existence theorem by characteristics, see \cite[Section~3.2.4]{Evans2010}, therefore applies. By \eqref{eq:BoundRHSChar}, its coefficients and all their derivatives are uniformly bounded for $|b|\leq1/2$. Thus the existence time and the smooth dependence estimates can be chosen independently of $s_0$. For $\theta_0$ sufficiently small we obtain a  unique smooth solution on $Q_{2c}$ satisfying
    \begin{equation}
        \label{eq:SolOnCyl}
        |b|\leq C\hat t,\qquad ||b||_{C^k(Q_{2c})}\leq C_k.
    \end{equation}
    Consider the set
    \[
    \{(s,t):|s-s_0|<c\sigma_\partial(s_0),\;0\le t\le \theta_0\sigma_\partial(s)\}.
    \]
    By \eqref{eq:Unifcs0}, the rescaling maps this set into $Q_{2c}$, where the solution satisfies \eqref{eq:SolOnCyl}. Uniqueness along characteristics patches the local solutions to a smooth solution $b$ for \eqref{PDEforb} on $\mc C_{\theta_0}$.\\
    \textbf{Step 2: Estimates.} From \eqref{eq:SolOnCyl} we have
     \begin{equation}
        \label{eq:Boundsb}
        |b|\leq C\theta,\qquad\sigma_\partial^j |D^j b|\leq C_j,\quad j\geq 0.
    \end{equation}
    At this point we have constructed $b$ and hence $X_c$. It remains to define the potential $B_c$ with $\curl B_c=X_c$. This can be achieved by solving the initial value problem
    \begin{equation}
        \label{eq:ODEBc}
        \begin{cases}
            \partial_tB_c-(B_c\cdot\nabla)\nu=\nu\times X_c,\\
            B_c(s,0)=0.
        \end{cases}
    \end{equation}
    Since this is a linear ODE with smooth coefficients it is clear that a smooth solution exists. We check that $\curl B_c=X_c$. Set $h=(\curl B_c-X_c)\cdot\nu$. Since $B_c=0$ on $\partial N$, its tangential derivatives vanish there. From \eqref{eq:ODEBc} one checks that $B_c\cdot\nu=0$. Hence by a simple computation
    \begin{align*}
    \nu\times\curl B_c&=\partial_tB_c-(B_c\cdot\nabla)\nu=\nu\times X_c,\nonumber\\
    \curl B_c-X_c
    &=h\nu.
    \end{align*}
    Note that $h$ solves
    \begin{equation}
    \begin{cases}
        \partial_t h-h\operatorname{div}\nu=-\operatorname{div}(h\nu)
    =-\operatorname{div}(\curl B_c-X_c)=0,\\
    h(s,0)=\left.(\curl B_c)\cdot\nu\right|_{t=0}+b(s,0)=0.
    \end{cases} 
    \end{equation}
    Consequently, $h\equiv 0$ and therefore
    \begin{equation*}
        \curl B_c=X_c.
    \end{equation*}
    By \eqref{eq:BoundRHSChar} and \eqref{eq:SolOnCyl}, the $C^k$ norms of $X_c$ in the rescaled coordinates are uniformly bounded. Writing down the ODE for $B_c$ in rescaled coordinates as above one obtains from standard ODE estimates that
    \[
    ||X_c||_{C^k(Q_{2c})}+ \sigma_\partial(s_0)^{-1}||B_c||_{C^k(Q_{2c})}\leq C_k.
    \]
    Using $\sigma_\partial(s_0)\simeq \sigma_\partial(s)$, this shows
     \[
    \sigma_\partial^k|D^k X_c|+\sigma_\partial^{k-1}|D^kB_c|\leq C_k.
    \]
    On $\mathcal C_{\theta_0}$, we have $\sigma\simeq\sigma_\partial$, uniformly for sufficiently small $\theta_0$. Hence the preceding estimates give the norm bound in \eqref{eq:BoundsXc}.\\
    At $t=0$, we have $X_c=T$, $\partial_tX_c=\sigma_\partial^{-1}(\nu\times T)-(\partial_tb)\nu$.
    Consequently,
    \begin{align}
    \label{eq:Helatt0}
    \left.\sigma_\partial X_c\cdot\curl X_c\right|_{t=0}
    &=\sigma_\partial\bigl(-[(\nu\cdot\nabla)X_c]\cdot(\nu\times T)+[((\nu\times T)\cdot\nabla)X_c]\cdot\nu
      \bigr) \nonumber\\
    &=\sigma_\partial\bigl(\partial_tX_c\cdot(\nu\times T)
      +[((\nu\times T)\cdot\nabla)X_c]\cdot\nu\bigr) \nonumber\\
    &=1+\sigma_\partial[((\nu\times T)\cdot\nabla)T]\cdot\nu \nonumber\\
    &=1-\sigma_\partial T\cdot[((\nu\times T)\cdot\nabla)\nu] \nonumber\\
    &\geq1-\sigma_\partial|D\nu|\geq\tfrac34,
    \end{align}
    using the definition of $\sigma_\partial$ in the final inequality. The derivative estimates already proved and $\partial_t\sigma_\partial=0$ imply
    \[
    \left|\partial_t\bigl(\sigma_\partial X_c\cdot\curl X_c\bigr)\right|
    \leq C\sigma_\partial^{-1}.
    \]
    Since $t\leq\theta_0\sigma_\partial$ on $\mc C_{\theta_0}$, integration in the $t$ direction gives together with \eqref{eq:Helatt0}, after decreasing $\theta_0$,
    \[
    \sigma_\partial X_c\cdot\curl X_c\geq\tfrac34-C\frac{t}{\sigma_\partial}\geq\tfrac34-C\theta_0\geq\tfrac12
    \quad\text{on }\mathcal C_{\theta_0}.
    \]
    Using again $\sigma\simeq\sigma_\partial$, we obtain, for some $\kappa_c>0$,
    \[
        q_{X_c} \geq\kappa_c\quad\text{on }\mathcal C_{\theta_0}.
    \]
    This completes the proof.
\end{proof}
Next we glue the already constructed fields to obtain a field in a neighborhood of $\p U$.
\begin{lemma}
    \label{lem:Matching}
    We use the notation of Lemmas \ref{lem:FormOnEnd}--\ref{lem:UnitCollar}. Define $K\subset N$ by
    \[
    K:=\mathcal C_{\theta_0}\cup\bigcup_{p \in \Sigma}U_\eps(p).
    \]
    For $\varepsilon>0$ sufficiently small and $\theta_0>0$ sufficiently small, there  exists $B_g\in C^\infty(K;\mathbb R^3)$ such that $X_g=\curl B_g$ is nowhere vanishing and
     \[
    X_g=\begin{cases}
        X_c,\quad \text{on }\mathcal C_{\theta_0/4}\\
        X_e,\quad \text{on }K\setminus \mathcal C_{\theta_0}.
    \end{cases}
    \]
    Moreover, we have
     \begin{equation}
        \label{eq:OmgBound}
        ||X_g||_{k;K}
        +||B_g||_{k,1;K}
        \leq C_k,\quad \text{for all }k\geq 0.
    \end{equation}
\end{lemma}

\begin{proof}
    After shrinking $U_\eps(p)$ if necessary, suppose that the construction of Lemma \ref{lem:FormOnEnd} is available on slightly larger concentric neighborhoods $\widetilde U_\eps(p)$, with $T=e_2$ on $\widetilde U_\eps(p)\cap\partial N$.\\
    We claim first that
     \begin{equation}
        \label{eq:betacandeFirst}
        (B_c-B_e)|_{t=0}=0,\quad D(B_c-B_e)|_{t=0}=0.
    \end{equation}
    Indeed, both  potentials vanish at $t=0$, so their tangential derivatives vanish  there as well. Since $n=-\nu$ and  $T=e_2$ on the boundary, their definitions give
    \[
    \partial_tB_e|_{t=0}=\frac{\partial_t\varphi}{|\nabla\varphi|}e_1
    =e_1=\nu\times T=\partial_tB_c|_{t=0},
    \]
     proving \eqref{eq:betacandeFirst}.\\
    For every $p \in \Sigma$, we have
    \[
    t\leq \theta_0|s|^{1+\delta},\quad
    r\simeq |s|,\quad
    \sigma_\partial=|s|^{1+\delta},\quad \text{on }\mc C_{\theta_0}\cap U_\eps(p).
    \]
    Differentiating the formula for  $B_e$ in Lemma~\ref{lem:FormOnEnd}, using \eqref{eq:FrameBound}, gives
     \[
    |D^jB_e|
    \leq C_jr^{(1+\delta)(1-j)}
    \leq C_j\sigma_\partial^{1-j},\qquad j\geq 0.
    \]
    Together with \eqref{eq:BoundsXc}, this yields
     \begin{equation}
        \label{eq:DerBoundsbetae}
        |D^jB_e|+|D^jB_c|
    \leq C_j\sigma_\partial^{1-j},\qquad j\geq 0.
    \end{equation}
    We first extend  $B_e$ to $\mc C_{\theta_0}$. Choose a smooth cutoff $\psi$ equal to $1$ on a neighborhood of $\bigcup_{p\in\Sigma}U_\eps(p)$ and supported in the enlarged neighborhoods $\bigcup_{p \in \Sigma} \tilde U_\eps(p)$. On $\mc C_{\theta_0}$, replace  $B_e$ by
    \[
    \psi B_e+(1-\psi)B_c,
    \]
    where  $\psi B_e$ is extended by zero outside the larger neighborhoods. We retain the notation  $B_e$ for this extension. It agrees with the original  $B_e$ on each $U_\eps(p)$. Moreover, its difference from  $B_c$ is $\psi(B_e-B_c)$, so \eqref{eq:betacandeFirst} holds on $\p N$. The preceding derivative bounds also remain valid: derivatives of $\psi$ are bounded and supported away from
    $\Sigma$, and we may assume $\sigma_\partial\leq 1$.\\
    Since $\sigma_\partial(s)$ is constant along each normal  line, Taylor's theorem, \eqref{eq:betacandeFirst} and \eqref{eq:DerBoundsbetae} with $j=2$  give
    \begin{equation}
    \label{eq:Integratedbetacande}
    |D(B_c-B_e)(s,t)|\leq C\frac{t}{\sigma_\partial(s)},\qquad
    |(B_c-B_e)(s,t)|\leq C\frac{t^2}{\sigma_\partial(s)}.
    \end{equation}

    Let $\chi=\chi(t/(\theta_0\sigma_\partial))$ be a smooth cutoff satisfying
     \[
    0\leq\chi\leq 1,\qquad
    \chi=1\quad\text{for }t\leq\theta_0\sigma_\partial/4,\qquad
    \chi=0\quad\text{for }t\geq\theta_0\sigma_\partial/2.
    \]
     The chain rule and the bounds for $\sigma_\partial$ give
    \[
    |\nabla\chi|\leq\frac{C}{\theta_0\sigma_\partial},\qquad
    |D^j\chi|\leq C_{j,\theta_0}\sigma_\partial^{-j}.
    \]
     Define
     \[
    B_g=
    \begin{cases}
    B_e+\chi(B_c-B_e),
    &\text{on }\mathcal C_{\theta_0},\\
    B_e,
    &\text{on }K\setminus\mathcal C_{\theta_0}.
    \end{cases}
    \]
    This defines a smooth  vector potential because $\chi=0$ near $\p_N\mc C_{\theta_0}$ and the extended  $B_e$ agrees with the original  $B_e$ on every $U_\eps(p)$.
    Thus  $X_g=\curl B_g$ is smooth and  divergence-free on $K$. We obtain
     \[
    X_g=X_c\quad\text{on }\mathcal C_{\theta_0/4},
    \qquad
    X_g=X_e\quad\text{on }
    K\setminus\mathcal C_{\theta_0}.
    \]
    To prove \eqref{eq:OmgBound}, note that
     \[
    \curl B_g-\curl B_e
    =\chi\,\curl(B_c-B_e)
    +\nabla\chi\times(B_c-B_e),
    \]
    which is zero if $t\geq\theta_0\sigma_\partial/2$. For $t\leq\theta_0\sigma_\partial/2$,
    \eqref{eq:Integratedbetacande} gives
     \[
    |\curl B_g-\curl B_e|
    \leq C\frac{t}{\sigma_\partial}
    +\frac{C}{\theta_0\sigma_\partial}\frac{t^2}{\sigma_\partial}
    \leq C\theta_0.
    \]
    Consequently, on the original neighborhoods,
     \[
    |X_g-X_e|\leq C\theta_0,\qquad
    \bigl||X_g|-1\bigr|
    \leq C\theta_0+C\varepsilon^\delta,
    \]
    where the second inequality follows from Lemma \ref{lem:FormOnEnd}.\\
    For higher derivatives, the  product rule and the bounds for  $\chi,B_e,B_c$
    yield
    \[
    \sigma_\partial^{j-1}|D^jB_g|
    +\sigma_\partial^j|D^jX_g|
    \leq C_{j,\theta_0},\quad j\geq 0,\quad \text{in }\mc C_{\theta_0}.
    \]
    Outside $\mc C_{\theta_0}$, the  potential $B_e$ is
    unchanged, so together with the derivative bounds in Lemma \ref{lem:FormOnEnd} this proves \eqref{eq:OmgBound}.\\
    Finally, in $\mc C_{\theta_0}$,
     \[
    X_g-X_c
    =-(1-\chi)\curl(B_c-B_e)
    +\nabla\chi\times(B_c-B_e).
    \]
    Using \eqref{eq:Integratedbetacande} once more, we obtain
     \[
    |X_g-X_c|
    \leq C\theta_0.
    \]
    Since  $|X_c|=1$, choosing $\theta_0$ sufficiently small in Lemma \ref{lem:UnitCollar} ensures  $|X_g|\geq 1/2$ on $\mc C_{\theta_0}$. In $K \setminus \mc C_{\theta_0}$,  $X_g=X_e$ is nowhere vanishing for $\varepsilon$ sufficiently small, by Lemma \ref{lem:FormOnEnd}.
\end{proof}
The idea of the next argument is that given a  nowhere vanishing vector field on a compact  domain which is a curl near the boundary, by adding a finite number of the building blocks from the introduction we can obtain a nowhere vanishing curl, without changing the boundary data. This statement can also be shown using a relative $h$-principle. See,
for example, \cite[Theorem 2.1 and Lemma 2.4]{McDuff1987}. As the following argument uses the same idea as the proof of the main theorem, we give this more analytic proof here.
\begin{lemma}
    \label{lem:RelativeExactExtension}
    Let  $W\subset\mathbb R^3$ be a compact  domain with smooth boundary. Let $X\in C^\infty(W;\mathbb R^3)$ be nowhere vanishing on $W$ such that  $X=\curl B_K$ near $\p W$ for a smooth potential $B_K$ defined near $\p W$.\\
    Then there exists  $B\in C^\infty(W;\mathbb R^3)$ such that $\curl B$ is nowhere vanishing and  $B=B_K$ near $\p W$.
\end{lemma}

\begin{proof}
    Without loss of generality we can rescale such that $\min_W|X|=1$. Smoothly extend $B_K$ to $W$, without relabeling. Denote
    \[
    R=X-\curl B_K.
    \]
    If $R=0$, there is nothing to prove. In the other case choose a nonnegative
    $\chi\in C_c^\infty(W^\circ)$ such that $\chi>|R|$ on $\supp R$.\\
    Recall the building block $Y$ from the introduction:
    Let $V\in C^\infty(W;\mb R^3)$, $\xi\in\mb R^3$ with $|\xi|=1$
    and a nonnegative $\rho\in C^\infty_c(W)$. Provided
    $b=|\Pi_{\xi^\perp}V|>0$ near $\supp\rho$, denote
    $e=\Pi_{\xi^\perp}V/b$. If $\rho/b\leq\gamma$, define near $\supp\rho$
    \begin{equation}
        \label{eq:DefY}
        Y=V-be+\sqrt{b^2+\rho^2}
        \bigl(\cos(a\sin\vartheta)e+
        \sin(a\sin\vartheta)(\xi\times e)\bigr),
    \end{equation}
    and set $Y=V$ elsewhere. Then
    \[
    \overline Y=V,\qquad |Y|^2=|V|^2+\rho^2,\qquad
    (Y-V)\cdot\xi=0,\qquad |Y-V|\leq C\rho.
    \]
    For the precise construction, see Step 1 in Lemma \ref{lem:LocalPert}.\\
    Write $w=Y-V$ and write as in \eqref{eq:PotentialDefi}
    \[
    \mc P_\xi w=-\xi\times\p_\vartheta^{-1}w,\qquad
    A_\lambda(x)=\lambda^{-1}
    (\mc P_\xi w)(x,\lambda\xi\cdot x).
    \]
    From the chain rule,
    \begin{equation}
        \label{eq:PerturbationThroughALambda}
        V+\curl A_\lambda
        =Y(x,\lambda\xi\cdot x)+O(\lambda^{-1}),\qquad
        \supp A_\lambda\subset\supp\rho.
    \end{equation}
    Fix $J\in\mb N$ such that
    \[
    \frac{\|\chi\|_{C^0}}{\sqrt J}
    \leq\min\left\{\frac{\gamma}{8},\frac{1}{16C}\right\}.
    \]
    Starting with $X_0=X$, we construct inductively $X_1,\ldots,X_{2J}$.\\
    Suppose that $|X_{2j}|\geq1/2$. As in the cutoff construction
    in Lemma \ref{lem:CorrectionDirection}, it is simple to choose nonnegative
    $\rho_{2j+1},\rho_{2j+2}\in C^\infty_c(W^\circ)$ such that
    \begin{equation}
        \label{eq:AmplitudeProp}
        \rho_{2j+1}^2+\rho_{2j+2}^2=\frac{\chi^2}{J},\qquad
        \supp\rho_{2j+i}\subset\supp\chi,\qquad
        |\Pi_{\mathbf e_i^\perp}X_{2j}|\geq\frac14
        \quad\text{on }\supp\rho_{2j+i}.
    \end{equation}
    In particular, $||\rho_{2j+i}||_{C^0}\leq||\chi||_{C^0}/\sqrt J$.
    Hence, by the choice of $J$, we can construct $Y$ with data $V=X_{2j}$, $\rho=\rho_{2j+1}$, $\xi=\mathbf e_1$,
    obtaining $A_\lambda$ as above. We define
    \[
    X_{2j+1}=X_{2j}+\curl(A_\lambda),
    \]
    with $\lambda$ large enough such that
    \begin{equation}
        \label{eq:ChangeFirstCorr}
        ||X_{2j+1}-X_{2j}||_{C^0}
        \leq\frac{C||\chi||_{C^0}}{\sqrt J}+\frac1{16}
        \leq\frac18,
    \end{equation}
    and $\|E_{2j+1}\|_{C^0}\leq1/(4J)$, where
    \[
    E_r=|X_r|^2-|X_{r-1}|^2-\rho_r^2.
    \]
    Hence, on $\supp\rho_{2j+2}$,
    \[
    |\Pi_{\mathbf e_2^\perp}X_{2j+1}|
    \geq\frac14-\frac18=\frac18.
    \]
    Thus, we can make a second perturbation, this time with data
    $V=X_{2j+1}$, $\rho=\rho_{2j+2}$, $\xi=\mathbf e_2$,
    choosing again its frequency sufficiently large so that
    $\|E_{2j+2}\|_{C^0}\leq1/(4J)$ and obtaining $X_{2(j+1)}$.
    Summing over the steps already completed gives
    \[
    |X_{2(j+1)}|^2
    \geq |X|^2-\sum_{r=1}^{2(j+1)}||E_r||_{C^0}
    \geq1-\frac{2(j+1)}{4J}\geq\frac12.
    \]
    In particular, $|X_{2(j+1)}|\geq1/2$. Using $X_{2(j+1)}$ we compute the next pair of cutoffs as in \eqref{eq:AmplitudeProp}.
    Iterating this procedure we obtain
    \[
    X_n=X_{n-1}+\curl A_n,\qquad
    |X_n|^2=|X_{n-1}|^2+\rho_n^2+E_n,\qquad
    \|E_n\|_{C^0}\leq\frac1{4J},
    \]
    with $\supp A_n\subset\supp\chi$. Moreover,
    \[
    |X_n|^2
    =|X|^2+\sum_{r=1}^n\rho_r^2+\sum_{r=1}^nE_r.
    \]
    After $2J$ steps,
    \begin{equation}
        \label{eq:InequalityLength}
        |X_{2J}|^2\geq\frac12+\chi^2>|R|^2
        \quad\text{on }\supp R.
    \end{equation}
    Define
    \[
    B=B_K+\sum_{n=1}^{2J}A_n.
    \]
    Then $\curl B=X_{2J}-R$ is nowhere vanishing on $\supp R$
    by \eqref{eq:InequalityLength} and clearly also outside $\supp R$.
    Moreover, since all corrections are supported in
    $\supp\chi\subset\subset W^\circ$, we have $B=B_K$ near $\p W$.
    This completes the proof.
\end{proof}

 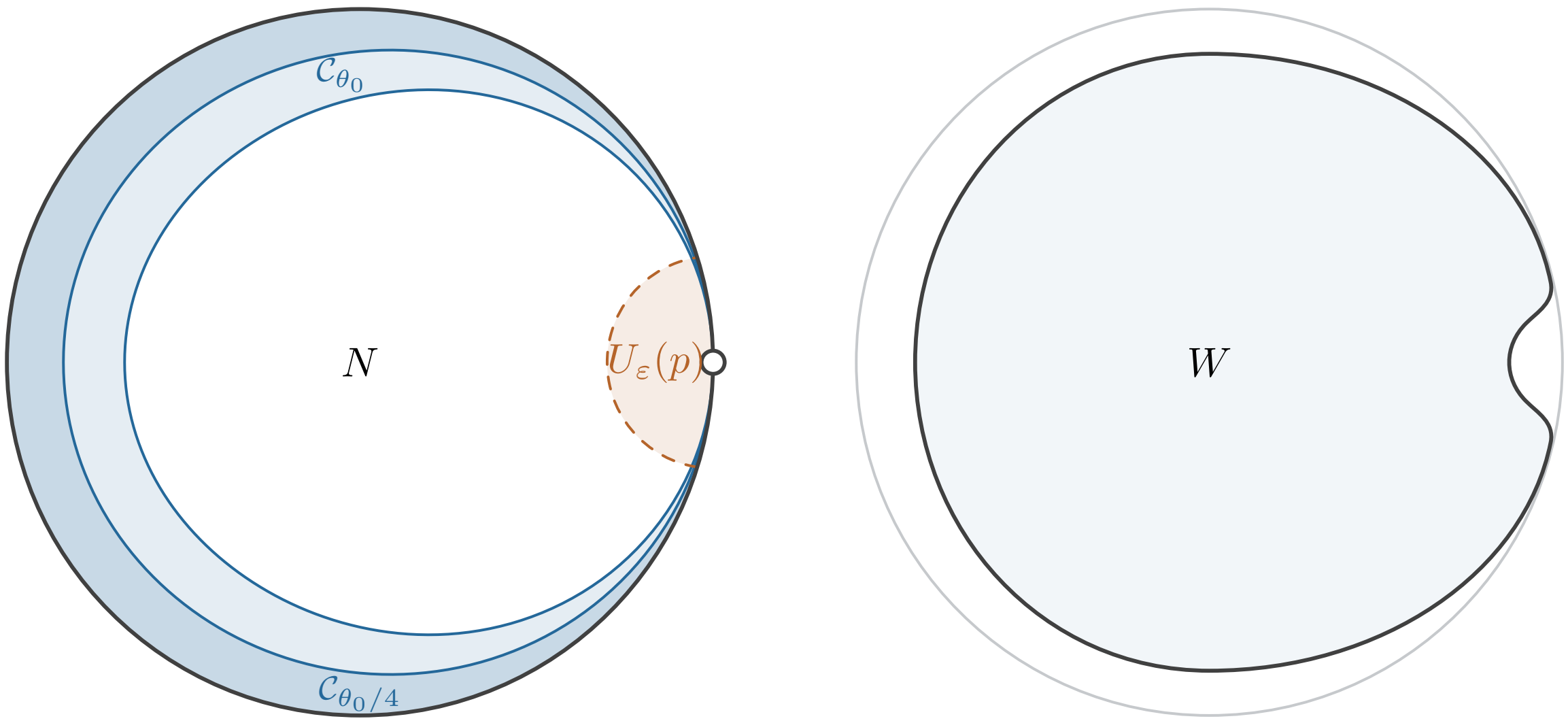
\begin{figure}[htbp]
\centering
 \resizebox{0.47\textwidth}{!}{%
\begin{tikzpicture}[font=\small,line cap=round,line join=round]
    \path[use as bounding box]
        (-3.2,-3.15) rectangle (3.2,3.15);

    \begin{scope}
        \clip (0,0) circle[radius=3];
        \fill[errororange!12] (3,0) circle[radius=0.9];
    \end{scope}

    \foreach \a/\shade in {1.0/12,0.48/25} {
        \fill[profileblue!\shade,even odd rule]
            (0,0) circle[radius=3]
            plot[domain=0:360,samples=181]
                ({(3-\a*(1-cos(\x/2)^4))*cos(\x)},
                 {(3-\a*(1-cos(\x/2)^4))*sin(\x)})
            -- cycle;

        \draw[profileblue,line width=0.65pt]
            plot[domain=0:360,samples=181]
                ({(3-\a*(1-cos(\x/2)^4))*cos(\x)},
                 {(3-\a*(1-cos(\x/2)^4))*sin(\x)})
            -- cycle;
    }

    \begin{scope}
        \clip (0,0) circle[radius=3];
        \draw[errororange,dashed,line width=0.65pt]
            (3,0) circle[radius=0.9];
    \end{scope}

    \draw[black!75,line width=0.95pt]
        (0,0) circle[radius=3];
    \draw[fill=white,draw=black!75,line width=0.95pt]
        (3,0) circle[radius=2.8pt];

    \node at (0,0) {$N$};
    \node[text=errororange] at (2.52,0) {$U_\eps(p)$};
    \node[font=\scriptsize,text=profileblue,inner sep=0pt]
    at (-0.15,2.44) {$\mathcal C_{\theta_0}$};
\node[font=\scriptsize,text=profileblue,inner sep=0pt]
    at (0,-2.82) {$\mathcal C_{\theta_0/4}$};
\end{tikzpicture}%
}\hfill
\resizebox{0.47\textwidth}{!}{%
\begin{tikzpicture}[font=\small,line cap=round,line join=round]
    \path[use as bounding box]
        (-3.2,-3.15) rectangle (3.2,3.15);

    \draw[figuregray!55,line width=0.65pt]
        (0,0) circle[radius=3];

    \path[fill=profileblue!6,draw=black!75,line width=0.95pt]
        (0,2.62)
        .. controls (1.48,2.62) and (2.68,1.78)
            .. (2.90,0.68)
        .. controls (2.934,0.51) and (2.79,0.43)
            .. (2.69,0.33)
        .. controls (2.50,0.14) and (2.50,-0.14)
            .. (2.69,-0.33)
        .. controls (2.79,-0.43) and (2.934,-0.51)
            .. (2.90,-0.68)
        .. controls (2.68,-1.78) and (1.48,-2.62)
            .. (0,-2.62)
        .. controls (-1.42,-2.62) and (-2.50,-1.43)
            .. (-2.50,0)
        .. controls (-2.50,1.43) and (-1.42,2.62)
            .. (0,2.62)
        -- cycle;

    \node at (0,0) {$W$};
\end{tikzpicture}%
}
\caption{Collars and singularity (left), the compact domain $W$ (right)}
\label{fig:CollarsAndCompactCore}
\end{figure}

\begin{proof}[Proof of Proposition \ref{prop:Prep}]
We split the proof into three steps.\\
\textbf{Step 1: Filling.}
Let $\theta_0>0$ be as in Lemma \ref{lem:Matching}. In the first step, we construct a nowhere vanishing $X_f=\curl B_f$ with $X_f\in C^\infty(N;\mathbb R^3)$ such that
\[
    X_f|_{\mathcal C_{\theta_0/4}}=X_c|_{\mathcal C_{\theta_0/4}}.
\]
Fix $\varepsilon>0$ small enough such that the sets $U_\eps(p)$ are pairwise disjoint and $\sigma_\partial(s)=|s|^{1+\delta}$ on $U_\eps(p)\cap\p N$.
Pick a smooth compact domain $W\subset\mathring N$ and set
\[
    \tilde K=N\setminus\mathring W,
\]
such that
\[
    \begin{gathered}
    \mathcal C_{\theta_0/4}\cup
    \bigcup_{p\in\Sigma}U_{\varepsilon/3}(p)
    \subset\tilde K\subset K,\\
    \p W\subset\{0\leq\theta<3\theta_0/4\}\cup\bigcup_{p\in\Sigma}U_{2\varepsilon/3}(p)
    \subset\operatorname{int}_N K,
    \end{gathered}
\]
where $K$ is as in Lemma \ref{lem:Matching}. See Figure \ref{fig:CollarsAndCompactCore} for an illustration of the situation. \\
In particular, Lemma \ref{lem:Matching} gives the smooth nowhere vanishing field $X_g=\curl B_g$ on $\tilde K$ and on a neighborhood of $\p W$. We first extend $X_g$ to a smooth nowhere vanishing vector field $\widetilde X_g$ on $W$ satisfying
\[
    \widetilde X_g=X_g\quad\text{near }\p W.
\]
To do so, view $f=X_g/|X_g|$ as a map  $\p W\to\mathbb S^2$. We will show that its restriction to each boundary component is null-homotopic. For this, we construct a nowhere vanishing  vector field $V$  perpendicular to $f$ on $\p W$.\\
Let  $n=\nabla\varphi/|\nabla\varphi|$ on each $U_\eps(p)$, as in Lemma \ref{lem:SmoothFrameInd}. Consider
\[
    \mathcal A=\p W\cap\{0\leq\theta<3\theta_0/4\},\qquad\mathcal B=\p W\cap\bigcup_{p\in\Sigma}U_\eps(p).
\]
By construction, $\mathcal A\cup\mathcal B=\p W$. Hence, $\nu$ is defined on $\mathcal A$, and $n$ is defined on $\mathcal B$. We now verify that the projections of  $-\nu$ and $n$ are nonzero on their respective domains.\\
 The construction in Lemma \ref{lem:UnitCollar} and
\eqref{eq:Boundsb} give
 \begin{equation*}
    |X_c\cdot\nu|=|b|\leq C\theta\leq C\theta_0.
\end{equation*}
Moreover, the proof of Lemma \ref{lem:Matching} gives
 $|X_g|\geq 1/2$ and
$|X_g-X_c|\leq C\theta_0$ throughout $\mc C_{\theta_0}$.
Consequently,
 \begin{equation}
    \label{eq:fnu}
    |f\cdot\nu|
    \leq C\theta_0
    \quad\text{on }\mathcal A.
\end{equation}
On $U_\eps(p)$, the  leading term in \eqref{eq:FormBe} is orthogonal to $n$, while the remainder has norm at most $Cr^\delta$. Hence  $|X_e\cdot n|\leq Cr^\delta\leq C\varepsilon^\delta$. The construction of Lemma \ref{lem:Matching} changes  $X_e$ by at most  $C\theta_0$ on $U_\eps(p)\cap\mathcal C_{\theta_0}$, and leaves it unchanged outside $\mc C_{\theta_0}$. Since  $|X_g|$ is uniformly bounded
away from zero, we obtain
 \begin{equation}
    \label{eq:fn}
    |f\cdot n|
    \leq C(\theta_0+\varepsilon^\delta)
    \quad\text{on }\mathcal B.
\end{equation}
Choose $\varepsilon$ small enough and then $\theta_0$ small enough such that both right hand sides of \eqref{eq:fnu}--\eqref{eq:fn} are at most $1/2$. Since $\nu$ and $n$ are unit, we obtain
 \begin{equation}
    \label{eq:ProjLowBd}
    |\Pi_{f^\perp}\nu|^2=1-(f\cdot\nu)^2
    \geq\frac{3}{4}\quad\text{on }\mathcal A,
    \qquad
    |\Pi_{f^\perp}n|^2=1-(f\cdot n)^2
    \geq\frac{3}{4}\quad\text{on }\mathcal B.
\end{equation}
On $\mathcal A\cap\mathcal B$, the identity  $n(s,0)=-\nu(s,0)$,
the bound  $|Dn|\leq Cr^{-1}$ from
Lemma \ref{lem:SmoothFrameInd}, and  $\partial_t\nu=0$ give
\[
    |n(s,t)+\nu(s,t)|
    \leq\int_0^t|\partial_un(s,u)|\,du
    \leq C\frac{t}{|s|}
    \leq C\theta_0|s|^\delta
    \leq C\theta_0\varepsilon^\delta,
\]
using that $\sigma_\partial(s)=|s|^{1+\delta}$ near each $p \in \Sigma$. Since orthogonal projection has norm one, it follows that
 \begin{equation}
    \label{eq:projectionsClose}
    |\Pi_{f^\perp}n+\Pi_{f^\perp}\nu|
    \leq |n+\nu|\leq C\theta_0\varepsilon^\delta
    \quad\text{on }\mathcal A\cap\mathcal B.
\end{equation}
Choose a smooth partition of unity $(\chi,1-\chi)$ on
$\p W$ subordinate to $(\mathcal A,\mathcal B)$, so that
\[
    0\leq\chi\leq1,\qquad
    \supp\chi\subset\mathcal A,\qquad
    \supp(1-\chi)\subset\mathcal B.
\]
Define on $\p W$
 \[
    V=-\chi\Pi_{f^\perp}\nu+(1-\chi)\Pi_{f^\perp}n,
\]
where we extend the summands by zero outside $\mathcal A$ and $\mathcal B$, respectively. The support conditions make $V$ smooth across the boundaries of these sets. We have  $V\cdot f=0$ everywhere, and on $\mc A\cap \mc B$ we have by \eqref{eq:projectionsClose}
 \[
    |V|\geq |\Pi_{f^\perp}\nu|
    -(1-\chi)|\Pi_{f^\perp}n+\Pi_{f^\perp}\nu|
    \geq\frac{\sqrt3}{2}-C\theta_0\varepsilon^\delta>0.
\]
Together with \eqref{eq:ProjLowBd}, this shows that $V$ is nowhere vanishing for $\varepsilon$ and $\theta_0$ sufficiently small.\\
The homotopy
\[
    (x,s)\longmapsto
    \cos(\pi s)f(x)+\sin(\pi s)\frac{V(x)}{|V(x)|},
    \qquad 0\leq s\leq1,
\]
joins $f$ to $-f$. Since the antipodal map on $\mb S^2$ has degree $-1$, the restriction of $f$ to each boundary component has degree zero. The Hopf degree theorem therefore implies that each restriction is null-homotopic.\\
From this it is easy to construct  $\widetilde X_g\in C^\infty(W;\mathbb R^3)$ such that
\[
    |\widetilde X_g|>0\quad\text{on }W,
    \qquad
    \widetilde X_g=X_g\quad\text{near }\p W.
\]
 Apply Lemma \ref{lem:RelativeExactExtension} to  $\widetilde X_g$, with the prescribed  potential $B_g$ near $\p W$. We obtain a smooth  potential $B_f$ on $W$ such that
 \[
    |\curl B_f|>0\quad\text{on }W,
    \qquad
    B_f=B_g\quad\text{near }\p W.
\]
Extend $B_f$ to $N$ by setting $B_f=B_g$ outside $W$.
The two definitions agree near $\p W$, so  $B_f$ is smooth. Moreover,  $X_f=\curl B_f$ is nowhere
vanishing and $X_f=X_g=X_c$ on $\mathcal C_{\theta_0/4}\subset\tilde K$. The estimates in
Lemmas \ref{lem:FormOnEnd}--\ref{lem:Matching} give uniform upper and positive lower bounds for  $|X_g|$ on $\tilde K$. On the compact set $W$, such bounds follow from smoothness and nonvanishing. Hence there are constants  $m,M>0$ such that
\begin{equation}
    \label{eq:prep-filled}
    0<m\leq|X_f|\leq M\quad\text{on }N,
    \qquad
    X_f|_{\mathcal C_{\theta_0/4}}=X_c|_{\mathcal C_{\theta_0/4}}.
\end{equation}
For these fixed choices, the derivative bounds in Lemmas \ref{lem:FormOnEnd}--\ref{lem:Matching}, together with compactness of $W$ and $\sigma\simeq\sigma_\partial$ in $\mathcal C_{\theta_0/4}$ near every $p \in \Sigma$, also give
 \[
    ||X_f||_k+||B_f||_{k,1}<\infty
    \quad\text{for }k\geq0.
\]
\textbf{Step 2: Rescaling.}
Next, we want to rescale the field such that it is strictly short away from $\mathcal C_{\theta_0/4}$. Since we want the defect to satisfy certain derivative estimates the rescaling has to be done with some care.\\
Fix
 \[
    0<\eta<\min\{1,(2M)^{-1}\}.
\]
We make the Ansatz
 \begin{equation*}
    B_0=\eta B_f+(a-\eta)B_c,\qquad X_0=\curl B_0=\eta X_f+\curl((a-\eta)B_c),
\end{equation*}
where  $a=a(\theta)$ will be chosen below, with  $a=1$ near $\theta=0$ and  $a=\eta$ for $\theta\geq\delta_c$, for some $0<\delta_c<\theta_0$. Here  $(a-\eta)B_c$ is extended
by zero outside $\mathcal C_{\theta_0}$. Where  $X_f=X_c$, this gives
 \begin{equation}
    \label{eq:Om0Comp}
    X_0=\curl(aB_c)=aX_c+a'\nabla\theta\times B_c.
\end{equation}
Since  $B_c$ vanishes at $\theta=0$, we write
 \[
    \Xi=\frac{\nabla\theta\times B_c}{\theta},
\]
so that the second term in \eqref{eq:Om0Comp} becomes  $\theta a'\Xi$. It is easy to see that $\Xi$ extends smoothly to $\theta=0$, with
 \[
    \Xi|_{\theta=0}=X_c|_{\theta=0},
    \qquad \Xi-X_c=O(\theta),
    \qquad
    ||\Xi||_{k;\mathcal C_{\theta_0}}\leq C_k,
    \quad \forall k\geq0.
\]
Choose $\delta_c>0$ sufficiently small that
 \begin{equation}
    \label{eq:prep-collar}
    X_f=X_c,\qquad |X_c|=1,\qquad
    X_c\cdot\Xi\geq\frac{3}{4},\qquad |\Xi|\leq2
    \quad\text{for }0\leq\theta\leq\delta_c.
\end{equation}
We will choose  $a$ nonincreasing, with  $|\theta a'|/a$ small, to control the derivative term in \eqref{eq:Om0Comp}. To obtain the required derivative estimates of the defect, choose a smooth nondecreasing cutoff $h:\mathbb R\to[0,1]$ such that
\begin{equation}
    \label{eq:prep-h}
    h(z)=
    \begin{cases}
        0, & z\leq0,\\
        z^2e^{-1/z}, & 0<z<\eps_h,\\
        1, & z\geq1,
    \end{cases}
\end{equation}
for some $\eps_h>0$ sufficiently small. The choice near zero gives $h'(z)=(1+2z)e^{-1/z}$, so both $h$ and $h'$ contain an exponential factor. We use a logarithmic change of variable to make  $|\theta a'|/a$ small. Let
\[
    L\geq\max\left\{1,
        \frac{8(1-\eta)||h'||_\infty}{\eta}\right\},
    \qquad
    \theta_*=\delta_c e^{-L},
    \qquad
    \tau=\frac{1}L\log\left(\frac{\theta}{\theta_*}\right),
\]
and define
 \[
    a(\theta)=1-(1-\eta)h(\tau),
\]
with  $a=1$ at $\theta=0$. These choices imply
 \begin{equation}
    \label{eq:prep-A}
    \begin{gathered}
        \eta\leq a\leq1,\qquad a'\leq0,\\
        a=1 \quad\text{for }\theta\leq\theta_*,
        \qquad
        a=\eta \quad\text{for }\theta\geq\delta_c,\\
        0\leq-\frac{\theta a'}{a}
          =\frac{(1-\eta)h'(\tau)}{La}
          \leq\frac{1}{8}.
    \end{gathered}
\end{equation}
With this choice of  $a$, the extension of  $(a-\eta)B_c$ by zero is smooth, since  $a-\eta$ vanishes to infinite order at $\theta=\delta_c$. By \eqref{eq:prep-collar},
 \begin{equation}
    \label{eq:prep-rescaled}
    X_0=
    \begin{cases}
        \curl(aB_c)=aX_c+\theta a'\Xi,
            & 0\leq\theta\leq\delta_c,\\
        \eta X_f,
            & \text{else}.
    \end{cases}
\end{equation}
\textbf{Step 3: Verification of the properties.}
The  field $X_0$ is smooth and  is a curl by construction. We now prove the properties claimed in the proposition.\\
For \ref{prop:it1}, on $\mathcal C_{\delta_c}$, \eqref{eq:prep-collar}--\eqref{eq:prep-rescaled} give
 \begin{equation}
    \label{eq:prep-norm}
    \begin{aligned}
        |X_0|^2
            &=a^2+2a\theta a'(X_c\cdot\Xi)
              +(\theta a')^2|\Xi|^2\\
            &\le a^2+\frac{3}{2}a\theta a'+4(\theta a')^2\\
            &\le a^2+a\theta a'
            \le a^2.
    \end{aligned}
\end{equation}
Here we used  $\theta a'\le0$ and $|\theta a'|\le a/8$. Similarly,
 \[
    |X_0|
    \ge a-|\theta a'|\,|\Xi|
    \ge a+2\theta a'
    \ge\frac{3}{4}a.
\]
On $N\setminus\mathcal C_{\delta_c}$, \eqref{eq:prep-filled} and \eqref{eq:prep-rescaled} give  $\eta m\le|X_0|\le\eta M<1/2$. Consequently, with $c_0=\eta\min\{m,3/4\}>0$,
 \[
    c_0\le|X_0|\le1
    \quad\text{on }N.
\]
Set $\mathcal C=\mathcal C_{\theta_*}$. By \eqref{eq:prep-A} and \eqref{eq:prep-rescaled},
 \begin{equation}
    \label{eq:prep-unit-collar}
    X_0|_{\mathcal C}=X_c|_{\mathcal C},
    \qquad
    |X_0|=1\quad\text{on }\mathcal C.
\end{equation}
Since $h(z)>0$ for $z>0$, we have  $a(\theta)<1$ for $\theta>\theta_*$. Thus \eqref{eq:prep-norm}, together with  $\eta M<1/2$ on $N\setminus\mathcal C_{\delta_c}$, gives
 \[
    |X_0|<1\quad\text{on }N\setminus\mathcal C.
\]
Hence $\{\Delta=0\}=\mathcal C$. It is simple to see that this implies
\[
    \dist_\sigma(\supp\Delta,\partial N)
    \geq\frac{\theta_*}{C}=:d_*>0.
\]
This completes the proof of \ref{prop:it1}.\\
For \ref{prop:it2}, \eqref{eq:prep-unit-collar} gives
 \[
    X_0\cdot\nu=X_c\cdot\nu=0,
    \qquad
    X_0=X_c=T\quad\text{on }\partial N.
\]
Hence the indices at points of $\Sigma$ are exactly those of $T$ in Corollary~\ref{coro:TangentFieldWithCorrectIndex}.\\
For \ref{prop:it3}, Lemma~\ref{lem:UnitCollar} allows us to choose $\kappa_0>0$ such that  $q_{X_c}\ge2\kappa_0$ on $\mathcal C_{\delta_c}$. On this set, \eqref{eq:prep-rescaled} gives
 \[
    X_0-X_c=(a-1)X_c+\theta a'\Xi.
\]
Both $(a-1)$ and $\theta a'$ vanish to every order at $\theta=\theta_*$. From the derivative bounds for $\theta$, $X_c$ and $\Xi$, this implies that
 \[
    ||X_0-X_c||_{1;\{\theta_*\le\theta\le \tilde \theta\}}
    \rightarrow 0
    \qquad\text{as }\tilde \theta\downarrow\theta_*.
\]
By continuity of $q_{X}$ in this norm and \eqref{eq:prep-unit-collar}, we may choose $\tilde \theta\in(\theta_*,\delta_c)$ sufficiently close to $\theta_*$ that
 \[
    q_{X_0}\ge\kappa_0
    \quad\text{on }\mathcal C_{\tilde\theta}.
\]
Moreover, monotonicity of  $a$ and \eqref{eq:prep-norm} give
 \[
    \Delta
    \ge\delta_0:=\min\left\{
        1-a(\tilde \theta)^2,\,
        1-\eta^2M^2
    \right\}>0
    \quad\text{on }N\setminus\mathcal C_{\tilde\theta},
\]
where the second term follows from \eqref{eq:prep-filled} and \eqref{eq:prep-rescaled}. Since $\tau_0(X_0)=|q_{X_0}|$, replacing $\kappa_0$ by $\min\{\kappa_0,\delta_0\}$ gives $\Delta+\tau_0(X_0)\geq\kappa_0$ on $N$. This proves \ref{prop:it3}.\\
For \ref{prop:it4}, it is enough to work on a neighborhood of
$\partial_N\mathcal C=\{\tau=0\}$. For $0<\tau<\eps_h$, the definition of $h$ yields
\[
    a=1-(1-\eta)\tau^2e^{-1/\tau},
    \qquad
    \theta a'=-\frac{1-\eta}{L}(1+2\tau)e^{-1/\tau}.
\]
Since $a\leq1$, \eqref{eq:prep-norm} implies together with $a\ge \eta$ that
\begin{equation}
    \label{eq:LowerBoundDelta}
    \Delta\geq1-a^2-a\theta a'
    \geq-a\theta a' \ge\frac{\eta(1-\eta)}{L}e^{-1/\tau}.
\end{equation}
On the other hand, substituting the expressions above into
\[
    \Delta=1-a^2-2a\theta a'(X_c\cdot\Xi)-(\theta a')^2|\Xi|^2
\]
and using the weighted derivative bounds for $\tau$, $X_c$ and $\Xi$ gives
\[
    \sigma^k|D^k\Delta|\leq C_k e^{-1/\tau}(1+\tau^{-2k}),\qquad k\geq0.
\]
The chain rule and \eqref{eq:LowerBoundDelta} hence yield
\[
    \sigma^k|D^k\sqrt{\Delta}|\leq C_k e^{-1/(2\tau)}(1+\tau^{-2k})\leq C_{k,\eps}(\sqrt{\Delta})^{1-\eps},\qquad k\geq0,\quad 0<\eps<1,
\]
where we absorbed the powers of $\tau^{-1}$ into the exponential. In particular, this shows that $\sqrt{\Delta}$ extends smoothly by zero across $\partial_N\mathcal C$. Away from a fixed neighborhood of $\partial_N\mathcal C$, $\Delta$ is either identically zero or positive so smoothness and the derivative bounds are clear. This proves \ref{prop:it4}.\\
Finally, for \ref{prop:it5}, the functions  $a$ and $\theta a'$ are smooth functions of $\theta$ on
$[0,\delta_c]$. The derivative bounds for $\theta$ and the chain rule therefore give
 \[
    ||a||_{k;\mathcal C_{\delta_c}}
    +||\theta a'||_{k;\mathcal C_{\delta_c}}<\infty,
    \qquad k\ge0.
\]
The product rule, \eqref{eq:prep-rescaled}, and the bounds for  $X_c$, $\Xi$ and  $X_f$ now yield
\[
    ||X_0||_k<\infty,\qquad k\ge0.
\]
This completes the proof.
\end{proof}

\section{Local Perturbations}
\label{sec:LocPert}

In the following we work  on $B_2(0)\subset\mathbb R^3$. Fix a unit vector $\xi\in\mathbb R^3$ and denote $\varphi(x)=x\cdot\xi$. Let $V\in C^\infty(B_2(0);\mathbb R^3)$ be nowhere vanishing and $\rho\in C_c^\infty(B_2(0))$ be nonnegative. Suppose that for some $\mu\geq2$, we have for every $k\geq0$ and $0<\eps<1$
 \begin{equation}
    \label{eq:AmplitudeBounds}
    |D^k\rho|\leq C_{k,\eps}\mu^k\rho^{1-\eps},\quad\text{on }B_2(0).
\end{equation}
and for each $k\geq0$
\begin{equation}
    \label{eq:DataMuBound}
    ||V||_{C^k},||\rho||_{C^k}\leq C_k\mu^k.
\end{equation}
Finally we suppose
\begin{equation}
    \label{eq:Transversality}
    |\Pi_{\xi^\perp}V|>c>0,\quad\text{in a neighborhood of }\supp\rho.
\end{equation}
In the application of the perturbations later, $B_2(0)$ is the domain of a chart $\Phi_i(y)=x_i+\ell_i y$ from Lemma \ref{lem:BoundCharts}. We use the normalizations \eqref{eq:ChartNormalizations}, namely
\begin{equation}
    \label{eq:RescalingSec4}
    X\mapsto X\circ\Phi_i,\qquad B\mapsto\ell_i^{-1}B\circ\Phi_i,\qquad \widehat\sigma=\frac{\sigma\circ\Phi_i}{\ell_i}.
\end{equation}
We do not relabel these rescaled fields. In these coordinates, the quantities from Definition \ref{def:Nondegen} are
\[
q_X=\widehat\sigma X\cdot\curl X,\qquad
H_\alpha^X=\widehat\sigma\Pi_{X^\perp}\mathbf e_\alpha,\qquad
\tau_j(X)^2=\sum_{|I|\leq j}|H_I^Xq_X|^2.
\]
Thus $q_X$ and $\tau_j(X)$ are the original quantities composed with $\Phi_i$. All $C^k$ norms below are taken on $B_2(0)$ and correspond to the weighted norms by \eqref{eq:Normequivalence}.
\begingroup
\begin{proposition}[Step]
    \label{prop:LocalPert}
    There exists $\gamma>0$ such that, if
    \begin{equation}
        \label{eq:SmallnessLocPert}
        \sup_{x\in\supp\rho}
        \frac{\rho(x)}{|\Pi_{\xi^\perp}V(x)|}\leq\gamma,
    \end{equation}
    then the following holds. Fix $m,L\in\mathbb N_0$, $\beta\geq8(m+2)$ and $\kappa,\rho_*>0$. For sufficiently large $\lambda=\mu^\beta$, there exists
    \[
        B_\lambda\in C_c^\infty(B_2(0);\mathbb R^3),
        \qquad \supp B_\lambda\subset\supp\rho,
    \]
    such that $X_\lambda=V+\curl B_\lambda$ is nowhere vanishing and we have the following:
    \begin{enumerate}
        \item The length error satisfies
        \begin{equation}
            \label{eq:LocalPertError}
            \big\||X_\lambda|^2-|V|^2-\rho^2\big\|_{C^m}
            \leq C_{L,m}\lambda^m
            \left(\frac{\mu}{\lambda}\right)^{L+1}.
        \end{equation}
        \item We have
        \begin{equation}
            \label{eq:LocalPertC0}
            ||X_\lambda-V||_{C^0}
            \leq C||\rho||_{C^0}+C_L\frac{\mu}{\lambda},
        \end{equation}
        and, for every $k\geq0$,
        \begin{equation}
            \label{eq:LocalPertDerivatives}
            ||X_\lambda||_{C^k}\leq C_{L,k}\lambda^k.
        \end{equation}
        \item For every $0\leq j\leq m$,
        \begin{equation}
            \label{eq:LocalPertPreservation}
            \tau_j(V)(x)\geq\kappa\mu^{-j}
            \quad\Longrightarrow\quad
            \tau_{j+2}(X_\lambda)(x)\geq c_1\lambda^{-(j+2)}.
        \end{equation}
        Moreover,
        \begin{equation}
            \label{eq:LocalPertCreation}
            \rho(x)\geq\rho_*
            \quad\Longrightarrow\quad
            \tau_1(X_\lambda)(x)\geq c_*>0.
        \end{equation}
    \end{enumerate}
    The constant $\gamma$ and the constant $C$ in \eqref{eq:LocalPertC0} are universal. The other constants may depend on $c,m,\beta,L,\kappa,\rho_*$ and the bounds \eqref{eq:AmplitudeBounds}--\eqref{eq:DataMuBound}, but are independent of $\mu$ and of the chart.
\end{proposition}
\endgroup
We split the proof of the proposition into two lemmas. In Lemma \ref{lem:LocalPert} we construct the perturbation and derive the error estimates, in Lemma \ref{lem:FiniteTypeAfterPer} we compute the change in the nondegeneracy condition.
\begin{remark}
    \label{rem:MotivationPerturbation} 
    Recall the standard building block $Y$ from the introduction. Let $w=Y-V$. If we perturb as in the introduction by 
    \[
    X_\lambda=V+\curl\left(\frac{1}{\lambda}\mc P_\xi w\right),
    \]
    we obtain
    \begin{equation*}
        |X_\lambda|^2=|V|^2+\rho^2+E,
    \end{equation*}
    where $E=O(\lambda^{-1})$. In Lemma \ref{lem:LocalPert} we first make the construction of $Y$ precise. Afterwards we obtain through successive corrections to $w$ an error $E$ as above which is of higher order in $\lambda^{-1}$. Let $w_j$ be the perturbation and $E_j$ be the corresponding length error after $j$ corrections. Writing $\widetilde E_j$ for the coefficient of $\lambda^{-j-1}$ in $E_j$, the idea behind the corrections $h_j$ is to pick them such that
     \begin{equation}
        \label{eq:Propertiesofhj}
        2Y\cdot h_j =\widetilde E_j,\quad \overline h_j=0,\quad h_j\cdot \xi =0.
    \end{equation}
    The first property is exactly what cancels the leading term of the error $E_j$ at each step if we replace the perturbation $w_j$ by $w_j-\lambda^{-j-1}h_j$. For a related construction using successive corrections to cancel the leading error terms, see Cao, Hirsch and Inauen \cite[Section~2.2]{CaoHirschInauen2025}.
    \begin{figure}[htbp]
\centering
%
 \resizebox{0.48\textwidth}{!}{%
\begin{tikzpicture}[
    font=\large,
    >=Latex,
    line cap=round,
    line join=round
]
    \path[use as bounding box]
        (-0.1,-4.6) rectangle (9.9,4.15);

    \coordinate (O) at (1.8,0);
    \coordinate (V) at (5.35440687,0);
    \coordinate (Y) at
        ({1.8+4.2*cos(27)},{4.2*sin(27)});

    \draw[figuregray!65,line width=0.65pt]
        ({1.8+4.2*cos(-60)},{4.2*sin(-60)})
        arc[start angle=-60,end angle=60,radius=4.2];

    \draw[profileblue,line width=2.2pt]
        ({1.8+4.2*cos(-45.836624)},
         {4.2*sin(-45.836624)})
        arc[
            start angle=-45.836624,
            end angle=45.836624,
            radius=4.2
        ];

    \draw[
        profileblue!70,
        line width=0.9pt,
        -{Latex[length=2.3mm]}
    ]
        (O) -- (Y)
        node[
            pos=0.47,above,sloped,
            text=black,inner sep=5pt
        ] {$\sqrt{b^2+\rho^2}$};

    \draw[
        black!75,
        line width=0.9pt,
        -{Latex[length=2.3mm]}
    ]
        (O) -- (V);

    \fill[black!75] (O) circle[radius=1.5pt];
    \node[below left,inner sep=5pt] at (O)
        {$(V\cdot\xi)\xi$};

    \fill[black] (V) circle[radius=2.1pt];
    \node[below left,inner sep=7pt] at (V)
        {$V$};

    \fill[profileblue] (Y) circle[radius=2.4pt];
    \node[right,inner sep=8pt,text=profileblue] at (Y)
        {$Y(x,\vartheta)$};
\end{tikzpicture}%
}%
\hfill
%
\resizebox{0.48\textwidth}{!}{%
\begin{tikzpicture}[
    font=\large,
    >=Latex,
    line cap=round,
    line join=round
]
    \path[use as bounding box]
        (-0.1,-4.6) rectangle (9.9,4.15);

    \draw[
        figuregray,
        line width=0.6pt,
        -{Latex[length=2mm]}
    ]
        (0.8,-1.1) -- (0.8,2.85);

    \draw[
        figuregray,
        line width=0.6pt,
        -{Latex[length=2mm]}
    ]
        (0.8,0) -- (9.2,0)
        node[right,text=black] {$\vartheta$};

    \node[below left,inner sep=4pt] at (0.8,0) {$0$};
    \node[below,inner sep=6pt] at (8.8,0) {$2\pi$};

    \draw[errororange,line width=1.2pt]
        plot[
            domain=0:360,
            samples=181,
            smooth,
            variable=\t
        ]
        ({0.8+8*\t/360},
         {1.3*(0.45+1.12*sin(\t)-0.22*cos(2*\t))});

    \draw[profileblue,line width=1.4pt]
        plot[
            domain=0:360,
            samples=181,
            smooth,
            variable=\t
        ]
        ({0.8+8*\t/360},
         {0.11+0.20*sin(2*\t)-0.08*cos(\t)});

    \node[text=errororange] at (2.8,2.7) {$E_j$};
    \node[
        text=profileblue,
        fill=white,
        inner sep=3pt
    ] at (5.5,0.75) {$E_{j+1}$};

    \node at (4.9,-2.2)
        {$w_{j+1}=w_j-\lambda^{-j-1}h_j,\qquad O(\lambda^{-j-1})
          \;\longrightarrow\;
          O(\lambda^{-j-2})$};
\end{tikzpicture}%
}

\caption{Leading profile (left) and error (right).}
\label{fig:LocalPert}
\end{figure}
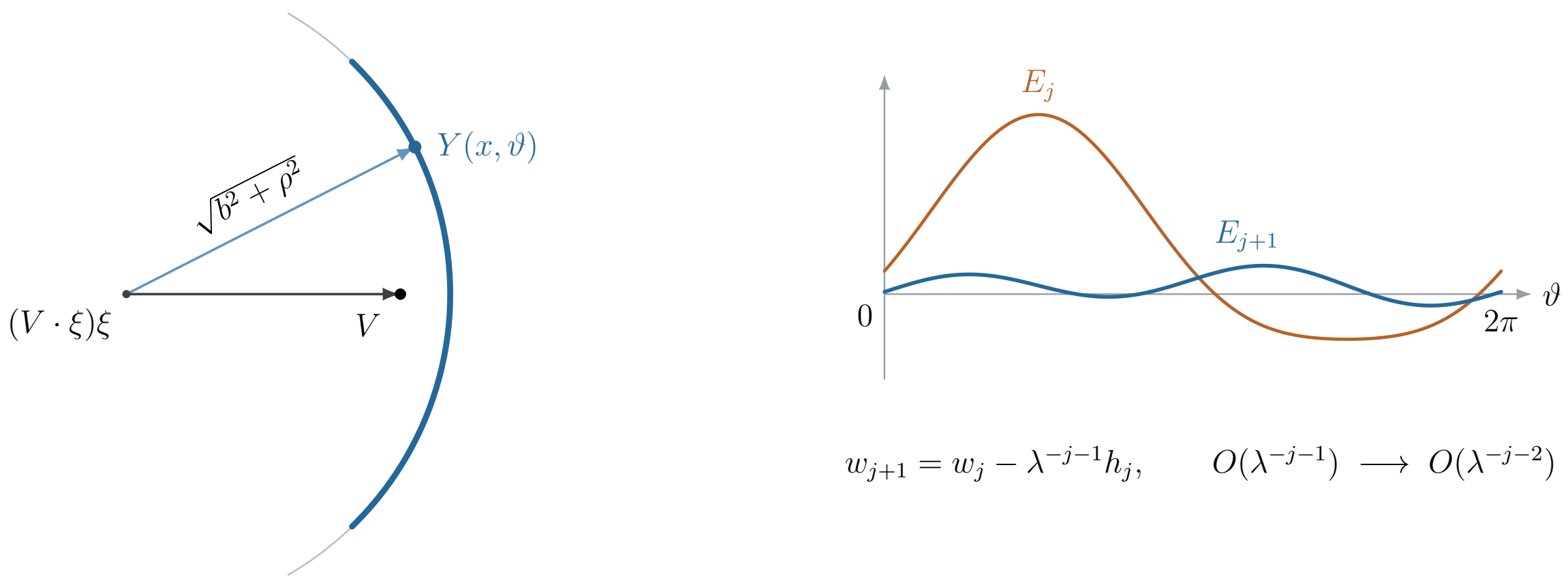
\end{remark}

\begin{lemma}
    \label{lem:LocalPert} 
    There exists $\gamma>0$ such that, under \eqref{eq:SmallnessLocPert}, the following holds. For every $L \in \mb N_0$ and sufficiently large $\lambda\geq 2\mu$, there exists  $B_\lambda\in C^\infty(B_2(0);\mb R^3)$ with $\supp B_\lambda \subset \supp {\rho}$ such that
    \[
    X_\lambda=V+\curl B_\lambda,\quad |X_\lambda|^2=|V|^2+\rho^2+E_\lambda.
    \]
    For every $k \geq 0$ we have
     \[
    ||E_\lambda||_{C^k}\leq C_{L,k}\lambda^k \left(\frac{\mu}{\lambda}\right)^{L+1},\quad ||X_\lambda||_{C^k}\leq C_{L,k}\lambda^k.
    \]
    We can describe the perturbation more precisely. There are smooth  $Y(x,\vartheta),R_\lambda(x,\vartheta)\in \mb R^3$, which are $2\pi$-periodic in $\vartheta$ such that
     \[
    X_\lambda(x)=Y(x,\lambda\varphi(x))+R_\lambda(x,\lambda\varphi(x)),\quad |Y-V|\leq C\rho.
    \]
    For every $0<\eps<1$ and $k,l\ge 0$ we have
     \[
    \begin{aligned}
        |D_x^k\p_\vartheta^l (Y-V)|\leq C_{k,l,\eps}\mu^k {\rho}^{1-\eps},\\
        |D_x^k\p_\vartheta^l R_\lambda|\leq C_{L,k,l,\eps}\mu^k (\mu/\lambda){\rho}^{1-\eps}.
    \end{aligned}
    \]
\end{lemma}

\begin{proof}
    Let $f(x,\vartheta)$ be $2\pi$-periodic in $\vartheta$. We write
    \[
    \overline f(x)=\frac{1}{2\pi}\int_0^{2\pi}f(x,\vartheta)\,d\vartheta.
    \]
    \textbf{Step 1: Leading profile.}  We use the following formulas near $\supp{\rho}$, where $b=|\Pi_{\xi^\perp}V|>c$, and set $a=0$ and $Y=V$ elsewhere. Denote
     \[
    e_1=\frac{V-(V\cdot\xi)\xi}{b}.
    \]
    Define for some smooth function  $a \in C^\infty(B_2(0))$ which we choose below, the profile
     \begin{equation}
        \label{eq:AnsatzY}
        Y=V-be_1 +\sqrt{b^2+\rho^2}\big(\cos(a\sin \vartheta)e_1+\sin(a\sin \vartheta)(\xi \times e_1) \big).
    \end{equation}
    Following Remark \ref{rem:MotivationPerturbation}, we want the function $a$ to satisfy
    \begin{equation}
        \label{eq:avga}
        \overline{\cos (a\sin \vartheta)}=\frac{b}{\sqrt{b^2+\rho^2}}.
    \end{equation}
     To solve \eqref{eq:avga}, note that
    \[
        \overline{\cos(a\sin\vartheta)}
        =1-\frac{a^2}{4}+O(a^4),
        \qquad
        \frac{b}{\sqrt{b^2+\rho^2}}
        =1-\frac12(\rho/b)^2+O((\rho/b)^4).
    \]
    In particular, the left-hand side of \eqref{eq:avga} is analytic in $a^2$, with derivative $-1/4$ at zero. Applying the implicit function theorem with $a^2$ as the unknown hence gives
    \[
        a^2=(\rho/b)^2\bigl(2+O((\rho/b)^2)\bigr),
    \]
    where the factor in parentheses is smooth in $(\rho/b)^2$ and positive near zero. Taking its square root and multiplying by $\rho/b$ yields a smooth
    \begin{equation}
        \label{eq:Choicea}
        a=\sqrt{2}\,\frac{\rho}{b}
          +O((\rho/b)^3),
    \end{equation}
    such that \eqref{eq:avga} holds. Pick $\gamma>0$ sufficiently small  such that the construction above is defined for $0\leq \rho/b\leq\gamma$ and we have $0\leq a\leq 1/2$. With this choice we obtain from \eqref{eq:avga} and the Ansatz for $Y$ that
     \begin{equation}
        \label{eq:propertiesY}
        \overline Y=V,\quad |Y|^2=|V|^2+\rho^2,\quad (Y-V)\cdot \xi=0.
    \end{equation}
    Moreover, from \eqref{eq:Choicea} it is easy to see that
     \begin{equation}
        \label{eq:boundDiffYV}
        |Y-V|\leq C\rho,\quad 2Y\cdot e_1\geq b >0.
    \end{equation}
    The constants $\gamma$ and $C$ are universal.\\
    Recall $\mc P:=\mc P_\xi$ from the introduction. Suppose  $w(x,\vartheta) \in \mb R^3$ is smooth, $2\pi$-periodic in $\vartheta$ and satisfies
     \begin{equation}
        \label{eq:Propw}
        w(x,\vartheta)\cdot \xi =0,\quad \overline{w}=0.
    \end{equation}
    Denote
     \[
    \mc P w=-\xi\times\partial^{-1}_\vartheta w,
    \]
    where $\partial^{-1}_\vartheta$ is the unique average $0$ primitive in $\vartheta$. Note that from \eqref{eq:Propw} and the chain rule one has
     \begin{equation}
        \label{eq:PerturbationThroughCurl}
        \curl(\lambda^{-1}(\mc P w)(x,\lambda\varphi(x)))=(w+\lambda^{-1}\curl_x (\mc P w))(x,\lambda\varphi(x)),
    \end{equation}
    where $\curl_x$ acts only on the spatial variables with $\vartheta$ fixed.\\
    \textbf{Step 2: Induction.} Let $w_0=Y-V$. We construct inductively $w_j$ for $0\leq j\leq L$ with the following properties. If we write
    \[
    X_j=V+w_j+\lambda^{-1}\curl_x(\mc P w_j),\quad E_j=|X_j|^2-|V|^2-\rho^2,
    \]
    we claim that
     \begin{equation}
        \label{eq:Indwjbas}
        w_j=Y-V+O(\lambda^{-1}),\quad w_j\cdot \xi =0,\quad \overline{w_j}=0,\quad E_j=O(\lambda^{-j-1}).
    \end{equation}
    Moreover, we construct $w_j$ as a polynomial of degree $j$ in $\lambda^{-1}$,
    \begin{equation}
        \label{eq:Expansionwj}
        w_j=\sum_{r=0}^j\lambda^{-r}f_r,\qquad f_0=Y-V.
    \end{equation}
    We claim for $0<\eps<1$ 
     \begin{equation}
        \label{eq:IndwjEst}
        |D_x^k\p_\vartheta^l f_r|\leq C_{j,r,k,l,\eps}\mu^{r+k}{\rho}^{1-\eps},\quad \text{on }\{{\rho}>0\}.
    \end{equation}
    Note that the formula for $X_j$ together with the Leibniz rule give the same bounds for the coefficients of $X_j-V$ and $E_j$ after expanding in $\lambda^{-1}$.\\
    For $0<\eps'<1$, we use on $\{{\rho}>0\}$ the following bounds
     \begin{equation}
        \label{eq:CutoffBounds}
        |D^k \rho|\leq C_{k,\eps'}\mu^k {\rho}^{1-\eps'},\quad |D^k (\rho^{-1})|\leq C_{k,\eps'}\mu^k {\rho}^{-1-\eps'},
    \end{equation}
    where the first bound is \eqref{eq:AmplitudeBounds}, and the second follows by a simple computation from \eqref{eq:AmplitudeBounds}.\\
    Let $j=0$. Since $w_0=Y-V$, we can use the definition of $Y$, together with the choice of $a$ in \eqref{eq:Choicea}, to see that $w_0$ factors as $\rho$ times a smooth function, which has mixed derivatives of order $k,l$ bounded by $C_{k,l} \mu^k$ due to \eqref{eq:DataMuBound}. Hence, together with $\eqref{eq:CutoffBounds}$ we obtain \eqref{eq:IndwjEst}. The statements \eqref{eq:Indwjbas} follow directly from \eqref{eq:propertiesY}.\\
    Next suppose that the induction hypotheses hold for $j<L$. Let $\widetilde E_j$ be the coefficient of $\lambda^{-j-1}$ in $E_j$, so that
    \[
    E_j=\lambda^{-j-1}\widetilde E_j+O(\lambda^{-j-2}).
    \]
    We construct the next perturbation by cancelling this coefficient. To do so, differentiate \eqref{eq:propertiesY} in the amplitude parameter $\rho$. This gives
     \begin{equation}
        \label{eq:DiffPropY}
        2Y\cdot \p_\rho Y=2\rho,\quad \overline{\p_\rho Y}=0,\quad \p_\rho Y\cdot \xi =0.
    \end{equation}
    Define on $\{\rho>0\}$
     \begin{equation}
        \label{eq:choicehj}
        c_j=\frac{\overline{\widetilde E_j/(2Y\cdot e_1)}}{\overline{(2Y\cdot e_1)^{-1}}},\quad h_j=\frac{c_j}{2\rho}\p_\rho Y+\frac{\widetilde E_j-c_j}{2Y\cdot e_1}e_1.
    \end{equation}
    From \eqref{eq:DiffPropY}, we have
     \begin{equation}
        \label{eq:Prophj}
        h_j\cdot \xi =0,\quad \overline{h_j}=0,\quad 2Y\cdot h_j =\widetilde E_j.
    \end{equation}
    Define $w_{j+1}=w_j-\lambda^{-j-1}h_j$, so that $f_{j+1}=-h_j$. We first verify \eqref{eq:IndwjEst} for $h_j$. Since $\overline{X_j-V}=0$, it follows that
    \[
    \overline{E_j}=\overline{|X_j-V|^2}-\rho^2.
    \]
    Using the expansion \eqref{eq:Expansionwj}, we write
    \[
    X_j-V=f_0+\sum_{r=1}^j\lambda^{-r}
    \big(f_r+\curl_x(\mc P f_{r-1})\big)
    +\lambda^{-j-1}\curl_x(\mc P f_j).
    \]
    Taking the coefficient of $\lambda^{-j-1}$ in the identity for $\overline{E_j}$, since there is no contribution from $-\rho^2$ we obtain
    \[
    \overline{\widetilde E_j}
    =2\overline{f_0\cdot\curl_x(\mc P f_j)} +\sum_{r=1}^j\overline{\big(f_r+\curl_x(\mc P f_{r-1})\big)\cdot
    \big(f_{j+1-r}+\curl_x(\mc P f_{j-r})\big)}.
    \]
    In each product above, the powers of $\mu$ add to $j+1$, and each factor has the weight $\rho^{1-\eps/2}$ by \eqref{eq:IndwjEst}. Thus Leibniz rule gives
    \begin{equation}
        \label{eq:AveragedCoefficientEst}
        |D^k\overline{\widetilde E_j}|
        \leq C_{j,k,\eps}\mu^{j+1+k}\rho^{2-\eps}.
    \end{equation}
    Since $b$ is independent of $\vartheta$ and $V\cdot e_1=b$, we have
    \begin{align}
        \label{eq:CoefficientPart}
        \overline{\frac{\widetilde E_j}{2Y\cdot e_1}}
        &=\overline{\widetilde E_j
        \left(\frac{1}{2Y\cdot e_1}-\frac{1}{2b}\right)}
        +\frac{\overline{\widetilde E_j}}{2b} \nonumber \\
        &=-\overline{\frac{\widetilde E_j (Y-V)\cdot e_1}
        {2b(Y\cdot e_1)}}+\frac{\overline{\widetilde E_j}}{2b}.
    \end{align}
    The coefficient estimate for $\widetilde E_j$ and the factorization of $Y-V$ as $\rho$ times a smooth function give the same derivative bound for the first term of \eqref{eq:CoefficientPart}. The second term is controlled by \eqref{eq:AveragedCoefficientEst}, as $b>c$. The denominator in the definition of $c_j$ is bounded away from $0$ and has derivatives of order $k$ bounded by $C_k\mu^k$. Using Leibniz rule once again and \eqref{eq:CutoffBounds}, we obtain from the definition of $c_j$ that
     \[
    \left|D^k \left(\frac{c_j}{\rho}\right)\right|\leq C_{j,k,\eps}\mu^{j+1+k}{\rho}^{1-\eps}.
    \]
    The derivatives of $\p_\rho Y$, $e_1$ and $(2Y\cdot e_1)^{-1} $ are due to \eqref{eq:DataMuBound} bounded by $C_{k,l}\mu^k$. From \eqref{eq:choicehj} we therefore see that $h_j$ satisfies \eqref{eq:IndwjEst} with $r=j+1$, which proves the required bound for $f_{j+1}=-h_j$. We extend $h_j$ by zero on $\{\rho=0\}$. The preceding estimates show that this extension is smooth. Thus $w_{j+1}=w_j-\lambda^{-j-1}h_j$ is smooth on $B_2(0)$ and supported in $\supp\rho$.\\
    It remains to check \eqref{eq:Indwjbas}. Since $h_j$ is independent of $\lambda$, the correction $\lambda^{-j-1}h_j$ is $O(\lambda^{-j-1})$. In particular, the constant coefficient of $w_{j+1}$ is $Y-V$, which proves
    \[
    w_{j+1}=Y-V+O(\lambda^{-1}).
    \]
    The statements  $w_{j+1}\cdot \xi=0 $ and $\overline{w_{j+1}}=0$ follow immediately from \eqref{eq:Prophj}. Finally, we use \eqref{eq:Prophj} to expand
     \begin{align*}
        E_{j+1}&=E_j-\lambda^{-j-1}\widetilde E_j-2\lambda^{-j-1}(X_j-Y)\cdot h_j\\
        &\quad-2\lambda^{-j-2}X_j\cdot\curl_x(\mc P h_j)\\
        &\quad+\lambda^{-2j-2}|h_j+\lambda^{-1}\curl_x(\mc P h_j)|^2.
    \end{align*}
    By the choice of $\widetilde E_j$, the first difference is $O(\lambda^{-j-2})$. The next two terms are also $O(\lambda^{-j-2})$, and the last term is $O(\lambda^{-2j-2})=O(\lambda^{-j-2})$. This proves $E_{j+1}=O(\lambda^{-j-2})$ and completes the proof of the induction.\\
    \textbf{Step 3: Conclusion.} We apply the induction up until $j=L$. Denote
     \[
    B_\lambda=\lambda^{-1}(\mc P w_L)(x,\lambda\varphi(x)),\quad E_\lambda(x)=E_L(x,\lambda\varphi(x)).
    \]
    Set  $X_\lambda=V+\curl B_\lambda$. Applying \eqref{eq:PerturbationThroughCurl} with $w=w_L$ gives
     \[
    X_\lambda(x)=V(x)+\curl B_\lambda(x)
    =X_L(x,\lambda\varphi(x)).
    \]
    Hence, by the definition of $E_L$,
     \[
    |X_\lambda|^2=|V|^2+\rho^2+E_\lambda.
    \]
    As $\mc P$ preserves the support in $x$, we have  $\supp B_\lambda\subset\supp{\rho}$.
    The polynomial $E_L$ has degree at most $2L+2$ and contains only powers $\lambda^{-r}$ with $r\geq L+1$. Since its coefficients satisfy \eqref{eq:IndwjEst}, the chain rule gives
    \[
    ||E_\lambda||_{C^k}\leq C\lambda^k\sum_{r=L+1}^{2L+2}\left(\frac{\mu}{\lambda}\right)^r
    \leq C\lambda^k\left(\frac{\mu}{\lambda}\right)^{L+1}.
    \]
    This proves the error estimate. Applying the same argument to $X_L-V$ and adding the bounds for $V$ gives the estimate for  $X_\lambda$. The bound for $V-Y$ was already shown in \eqref{eq:boundDiffYV}. Finally, let $R_\lambda(x,\vartheta)=X_L(x,\vartheta)-Y(x,\vartheta)$. This contains by \eqref{eq:Indwjbas} only positive powers of $\lambda^{-1}$ and \eqref{eq:IndwjEst} gives immediately the derivative estimates. The estimates for $Y-V$ follow straight from the definition of $Y$ together with \eqref{eq:CutoffBounds}.
    Taking $k=l=0$ and $\eps=1/2$ in the remainder estimate and using \eqref{eq:DataMuBound}, we obtain
    \[
    ||X_\lambda-V||_{C^0}
    \leq C||\rho||_{C^0}+C_L\frac{\mu}{\lambda}.
    \]
    Moreover, $|Y|^2=|V|^2+\rho^2\geq c^2$ near $\supp\rho$. Thus $X_\lambda$ is nowhere vanishing if $\mu/\lambda$ is sufficiently small. This completes the proof.
\end{proof}

\begin{remark}
    Let us give a brief intuition why the perturbations of the previous lemma can be used to create the nondegeneracy condition if we increase the order of the condition by two. For this one has to understand how the perturbations change $q_{X_\lambda}$. Below we compute this in terms of the profile
     \begin{equation}
        \label{qXlamProf}
        Q=\frac{q_V}{\lambda}-Ga\cos\vartheta+E,
    \end{equation}
    where  $E$ is an error term. From \eqref{qXlamProf}, we see the key idea: While  $-Ga\cos\vartheta$ may vanish at some point, its derivative is  $Ga\sin\vartheta$. Because sine and cosine can not vanish simultaneously we should always be able to detect nondegeneracy, as long as we look at the term and its derivative together.
\end{remark}
We use the notation of Lemma \ref{lem:LocalPert}.
\begin{lemma}
    \label{lem:FiniteTypeAfterPer}
    Fix $m \in \mb N_0$,  $\beta \geq 8(m+2)$ and $L \in \mb N_0$. Let  $\lambda=\mu^\beta$. For sufficiently large $\lambda$, the field $X_\lambda$ constructed in Lemma \ref{lem:LocalPert} has the following properties.
    \begin{enumerate}
        \item \label{it1:type} If $\tau_m(V)(x)\geq c_0\mu^{-m}$, then
        \[
        \tau_{m+2}(X_\lambda)(x)\geq c_1\lambda^{-(m+2)}.
        \]
        \item \label{it2:type} Let $\rho_*>0$. If $\rho(x)\geq\rho_*$ then
        \[
        \tau_1(X_\lambda)(x)\geq c_*>0.
        \]
    \end{enumerate}
    The constants may depend on $m,\beta,L$ and the bounds \eqref{eq:DataMuBound}, \eqref{eq:AmplitudeBounds}, as well as on $c_0$ in \ref{it1:type} and on $\rho_*$ in \ref{it2:type}, but are independent of the chart $\Phi_i$.
\end{lemma}

\begin{proof}
    Recall that the perturbation vanishes outside $\supp {\rho}$. Hence it is enough to work near $\supp {\rho}$, otherwise there is nothing to prove.\\
    \textbf{Step 1: Helicity expansion.} From the proof of Lemma \ref{lem:LocalPert} we have the formula
     \[
    Y=V-be_1+\sqrt{b^2+\rho^2}\big(\cos(a\sin \vartheta)e_1+\sin(a\sin \vartheta)(\xi \times e_1)\big).
    \]
    Also from \eqref{eq:Choicea}
    \begin{equation}
        \label{eq:acomprhochi}
        a\simeq \rho/b \simeq {\rho}
    \end{equation}
    so we can assume that $0\leq a\leq a_0<1$. Moreover, note by using the cutoff bounds \eqref{eq:CutoffBounds} it holds for any $0<\eps<1$
     \begin{equation}
        \label{eq:BoundsDera}
        |D^k a|\leq C\mu^k a^{1-\eps}.
    \end{equation}
    Denoting
    \[
    X(x,\vartheta)=Y(x,\vartheta)+R_\lambda(x,\vartheta),\quad X_\lambda(x)=X(x,\lambda\varphi(x)),
    \]
    we have from Lemma \ref{lem:LocalPert} together with \eqref{eq:acomprhochi} the bounds
     \begin{equation}
        \label{eq:BoundsRandVY}
        \begin{aligned}
            |D_x^k\p_\vartheta^l(Y-V)|\leq C\mu^k a^{1-\eps}\\
            |D_x^k\p_\vartheta^l R_\lambda|\leq C\mu^k \frac{\mu}{\lambda}a^{1-\eps}.
        \end{aligned}
    \end{equation}
    We define the function
     \[
    Q(x,\vartheta)=\widehat\sigma X\cdot (\xi\times \p_\vartheta X+\lambda^{-1}\curl_x X),
    \]
    so that one computes
    \[
    Q(x,\lambda\varphi(x))=\frac{q_{X_\lambda}}{\lambda}.
    \]
    One can compute directly from the definition that
     \[
    Y\cdot(\xi\times\p_\vartheta Y) =-(b^2+\rho^2)a \cos \vartheta.
    \]
    Hence, if we set  $G=\widehat\sigma(b^2+\rho^2)\geq C^{-1}>0$, we expand
     \begin{equation}
        \label{eq:ExpansionForQ}
        Q=\frac{q_V}{\lambda}-Ga\cos \vartheta+E.
    \end{equation}
    We claim that  $E$ satisfies
    \begin{equation}
        \label{eq:EstForB}
        |D_x^k \p_\vartheta^l E|\leq C \mu^k\frac{\mu}{\lambda}a^{1-\eps}.
    \end{equation}
    Indeed, one has
    \[
    E=\widehat\sigma\left(
        R_\lambda\cdot(\xi\times\p_\vartheta X)+Y\cdot(\xi\times\p_\vartheta R_\lambda)
        +\lambda^{-1}(X-V)\cdot\curl_x X
        +\lambda^{-1}V\cdot\curl_x(X-V)
    \right).
    \]
    Using $X-V=(Y-V)+R_\lambda$, the product rule, \eqref{eq:DataMuBound}, \eqref{eq:BoundsRandVY} and the bounds for $\widehat\sigma$ we obtain \eqref{eq:EstForB}.\\
    \textbf{Step 2: Horizontal Vector Field.} Next we introduce the horizontal vector field which we use to detect the changes in $\tau_m$. We define
     \[
    W(x,\vartheta)=\frac{\Pi_{X(x,\vartheta)^\perp}\xi}{|\Pi_{X(x,\vartheta)^\perp}\xi|^2},\quad W_\lambda(x)=W(x,\lambda\varphi(x)).
    \]
    This is well-defined since the denominator is positive for large enough $\mu$. Clearly we have  $W\cdot X=0$. Moreover,
     \begin{equation}
        \label{eq:EstW}
        \begin{aligned}
            &|D_x^k\p_\vartheta^l W|\leq C\mu^k\quad k,l\geq 0\\
            &|D_x^k\p_\vartheta^l W|\leq C\mu^k a^{1-\eps}\quad k\geq0,l\geq 1.
        \end{aligned}
    \end{equation}
    The factor $a^{1-\eps}$ in the second estimate appears since every derivative in $\vartheta$ produces a $\vartheta$-derivative of $X$. Since $V$ is independent of $\vartheta$ we can write $\p_\vartheta X=\p_\vartheta (Y-V)+\p_\vartheta R_\lambda$ and the bound follows from \eqref{eq:BoundsRandVY}. Denote the operator
     \[
    \mc L =\p_\vartheta +\lambda^{-1}W\cdot\nabla_x.
    \]
    From the chain rule and the fact that $W\varphi =1$, one has
    \begin{equation}
        \label{eq:RepeatedLdiff}
        (\mc L^r Q)(x,\lambda\varphi(x))=\lambda^{-r-1}W_\lambda^r q_{X_\lambda}(x).
    \end{equation}
     Since $W_\lambda\perp X_\lambda$, we have
    \[
    W_\lambda=\sum_{\alpha=1}^3 c_\alpha H_\alpha^{X_\lambda},\qquad c_\alpha=\widehat\sigma^{-1}W_\lambda\cdot \mathbf e_\alpha,
    \]
    where the coefficients $c_\alpha$ have $j$-th derivatives of order $O(\lambda^j)$, since they depend on $X_\lambda$. Hence, we obtain by recalling that $W_\lambda$ is horizontal and combining \eqref{eq:RepeatedLdiff} with the product rule
     \begin{equation}
        \label{ineq:Lrtaur}
        |(\mc L^r Q)(x,\lambda\varphi(x))|
        \leq C_r\sum_{|I|\leq r}\lambda^{-|I|-1}|H_I^{X_\lambda}q_{X_\lambda}(x)|\leq C_r\tau_r(X_\lambda)(x).
    \end{equation}
    Moreover, we claim that for every $r\geq 0$ it holds
     \begin{equation}
        \label{eq:MainFiniteEst}
        |\mc L^r Q+Ga\p_\vartheta^r \cos \vartheta|\leq C\left(\frac{\mu}{\lambda}a^{1-\eps}+\left(\frac{\mu}{\lambda}\right)^{r+1}\right).
    \end{equation}
    Indeed, recall the splitting \eqref{eq:ExpansionForQ}. We go through each of the terms and explain why applying $\mc L^r$ produces one of the two bounds on the right hand side of \eqref{eq:MainFiniteEst}.\\
    For the term  $-Ga\cos \vartheta$, the term where only $\p_\vartheta$ derivatives are applied is exactly cancelled. Every other term uses  $\lambda^{-1}W\cdot\nabla_x$ at least once. Spatial derivatives cost at most $\mu$, derivatives of $a$ produce $a^{1-\eps}$ and hence by \eqref{eq:EstW} this term is $O((\mu/\lambda)a^{1-\eps})$.\\
    For  $E$ it is clear by \eqref{eq:EstForB} that
     \[
    \mc L^r E=O((\mu/\lambda)a^{1-\eps}).
    \]
    The final term $q_V/\lambda$ is independent of $\vartheta$. Terms in which every application of $\mc L$ uses a spatial derivative are $O((\mu/\lambda)^{r+1})$. For the other terms, any nonzero term which involves a $\vartheta$-derivative must differentiate $W$, which produces by \eqref{eq:EstW} the term $a^{1-\eps}$. Such a term must then also have a spatial derivative, so it is $O((\mu/\lambda)^2a^{1-\eps})$ which fits into the first error bound.\\
    \textbf{Step 3: Creation of nondegeneracy.} Suppose we are in the situation of \ref{it2:type}. By \eqref{eq:acomprhochi}, we have $a(x)\geq C^{-1}\rho_*>0$, with $C$ independent of $\mu$. Apply \eqref{eq:MainFiniteEst} with $r=0,1$. Since $a$ is uniformly bounded, we have
     \[
    Q=-Ga\cos \vartheta+O(\mu/\lambda),\quad \mc L Q=Ga\sin \vartheta+O(\mu/\lambda).
    \]
    Since $\max\{|\sin \vartheta|,|\cos \vartheta|\}\geq 1/\sqrt{2}$, we have
    \[
    \max\{|Q|,|\mc LQ|\}\geq C^{-1}\rho_*.
    \]
    Hence, from \eqref{ineq:Lrtaur} we have completed the proof of \ref{it2:type}.\\
    \textbf{Step 4: Preservation of nondegeneracy.} Let $\eps =1/(4(m+2))$. We distinguish between two cases.\\
    \textit{Case 1: }Suppose that $a(x)\leq \lambda^{-m-3/2}$. From the estimates \eqref{eq:BoundsRandVY} we have
     \[
    |D^j (X_\lambda - V)|\leq C\lambda^j a^{1-\eps},\quad 0\leq j\leq m+1.
    \]
    Hence, from the definition , we obtain
     \[
    \begin{aligned}
        \max_{|I|\leq m}
        |H_I^{X_\lambda}q_{X_\lambda}-H_I^Vq_V|
        &\leq C\lambda^{m+1}a^{1-\eps}
         \leq C\lambda^{-1/4}\\
        &=C\mu^{-m}\lambda^{m/\beta-1/4}
         \leq C\mu^{-m}\lambda^{-1/8}
         =o(\mu^{-m}),
    \end{aligned}
    \]
    where we used $\eps(m+3/2)<1/4$ and  $m/\beta\leq1/8$.\\
    Thus, if
    \[
    \tau_m(V)(x)\geq c_0 \mu^{-m}
    \]
    it follows that
    \[
    \tau_m(X_\lambda)(x)\geq \frac{1}{2}c_0\mu^{-m},
    \]
    which proves \ref{it1:type} in this case.\\
    \textit{Case 2: }Suppose that $a(x)\geq \lambda^{-m-3/2}$. Apply \eqref{eq:MainFiniteEst} with $r=m+1$, $r=m+2$. The bound in this case is
    \[
    C\left(\frac{\mu}{\lambda}a^{-\eps}+\left(\frac{\mu}{\lambda}\right)^{m+2}a^{-1}\right)\leq C(\lambda^{-1/2}+\lambda^{-3/8})\to 0,\quad \text{as }\lambda\to \infty.
    \]
    Hence, similar to Step 3, we get that
    \[
    \max\{|\mc L^{m+1}Q|,|\mc L^{m+2}Q|\}\geq C^{-1}a\geq C^{-1}\lambda^{-m-3/2}.
    \]
    Again from \eqref{ineq:Lrtaur} we prove \ref{it1:type} in this case. This completes the proof.
\end{proof}

\begin{proof}[Proof of Proposition \ref{prop:LocalPert}]
Pick $\gamma$ as in Lemma \ref{lem:LocalPert} and consider $X_\lambda$. Hence, the property for the support, \eqref{eq:LocalPertError} and \eqref{eq:LocalPertDerivatives} follow directly from Lemma \ref{lem:LocalPert}. Its proof gives \eqref{eq:LocalPertC0}, with a universal constant, and moreover shows that $X_\lambda$ is nowhere vanishing for sufficiently large $\lambda=\mu^\beta$, since $\mu/\lambda\to0$.
Apply Lemma \ref{lem:FiniteTypeAfterPer} \ref{it1:type} to $X_\lambda$ for each $j=0,\ldots,m$, with $c_0=\kappa$. Since $\beta\geq8(m+2)\geq8(j+2)$, we can do this. Taking the largest of the finitely many thresholds and the smallest of the positive lower constants gives \eqref{eq:LocalPertPreservation}. Lemma \ref{lem:FiniteTypeAfterPer} \ref{it2:type} gives \eqref{eq:LocalPertCreation}.
\end{proof}

\section{Nash Moser Iteration}
\label{sec:NM}
The next proposition is the central correction mechanism of our construction. In its proof we follow closely the ideas of Moser \cite{Moser1961}. See also \cite{Moser1966}.
\begin{proposition}
    \label{prop:NashMoser}
    Fix $m\in\mb N_0$. Let  $X_* \in C^\infty(N;\mathbb R^3)$ and write $r=|X_*|^2-1$. Suppose that for $d_*,\kappa>0$ we have
     \[
    \dist_\sigma(\supp r,\p N)\geq d_*,\quad \tau_m(X_*)\geq \kappa,\quad \text{on }N.
    \]
     Let $p=2m+2$ . There are  an integer $s=s(m)>p$ and constants $c,P>0$ depending only on $m$ such that if  $||X_*||_s<\infty$ and
    \begin{equation}
        \label{eq:smallnessOfT}
        ||r||_m\leq c\bigl(1+||X_*||_s+\kappa^{-1}+d_*^{-1}\bigr)^{-P},
    \end{equation}
    then there exists  $A \in C^p(N;\mathbb R^3)$ with $||A||_{p,1}<\infty$ such that
    \[
    |X_*+\curl A|=1\quad \text{on }N,\qquad A=0\quad \text{where }\dist_\sigma(x,\p N)<d_*/2.
    \]
\end{proposition}
For the proof of Proposition \ref{prop:NashMoser} we first have to construct a suitable inverse for the linearized equation. Write for $X \in C^\infty(N;\mathbb R^3)$
\[
P_X B =2X\cdot\curl B,\quad B \in C^\infty(N;\mathbb R^3),
\]
which is the linearization at $0$ of $A \mapsto |X+\curl A|^2$, $A  \in C^\infty(N;\mathbb R^3)$. The construction of the right inverse below uses exactly the idea we described prior to Definition \ref{def:Nondegen}. We restate the lemma for the convenience of the reader.

\begin{lemma}
\label{lem:inverse}
Suppose that $X \in C^\infty(N;\mathbb R^3)$ is nowhere vanishing and $\tau_m(X)>0$ in $N$. Then there exists a differential right inverse $L_X$ of order $m$ for $P_X$. In particular,
\[
    \supp L_X f\subset\supp f.
\]
The coefficients of $L_X$ depend smoothly on the derivatives of $X$ through order $2m+1$. For every $k\ge0$ and every $f\in C^\infty(N)$,
\begin{equation}
    \label{eq:inverseBound}
    ||L_X f||_{k,1}\le C_k\left(1+||X||_{2m+1+k}+|||X|^{-1}||_0+||\tau_m(X)^{-1}||_0\right)^{P_k}||f||_{m+k}.
\end{equation}
The constants $C_k,P_k>0$ depend only on $m$ and $k$.
\end{lemma}

\begin{proof}
Write in the terminology of Definition \ref{def:Nondegen} for simplicity $q=q_X$ and $H_a=H_a^X$. We first construct, for every multi-index $I$, a differential operator $L_I$ such that
\begin{equation}
    \label{eq:IndLI}
     P_X L_I f=(H_Iq)f.
\end{equation}
First, for $I=\emptyset$ define
 \[
    L_\emptyset f=\frac{\sigma f}{2}X.
\]
From the product rule we have
\begin{equation}
    \label{eq:InvIdeOd0}
    P_X L_\emptyset f =X\cdot\left(\nabla(\sigma f)\times X+\sigma f\curl X\right)
    =qf,
\end{equation}
which proves \eqref{eq:IndLI} in the case $I=\emptyset$. Define next
\[
    B_a=\frac{X\times H_a}{2|X|^2}.
\]
Since $H_a\perp X$, we have  $2B_a\times X=H_a$. The product rule therefore gives, for every $u \in C^\infty(N)$ that
 \begin{equation}
    \label{eq:NMProdRuleInv}
    P_X(uB_a)=H_a u+uP_X B_a.
\end{equation}
Suppose that for the sake of induction that  $P_X L_I f=(H_Iq)f$. Applying \eqref{eq:NMProdRuleInv} with $u=(H_Iq)f$ produces the desired term $(H_aH_Iq)f$, together with a term divisible by $H_Iq$. The latter can be cancelled using $L_I$. Accordingly, set
\[
    L_{aI}f=(H_Iq)fB_a-L_I\bigl(H_af+fP_X B_a\bigr).
\]
Indeed from \eqref{eq:NMProdRuleInv} and the induction hypothesis
 \begin{align*}
    P_X L_{aI}f=H_a\bigl((H_Iq)f\bigr)+(H_Iq)fP_X B_a-(H_Iq)\bigl(H_af+fP_X B_a\bigr)
    =(H_aH_Iq)f.
\end{align*}
This completes the induction.\\
Since
 \[
    \tau_m(X)^2=\sum_{|I|\leq m}(H_Iq)^2>0,
\]
we can now combine the $L_I$ to obtain the right inverse as
 \begin{equation}\label{eq:inverse}
    L_X f
    =\sum_{|I|\leq m}
       L_I\left(\frac{(H_Iq)f}{\tau_m(X)^2}\right).
\end{equation}
By construction,
 \[
    P_X L_X f
    =\sum_{|I|\leq m}
       \frac{(H_Iq)^2}{\tau_m(X)^2}f
    =f.
\]
The recursion shows that $L_I$ has order $|I|$. Hence  $L_X$ has order $m$ and, being a differential operator, satisfies
 \[
    \supp L_X f\subseteq\supp f.
\]
We next verify the dependence of its coefficients on  $X$. The coefficients of $L_\emptyset$ involve no derivatives of $X$. For a multiindex $I$, $H_Iq$ involves at most $|I|+1$ derivatives of $X$, while  $P_X B_a$ involves at most one derivative. Consequently, the recursion shows inductively that the coefficients of $L_I$ involve at most $|I|$ derivatives of $X$. The factor
 \[
    \frac{H_Iq}{\tau_m(X)^2}
\]
involves at most $m+1$ derivatives of $X$. Applying $L_I$ differentiates it at most $|I|$ times, so the coefficients of  $L_X$ depend smoothly on derivatives of  $X$ through order $m+1+|I|\leq 2m+1$.\\
Next we want to estimate $||L_X f||_{k,1}$. Writing $\sigma^{-1}L_Xf$ in terms of the weighted derivatives $\sigma^{|\alpha|}D^\alpha f$ by the construction above its coefficients are rational expressions in the weighted derivatives of $X$ through order $2m+1$, with denominators consisting of powers of $|X|$ and $\tau_m(X)$. Since derivatives of the weight are uniformly bounded by \eqref{eq:DerBoundSig}, we can ignore its contribution. The product and chain rules for the weighted derivatives therefore give \eqref{eq:inverseBound}.\\
\end{proof}
We also record the following consequence we use below. If $R\geq1$ and
\[
    ||X||_{2m+1+h}\leq M_hR^h\qquad(0\leq h\leq k),
\]
then counting the total number of derivatives falling on the coefficients of $L_Xf$ gives
\begin{equation}\label{eq:NMbootstrapCoeff}
    ||L_X f||_{k,1}
    \leq C_k\sum_{h=0}^k R^{k-h}||f||_{m+h}.
\end{equation}
From the structure of the coefficients described above it is clear that $C_k$ is bounded by a polynomial in $M_0,\ldots,M_k$, $|||X|^{-1}||_0$ and $||\tau_m(X)^{-1}||_0$.
\begin{proof}[Proof of Proposition \ref{prop:NashMoser}]
 Set $s=15p+2$ and write
\[
F(A)=|X_*+\curl A|^2-1.
\]
We write here
\[
F'(A)B:=P_{X_*+\curl A}B.
\]
The constants $C$, $C_k$ used below may change from line to line. These constants are bounded by polynomials in $||X_*||_s$, $\kappa^{-1}$, and $d_*^{-1}$, with coefficients and degrees depending only on $p$ and the  fixed geometry and weight. Throughout the proof, we need to choose $R_0$ sufficiently large at several steps. These  are finitely many polynomial lower bounds, so we may choose
\[
    R_0=C(1+||X_*||_s+\kappa^{-1}+d_*^{-1})^{P_1}\geq2
\]
to satisfy them all.\\
\textbf{Step 1: Estimates for the correction.} By \eqref{eq:smallnessOfT},  $|X_*|\ge1/2$. The quantities defining $\tau_m$ depend smoothly on the first $m+1$ derivatives of a  vector field. Since $p=2m+2\ge m+2$, we may choose $R_0$ so that
 \begin{equation}\label{eq:NMNeighbourhood}
||A||_{p,1}\le R_0^{-1}
\quad\Longrightarrow\quad
|X_*+\curl A|\ge\tfrac14,
\qquad
\tau_m(X_*+\curl A)\ge\tfrac\kappa2.
\end{equation}
It is not difficult to see that \eqref{eq:NMNeighbourhood} can be achieved by choosing $R_0$ larger than some lower bound which is polynomial in  $||X_*||_s,\kappa^{-1}$.\\
Suppose for the moment that  $A \in C^\infty(N;\mathbb R^3)$ and for some $R\ge1$,
 \[
||A||_{p+h,1}\le R_0^{-1}R^h
\qquad(0\le h\le14p).
\]
Set $X=X_*+\curl A$. We claim that
 \begin{equation}\label{eq:NMCorrectionBounds}
||L_X F(A)||_{0,1}\le C||F(A)||_m,
\qquad
||L_X F(A)||_{14p,1}\le CR^{14p}.
\end{equation}
The first estimate is \eqref{eq:inverseBound} with $k=0$, since $2m+1=p-1$ and \eqref{eq:NMNeighbourhood} controls the $C^0$ norms of  $|X|^{-1}$ and $\tau_m(X)^{-1}$.\\
For the second estimate, since $m\leq p-1$, the product rule gives
\begin{equation}
    \label{eq:NMbootstrapOm}
    ||X||_{p-1+h}+||F(A)||_{m+h}\le C_hR^h
    \qquad(0\le h\le14p).
\end{equation}
Taking $f=F(A)$ and $k=14p$ in \eqref{eq:NMbootstrapCoeff} gives the second estimate, with constants polynomial in $||X_*||_s$ and $\kappa^{-1}$ by \eqref{eq:NMNeighbourhood}.
This proves \eqref{eq:NMCorrectionBounds}. Note that the highest derivative of $X_*$ we used here has order $p-1+14p=15p-1<s$.\\
We collect two more estimates. First, we have directly from the definition
 \begin{equation}\label{eq:NMErrorIdentity}
F(A+B)=F(A)+F'(A)B+|\curl B|^2.
\end{equation}
Second, since $m+1=p/2$,  standard product estimates give
\begin{equation}\label{eq:NMErrorBounds}
||F'(A)B||_m\le C||B||_{p/2,1},
\qquad
|||\curl B|^2||_m\le C||B||_{p/2,1}^2,
\end{equation}
whenever  $||A||_{p,1}\le R_0^{-1}$.\\
\textbf{Step 2: Iteration.} We apply Lemma \ref{lem:Smoothing} with $\bar k=15p$ and $d=d_*/2$. The choices $(r,s)=(0,k)$ and
$(r,s)=(p/2,27p/2)$ give, for $R\ge C(1+d_*^{-1})$ and  $B \in C^\infty(N;\mathbb R^3)$ with $\supp B \subset \{x:\dist_\sigma(x,\partial N)\ge d_*/2\}$,
\begin{equation}\label{eq:NMSmoothing}
\begin{aligned}
||S_RB||_{k,1}
 &\le C_kR^k||B||_{0,1}
 &&(0\le k\le15p),\\
||(I-S_R)B||_{p/2,1}
 &\le CR^{-27p/2}||B||_{14p,1},\\
\supp S_RB
 &\subseteq\{x:\dist_\sigma(x,\supp B)\le C/R\}.
\end{aligned}
\end{equation}
Note that $p/2+27p/2=14p\le15p$.\\
We write
 \begin{equation}\label{eq:NMIteration}
\begin{gathered}
R_{j+1}=R_j^{3/2},\qquad A_0=0,\\
B_j=-L_{X_*+\curl A_j}F(A_j),\qquad
A_{j+1}=A_j+S_{R_{j+1}}B_j.
\end{gathered}
\end{equation}
We choose $R_0$ so that $R_1=R_0^{3/2}\ge C(1+d_*^{-1})$ and will inductively justify these definitions.\\
By Lemma \ref{lem:inverse}, $\supp B_j\subseteq\supp F(A_j)\subseteq\supp r\cup\supp A_j$. At every stage where the iteration is defined, \eqref{eq:NMSmoothing} therefore bounds the total support enlargement from $\supp r$ by $C\sum_{i=1}^jR_i^{-1}$. We choose $R_0$ so that
\begin{equation}\label{eq:NMSupportSum}
C\sum_{i=1}^\infty R_i^{-1}
\le2CR_0^{-3/2}<d_*/2.
\end{equation}
Since $\dist_\sigma(\supp r,\partial N)\ge d_*$, it follows that
\begin{gather}
\label{eq:NMbetasupp}
\supp B_j \subset \{x:\dist_\sigma(x,\partial N)\ge d_*/2\},\\
\label{eq:NMAlphaVanish}
A_j(x)=0 \quad\text{if}\quad \operatorname{dist}_{\sigma}(x,\partial N)<\frac{d_*}{2}.
\end{gather}
Thus Lemma \ref{lem:Smoothing} applies at the next step.
We now prove inductively that the right inverse defining $B_j$ in \eqref{eq:NMIteration} exists and it holds
\begin{equation}\label{eq:NMInduction}
||A_j||_{p+k,1}\le R_0^{-1}R_j^k
\quad(0\le k\le14p),
\qquad
||F(A_j)||_m\le R_j^{-4p}.
\end{equation}
At $j=0$, this follows from  $A_0=0$ and $||r||_m\le R_0^{-4p}$, which follows from \eqref{eq:smallnessOfT} with $P=4pP_1$ and $c=C^{-4p}$.
Suppose \eqref{eq:NMInduction} holds up to index $j$. By \eqref{eq:NMNeighbourhood}, the inverse defining  $B_j$ exists, and \eqref{eq:NMbetasupp} justifies smoothing. From \eqref{eq:NMCorrectionBounds} and \eqref{eq:NMSmoothing},
 \begin{equation}\label{eq:NMIncrementBounds}
\begin{aligned}
||S_{R_{j+1}}B_j||_{p/2,1}
 &\le CR_{j+1}^{p/2}||F(A_j)||_m
 \le CR_j^{-13p/4},\\
||S_{R_{j+1}}B_j||_{p+k,1}
 &\le CR_{j+1}^{p+k}||F(A_j)||_m
 \le CR_{j+1}^{k-5p/3}
 \qquad(0\le k\le14p).
\end{aligned}
\end{equation}
The right inverse gives  $F(A_j)+F'(A_j)B_j=0$. Using \eqref{eq:NMErrorIdentity}, we therefore have
 \[
F(A_{j+1})
=F'(A_j)(S_{R_{j+1}}-I)B_j
 +|\curl(S_{R_{j+1}}B_j)|^2.
\]
Consequently,
 \begin{align}
||F(A_{j+1})||_m
&\le C\left(||(S_{R_{j+1}}-I)B_j||_{p/2,1}
          +||S_{R_{j+1}}B_j||_{p/2,1}^2\right)\nonumber\\
&\le C\left(R_{j+1}^{-27p/2}||B_j||_{14p,1}
          +R_{j+1}^{p}||B_j||_{0,1}^2\right)\nonumber\\
&\le C\left(R_{j+1}^{-27p/2}R_j^{14p}
          +R_{j+1}^{p}||F(A_j)||_m^2\right)\nonumber\\
&=C\left(R_j^{-25p/4}+R_{j+1}^{p}||F(A_j)||_m^2\right),
\label{eq:NMErrorRecurrence}
\end{align}
where we applied first \eqref{eq:NMErrorBounds}, then the smoothing estimates \eqref{eq:NMSmoothing}, then the correction bounds \eqref{eq:NMCorrectionBounds}, and finally $R_{j+1}=R_j^{3/2}$.
Using \eqref{eq:NMInduction} and $R_{j+1}=R_j^{3/2}$ again, we obtain
\begin{align*}
||F(A_{j+1})||_m
&\le C\left(R_j^{-13p/2}+R_j^{-25p/4}\right)\\
&=C\left(R_j^{-p/2}+R_j^{-p/4}\right)R_j^{-6p}
 \le R_{j+1}^{-4p}.
\end{align*}
By the definition of $R_j$, the last inequality holds for all $j$ once $C(R_0^{-p/2}+R_0^{-p/4})\le1$.\\
For the remaining bound in \eqref{eq:NMInduction}, sum \eqref{eq:NMIncrementBounds} to obtain
 \[
||A_{j+1}||_{p+k,1}
\le CR_{j+1}^k\sum_{i=0}^jR_{i+1}^{-5p/3}
\le CR_0^{-5p/2}R_{j+1}^k
 \le R_0^{-1}R_{j+1}^k.
\]
Here we used $R_1=R_0^{3/2}$ and chose $R_0$ so that $CR_0^{-5p/2}\le R_0^{-1}$. This proves \eqref{eq:NMInduction} for each $j$. Since $p\ge2$, all the restrictions on $R_0$, including \eqref{eq:NMSupportSum}, are polynomial lower bounds in the original data. Thus as explained above at the start of the proof, one choice of  $C,P_1$, and hence of $c,P$ in \eqref{eq:smallnessOfT}, satisfies them all.\\
\textbf{Step 3: Conclusion.} The case $k=0$ of \eqref{eq:NMIncrementBounds} gives
\[
\sum_{j=0}^\infty ||A_{j+1}-A_j||_{p,1}=\sum_{j=0}^\infty ||S_{R_{j+1}}B_j||_{p,1}
\le C\sum_{j=0}^\infty R_{j+1}^{-5p/3}<\infty.
\]
It follows that there exists a potential $A$ on $N$ with  $||A||_{p,1}<\infty$ such that
\[
||A_{j}-A||_{p,1}\to 0,\quad \text{as }j\to \infty.
\]
It is easy to see that this convergence also holds in the ordinary $C^p$ norm on compact subsets of $N$. In particular,  $A\in C^p(N;\mathbb R^3)$.\\
By \eqref{eq:NMInduction},  $||F(A_j)||_m\to0$. Passing to the limit in  $F(A_j)$ gives 
\[
|X_*+\curl A|=1
\]
Finally, \eqref{eq:NMAlphaVanish} and the convergence  $A_j \to A$ show that
\begin{equation*}
    A(x)=0 \quad\text{if}\quad \operatorname{dist}_{\sigma}(x,\partial N)<\frac{d_*}{2},
\end{equation*}
which completes the proof.
\end{proof}

\section{Proof of the main theorem}
    \label{sec:ProofMain}
    Let $X_0$ be the vector field constructed in Proposition \ref{prop:Prep} and recall that
     \[
        \Delta=1-|X_0|^2.
    \]
    In particular, we have $c_0\leq|X_0|\leq 1$ for some $c_0>0$.
\begin{lemma}
    \label{lem:CorrectionDirection}
    Let $\Delta$ and $c_0$ be as above. There exist families $\mathcal F_1,\ldots,\mc F_{J_0}$ of pairwise disjoint balls from Lemma \ref{lem:BoundCharts}, each equipped with a direction $\xi_\nu\in\{\mathbf e_1,\mathbf e_2\}$, and $b_0>0$ such that the following holds.\\
    For every $J\in\mathbb N$ and every $Z\in C^\infty(N;\mathbb R^3)$ with
    \[
        \frac{c_0}{2}\leq |Z|\leq2,
    \]
    there exist nonnegative $\rho_1,\ldots,\rho_{J_0}\in C^\infty(N)$ such that $\rho_\nu \circ \Phi_a \in C_c^\infty(B_2(0))$ and
    \begin{equation}
        \label{eq:AmplitudeDecomposition}
        \sum_{\nu=1}^{J_0}\rho_\nu^2=\frac{\Delta}{J}
        \quad\text{on }N.
    \end{equation}
    For each $\nu$, we have
    \[
        \supp\rho_\nu\subset\supp\Delta\cap
        \bigcup_{U_a(2)\in\mathcal F_\nu}U_a(2),
        \qquad
        |\Pi_{\xi_\nu^\perp}Z|\geq2b_0
        \quad\text{on }\supp\rho_\nu.
    \]
    Moreover, if $\mu\geq2$ and $||Z||_k\leq M_k\mu^k$ for every $k\geq0$, then
    \begin{equation}
        \label{eq:AmplitudeMuBound}
        \sigma^k|D^k\rho_\nu|\leq C_{k,\eps}\mu^k\rho_\nu^{1-\eps},
        \qquad k\geq0,\quad0<\eps<1.
    \end{equation}
    The constants $C_{k,\eps}$ depend only on $\Delta$ and on $(M_k)_k$.
\end{lemma}

\begin{proof}
    From Proposition \ref{prop:Prep} we have that $\dist_\sigma(\supp \Delta, \p N)=d_*>0$.
    Thus, by Lemma \ref{lem:BoundCharts} there exist balls $U_i(2)$ and a partition $I=\bigsqcup_{q=1}^Q I_q$ such that the sets $U_i(1)$ cover a neighborhood of $\supp\Delta$, and the sets $U_i(2)$ with $i\in I_q$ are pairwise disjoint. For $q=1,\ldots,Q$ and $j=1,2$, enumerate
    \[
        \mathcal F_{2(q-1)+j}=\{U_i(2):i\in I_q\},\qquad \xi_{2(q-1)+j}=\mathbf e_j,
        \qquad J_0=2Q.
    \]
    For $Z$ as above, we have
    \begin{equation}
        \label{eq:QuantitativeIndependence}
        \sum_{j=1}^2|\Pi_{\mathbf e_j^\perp}Z|^2
        =|Z|^2+Z_3^2\geq\frac{c_0^2}{4}.
    \end{equation}
     Hence at least one of the summands  is bounded below by  $c_0^2/8$. Fix $b_0>0$ such that  $3b_0<c_0/(2\sqrt2)$. Let $\chi_i$ be a smooth cutoff, with  $\chi_i\equiv1$ on $U_i(1)$ and $\supp\chi_i\subset U_i(2)$, such that $\chi_i\circ \Phi_i$ is a fixed cutoff on $B_2(0)$. Choose a cutoff $h\in C^\infty([0,\infty))$ such that
     \[
        h\equiv0\quad\text{on }[0,(2b_0)^2],
        \qquad h\equiv1\quad\text{on }[(3b_0)^2,\infty).
    \]
     Define
    \[
        v_\nu(Z)=\sum_{U_i(2)\in\mathcal F_\nu} \chi_i h(|\Pi_{\xi_\nu^\perp}Z|^2).
    \]
    extending each summand by zero outside $U_i(2)$. As $(U_i(2))_i$ is locally finite,  these sums are smooth. At every point of $U_i(1)$, one of the corresponding cutoffs $v_\nu(Z)$ equals $1$
    \[
    \sum_{\nu=1}^{J_0}v_\nu(Z)^2 \geq 1,\quad \text{on }\bigcup_i U_i(1).
    \]
    Set
    \[
        \rho_\nu=\sqrt{\frac{\Delta}{J}}\,\frac{v_\nu(Z)}{\sqrt{\sum_{\nu=1}^{J_0}v_\nu(Z)^2}},\quad \text{on }\bigcup_i U_i(1),
    \]
    and extend by zero to $N$. Since $\sqrt\Delta$ is smooth and $\supp\Delta$ is contained in the chosen neighborhood, the extension is smooth. We easily have now \eqref{eq:AmplitudeDecomposition}, and the support and projection bounds follow from the choice of $\chi_i$ and $h$.\\
    Finally, suppose that $||Z||_k \le M_k\mu^k$. Clearly this gives
    \[
    |||\Pi_{\xi_\nu^\perp}Z|^2||_k\le C_k \mu^k.
    \]
    If one picks the cutoffs defining $\chi_i$ and $h$ using the standard smooth function $e^{-1/t}$ for $t>0$, extended by zero for $t\leq0$, they satisfy \eqref{eq:AmplitudeBounds} with $\mu=1$ and uniform constants. The chain and product rules  together with $\sigma\simeq\ell_i$ on $U_i(2)$ give
    \[
        \sigma^k|D^kv_\nu(Z)|\leq C_{k,\eps}\mu^kv_\nu(Z)^{1-\eps}.
    \]
    Combining these bounds with Proposition \ref{prop:Prep}\,\ref{prop:it4} and the Leibniz rule we easily have
    \[
        \sigma^k|D^k\rho_\nu| \leq C_{k,\eps}\mu^kJ^{-\eps/2}\rho_\nu^{1-\eps}
        \leq C_{k,\eps}\mu^k\rho_\nu^{1-\eps}.
    \]
    This proves \eqref{eq:AmplitudeMuBound}.
\end{proof}
We briefly describe how we add the perturbations to $X_0$. We work in \textit{steps} and \textit{stages}. A step consists of an application of Proposition \ref{prop:LocalPert}. We group $J_0$ steps into a stage, where at the start of each stage, we recompute the functions $\rho_{\nu}$ using the vector field from the start of the stage. Hence by increasing the number of stages $J$ we can ensure that each step makes a small enough perturbation such that the smallness condition of Proposition \ref{prop:LocalPert} is satisfied. Moreover, this make sure that the directions $\xi_\nu$ are sufficiently transversal with respect to the current field.
\[
\underbrace{
X_0\xrightarrow{\rho_1}X_1
\longrightarrow\cdots
\xrightarrow{\rho_{J_0}}X_{J_0}
}_{\substack{\text{stage }0\\
             \text{$\rho_\nu$ from }X_0}}
\underbrace{
\xrightarrow{\rho_{J_0+1}}X_{J_0+1}
\longrightarrow\cdots
\xrightarrow{\rho_{2J_0}}X_{2J_0}
}_{\substack{\text{stage }1\\
             \text{$\rho_\nu$ from }X_{J_0}}}
\longrightarrow\cdots\longrightarrow X_*
\]
\begin{proposition}[Stage]
    \label{prop:Stage}
    Let $J_0$ be as in Lemma \ref{lem:CorrectionDirection}. There exists $J_*\in\mathbb N$ such that the following holds. Let $J\in\mathbb N$ with $J\geq J_*$. Fix $\kappa>0$ and $m,L\in\mathbb N_0$ with
    \[
        m\geq2J_0,\qquad L\geq2m.
    \]
    Suppose that $Z\in C^\infty(N;\mathbb R^3)$ satisfies for some $\mu \ge 2$
    \[
        \frac{c_0}{2}\leq |Z|\leq2\quad\text{on }N,\qquad
        ||Z||_k\leq M_k\mu^k,\quad k\geq0.
    \]
    Set $\beta=8(m+2)$ and $\lambda=\mu^{\beta^{J_0}}$.\\
    For $\mu$ large enough, there exists
    \[
        B\in C^\infty(N;\mathbb R^3),\qquad \supp B\subset\supp\Delta,
    \]
    such that $\widetilde Z=Z+\curl B$ is nowhere vanishing and
    \begin{equation}
        \label{eq:StageError}
        |\widetilde Z|^2=|Z|^2+\frac{\Delta}{J}+E,
        \qquad
        ||E||_m\leq C_{L,m}\mu^{\beta(m-(L+1)/2)}.
    \end{equation}
    Moreover, for every $k\geq0$,
    \begin{equation}
        \label{eq:StageDerivatives}
        ||\widetilde Z||_k\leq C_{L,k}\lambda^k,
    \end{equation}
    and, for every $0\leq j\leq m-2J_0$ and $x\in N$,
    \begin{equation}
        \label{eq:StageNondegeneracy}
        \Delta(x)+\mu^j\tau_j(Z)(x)\geq\kappa
        \quad\Longrightarrow\quad
        \tau_{j+2J_0}(\widetilde Z)(x)\geq c_1\lambda^{-(j+2J_0)}.
    \end{equation}
\end{proposition}

\begin{proof}
    Let $b_0$ be as in Lemma \ref{lem:CorrectionDirection}, and let $C,\gamma$ be the constants from Proposition \ref{prop:LocalPert}. We pick $J_*\in\mb N$ large enough such that
    \[
        \frac{CJ_0}{\sqrt{J_*}}<\frac{b_0}{4},\qquad
        \frac{1}{b_0\sqrt{J_*}}\leq\gamma.
    \]
    Denote
    \[
        \lambda_0=\mu,\qquad \lambda_{\nu+1}=\lambda_\nu^{\beta}, \quad 0\leq\nu<J_0,
    \]
    so that $\lambda_1=\mu^\beta$ and $\lambda_{J_0}=\lambda$.\\
    Apply Lemma \ref{lem:CorrectionDirection} with $Z$ and $J$ to obtain $\rho_1,\ldots,\rho_{J_0}$ with
    \begin{equation}
        \label{eq:SumRhoBounds}
        0\leq\rho_\nu\leq J^{-1/2},\qquad
        \sum_{\nu=1}^{J_0}\rho_\nu^2=\frac{\Delta}{J}.
    \end{equation}
    Moreover, the functions $\rho_\nu$ satisfy \eqref{eq:AmplitudeMuBound}, and
    \[
        |\Pi_{\xi_\nu^\perp}Z|\geq2b_0
        \quad\text{on }\supp\rho_\nu.
    \]
    Set $Z_0=Z$. We construct $Z_1,\ldots,Z_{J_0}$ by successive applications of Proposition \ref{prop:LocalPert}. Suppose that, for $0\leq\nu<J_0$, we constructed nowhere vanishing $Z_0,\ldots,Z_{\nu}$ such that
    \begin{equation}
        \label{eq:StageIndDerivatives}
        ||Z_r||_k\leq C_{r,k,L}\lambda_r^k,
        \quad k\geq0,\quad 0\leq r\leq\nu,
    \end{equation}
    and for $1\leq r\leq\nu$
    \begin{align}
        \label{eq:StageIndPotential}
        Z_r&=Z_{r-1}+\curl B_r,\qquad
        \supp B_r\subset\supp\rho_r\subset\supp\Delta,\\
        \label{eq:StageIndIncrement}
        ||Z_r-Z_{r-1}||_0&\leq\frac{C}{\sqrt J}
        +C_{r,L}\frac{\lambda_{r-1}}{\lambda_r},\\
        \label{eq:StageIndLength}
        |Z_r|^2&=|Z_{r-1}|^2+\rho_r^2+E_r,\\
        \label{eq:StageIndError}
        ||E_r||_m&\leq C_{r,L}\lambda_r^m
        \left(\frac{\lambda_{r-1}}{\lambda_r}\right)^{L+1}.
    \end{align}
    For $\nu=0$, this holds by the assumptions on $Z_0=Z$.\\
    We verify that we can apply Proposition \ref{prop:LocalPert} at the next step. Summing \eqref{eq:StageIndIncrement} and using $\beta\geq2$, we have
    \begin{equation}
        \label{eq:ItBounds}
        ||Z_{\nu}-Z||_0
        \leq\sum_{r=1}^{\nu}||Z_r-Z_{r-1}||_0
        \leq\frac{CJ_0}{\sqrt J}+C_L\lambda_1^{-1/2}.
    \end{equation}
    Hence, for sufficiently large $\mu$, the choice of $J_*$ gives
    \[
        |\Pi_{\xi_{\nu+1}^\perp}Z_{\nu}|
        \geq2b_0-||Z_{\nu}-Z||_0\geq b_0
        \quad\text{on }\supp\rho_{\nu+1}.
    \]
    In particular, $|\Pi_{\xi_{\nu+1}^\perp}Z_{\nu}|\ge b_0/2$ near $\supp\rho_{\nu+1}$, which proves \eqref{eq:Transversality}. By \eqref{eq:SumRhoBounds} and the choice of $J_*$, we have $\rho_{\nu+1}/|\Pi_{\xi_{\nu+1}^\perp}Z_{\nu}|\leq1/(b_0\sqrt J)\leq\gamma$ which is \eqref{eq:SmallnessLocPert}. The bounds \eqref{eq:DataMuBound} for $Z_{\nu}$ and $\rho_{\nu+1}$ are clear. Similarly, \eqref{eq:AmplitudeBounds} follows from \eqref{eq:AmplitudeMuBound}.\\
    For every $U_i(2)\in\mathcal F_{\nu+1}$, we apply Proposition \ref{prop:LocalPert} with
    \begin{equation*}
        \mu=\lambda_{\nu},\qquad \lambda=\lambda_{\nu+1},\qquad V=Z_{\nu}\circ\Phi_i,\qquad \rho=\rho_{\nu+1}\circ\Phi_i,\qquad \xi=\xi_{\nu+1}.
    \end{equation*}
    We write $\widehat B_{\nu+1,i}$ for the potentials on $B_2(0)$. Recall that the balls within $\mathcal F_{\nu+1}$ are pairwise disjoint. Hence, extending $\ell_i\widehat B_{\nu+1,i}\circ\Phi_i^{-1}$ by zero outside $U_i(2)$ and taking the sum over $i$, we obtain $B_{\nu+1}\in C^\infty(N;\mathbb R^3)$ with
    \begin{equation}
        \label{eq:Suppbetai}
        \supp B_{\nu+1}\subset\supp\rho_{\nu+1}\subset\supp\Delta.
    \end{equation}
    By Proposition \ref{prop:LocalPert} and \eqref{eq:LocalPertDerivatives}, $Z_{\nu+1}=Z_{\nu}+\curl B_{\nu+1}$ is nowhere vanishing and
    \[
        ||Z_{\nu+1}||_k\leq C_{\nu+1,k,L}\lambda_{\nu+1}^k,
        \quad k\geq0.
    \]
    Together with \eqref{eq:Suppbetai}, this proves \eqref{eq:StageIndPotential} and \eqref{eq:StageIndDerivatives} for $r=\nu+1$. By \eqref{eq:LocalPertC0} and \eqref{eq:SumRhoBounds},
    \[
        ||Z_{\nu+1}-Z_{\nu}||_0
        \leq\frac{C}{\sqrt J}+C_{\nu+1,L}\frac{\lambda_{\nu}}{\lambda_{\nu+1}},
    \]
    which proves \eqref{eq:StageIndIncrement} for $r=\nu+1$. Moreover, by \eqref{eq:LocalPertError},
    \begin{equation}
        \label{eq:errorFormi}
        |Z_{\nu+1}|^2=|Z_{\nu}|^2+\rho_{\nu+1}^2+E_{\nu+1},
        \qquad
        ||E_{\nu+1}||_m\leq C_{\nu+1,L}\lambda_{\nu+1}^m
        \left(\frac{\lambda_{\nu}}{\lambda_{\nu+1}}\right)^{L+1}.
    \end{equation}
    Thus \eqref{eq:StageIndLength}--\eqref{eq:StageIndError} hold for $r=\nu+1$, completing the induction.
    Since we apply Proposition \ref{prop:LocalPert} only a finite number of times, the lower bound on $\mu$ required for an application can be satisfied by increasing $\mu$.\\
    Set
    \[
        B=\sum_{\nu=1}^{J_0}B_\nu,\qquad
        \widetilde Z=Z_{J_0},\qquad E=\sum_{\nu=1}^{J_0}E_\nu.
    \]
    The support property and \eqref{eq:StageDerivatives} follow immediately. Summing \eqref{eq:errorFormi} and using \eqref{eq:SumRhoBounds} gives
    \[
    |\widetilde Z|^2=|Z|^2+\frac{\Delta}{J}+E.
    \]
    Since $\beta\geq2$ and $L\geq2m$, we have
    \[
        ||E||_m\leq C_{L,m}\sum_{\nu=1}^{J_0}
        \lambda_\nu^{m-(L+1)/2}
        \leq C_{L,m}\lambda_1^{m-(L+1)/2}
        =C_{L,m}\mu^{\beta(m-(L+1)/2)},
    \]
    which proves \eqref{eq:StageError}.\\
    It remains to prove \eqref{eq:StageNondegeneracy}. Fix $0\leq j\leq m-2J_0$ and $x\in N$ with
    \begin{equation}
        \label{eq:ZNondegen}
        \Delta(x)+\mu^j\tau_j(Z)(x)\geq\kappa.
    \end{equation}
    If $\Delta(x)<\kappa/2$, \eqref{eq:ZNondegen} gives
    \[
        \tau_j(Z)(x)\geq\frac{\kappa}{2}\mu^{-j}.
    \]
    Applying \eqref{eq:LocalPertPreservation} successively at orders $j,j+2,\ldots,j+2J_0-2$ gives \eqref{eq:StageNondegeneracy} at $x$. All these orders are admissible since $j+2J_0\leq m$.\\
    If instead $\Delta(x)\geq\kappa/2$, then \eqref{eq:SumRhoBounds} gives $\nu\in\{1,\ldots,J_0\}$ such that
    \[
        \rho_\nu(x)\geq\sqrt{\frac{\kappa}{2JJ_0}}.
    \]
    By \eqref{eq:LocalPertCreation}, we have
    \[
        \tau_{j+2\nu}(Z_\nu)(x)\geq\tau_1(Z_\nu)(x)
        \geq c_*\lambda_\nu^{-(j+2\nu)}.
    \]
    Applying \eqref{eq:LocalPertPreservation} at every following step again gives \eqref{eq:StageNondegeneracy}. Similar as before since the set of admissible $j$ is finite, after possibly decreasing $c_*$ we get a constant $c_1>0$ which is uniform in $j$.
\end{proof}

\begin{proof}[Proof of Theorem \ref{thm:main}]
    Recall that $X_0$ in Proposition \ref{prop:Prep} has the following properties. There exist $c_0,\kappa_0,d_*>0$ with
     \[
    c_0\leq |X_0|\leq 1,\qquad \Delta+\tau_0(X_0)\geq\kappa_0,\quad\text{on }N.
    \]
    Moreover,
    \begin{equation}
        \label{eq:ErrorSupportd}
        \dist_{\sigma}(\supp\Delta,\p N)\geq d_*.
    \end{equation}
    Let $J_0$ be as in Lemma \ref{lem:CorrectionDirection}, and pick $J\geq J_*$ as in Proposition \ref{prop:Stage}. Denote
     \begin{equation}
        \label{eq:ParaChoice}
        N_*=JJ_0,\quad m=2N_*,\quad \beta=8(m+2).
    \end{equation}
    Let  $s=s(m)$ be as in Proposition \ref{prop:NashMoser}. Fix an integer $L\geq2m$, to be chosen sufficiently large in Step 2. For some $\Lambda \geq 1$ which we choose large enough later, we write
    \[
    \lambda_0=\Lambda^{1/\beta},\quad\lambda_{j+1}=\lambda_j^{\beta^{J_0}},\quad 0\leq j<J.
    \]
    We can assume that $\lambda_0 \geq 2$.\\
    \textbf{Step 1: Stages.} In the following we construct $X_J$ in $J$ stages by successive applications of Proposition \ref{prop:Stage}.\\
    Suppose that $X_0,X_1,\ldots,X_j$ have been constructed for some $0\leq j<J$ and satisfy, for $0\leq r\leq j$,
    \begin{align}
        \label{eq:IndHypoXrLength}
        &\frac{c_0}{2}\leq |X_r|\leq 2,\\
        \label{eq:IndHypoXrDer}
        &||X_r||_k\leq C_{r,k,L}\lambda_r^k,\quad k\geq0, \\
        \label{eq:NondegeneracyStage}
        &\Delta+\lambda_r^{2rJ_0}\tau_{2rJ_0}(X_r)\geq\kappa_r \quad\text{on }N,
    \end{align}
    and for $0\leq r<j$
    \begin{align}
        \label{eq:IndHypoLengthSplit}
        |X_{r+1}|^2&=|X_r|^2+\frac{\Delta}{J}+E_r,\\
        \label{eq:IndHypoErrorGrowth}
        ||E_r||_m&\leq C_{r,L}\lambda_r^{\beta(m-(L+1)/2)}
        \leq C_{r,L}\Lambda^{m-(L+1)/2}.
    \end{align}
    These bounds hold for $X_0$ by Proposition \ref{prop:Prep}.\\
    To obtain $X_{j+1}$ we apply Proposition \ref{prop:Stage} with $Z=X_j$, $\mu=\lambda_j$, $\kappa=\kappa_j$ and order $2jJ_0$ in \eqref{eq:StageNondegeneracy}. This gives $B_{j+1}$ with $\supp B_{j+1}\subset\supp\Delta$ such that
    \[
        X_{j+1}=X_j+\curl B_{j+1}.
    \]
    Hence $X_{j+1}$ is divergence-free. By \eqref{eq:StageDerivatives} and \eqref{eq:StageNondegeneracy}, we have
    \begin{align}
        \label{eq:XjGrowth}
        ||X_{j+1}||_k&\leq C_{j+1,k,L}\lambda_{j+1}^k,\quad k\geq0,\\
        \label{eq:XjNonDegene}
        \tau_{2(j+1)J_0}(X_{j+1})&\geq\kappa_{j+1}\lambda_{j+1}^{-2(j+1)J_0}.
    \end{align}
    In particular, this means that \eqref{eq:IndHypoXrDer}--\eqref{eq:NondegeneracyStage} hold for $X_{j+1}$.\\
    By \eqref{eq:StageError} one has
    \begin{align}
        \label{eq:LengthXj1}
        |X_{j+1}|^2&=|X_j|^2+\frac{\Delta}{J}+E_j,\\
        \label{eq:ErrorEj}
        ||E_j||_m&\leq C_{j,L}\lambda_j^{\beta(m-(L+1)/2)} \leq C_{j,L}\Lambda^{m-(L+1)/2},
    \end{align}
    where we used $L\geq2m$ and $\lambda_j^\beta\geq\Lambda$. Thus \eqref{eq:IndHypoLengthSplit}--\eqref{eq:IndHypoErrorGrowth} hold for $r=j$. It remains to check \eqref{eq:IndHypoXrLength}. The support property of $B_{j+1}$ also gives $\supp E_j\subset\supp\Delta$. Summing \eqref{eq:IndHypoLengthSplit} for $r=0,\ldots,j$, we have
    \[
        |X_{j+1}|^2=|X_0|^2+\frac{j+1}{J}\Delta+\sum_{r=0}^j E_r.
    \]
    Since $0\leq(j+1)\Delta/J\leq\Delta$, it follows that
    \[
        c_0^2-\sum_{r=0}^j||E_r||_0\leq |X_{j+1}|^2
        \leq1+\sum_{r=0}^j||E_r||_0.
    \]
    Thus, for sufficiently large $\Lambda$,
    \begin{equation}
        \label{eq:OmiLengthComp}
        \frac{c_0}{2}\leq |X_{j+1}|\leq2.
    \end{equation}
    Since we only have a finite number of stages, we can choose $\Lambda$ large enough such that all applications of Proposition \ref{prop:Stage} and \eqref{eq:OmiLengthComp} hold. This completes the proof of \eqref{eq:IndHypoXrLength} for $X_{j+1}$ and hence the proof of the induction.\\
    By construction we have
    \[
    \Delta_\Lambda:=|X_J|^2-1=|X_0|^2+\Delta+\sum_{j=0}^{J-1}E_j-1=\sum_{j=0}^{J-1}E_j
    \]
    and
    \begin{equation}
        \label{eq:FinalError}
        ||\Delta_\Lambda||_m\leq C_L\Lambda^{m-(L+1)/2},\quad
        \supp \Delta_\Lambda\subset\supp\Delta.
    \end{equation}
    In particular, by \eqref{eq:ErrorSupportd}
    \begin{equation}
        \label{eq:rLambdaSupp}
        \dist_{\sigma}(\supp \Delta_\Lambda,\p N)\geq d_*.
    \end{equation}
    By \eqref{eq:XjNonDegene} with $j+1=J$, together with $m=2JJ_0$
    \begin{equation}
        \label{eq:TypeCondNash}
        \tau_m(X_J)\geq c_L\lambda_J^{-m},\quad\text{on }N.
    \end{equation}
    \textbf{Step 2: Nash--Moser Correction.} With $\kappa=c_L\lambda_J^{-m}$, the bounds for $X_J$ in \eqref{eq:IndHypoXrDer}, together with $s>m$ and $\lambda_J=\Lambda^{\beta^{N_*-1}}$, give
    \[
    1+||X_J||_s+\kappa^{-1}+d_*^{-1}
    \leq C_L\lambda_J^s
    =C_L\Lambda^{s\beta^{N_*-1}}.
    \]
    Let $c,P$ be the constants from Proposition \ref{prop:NashMoser}. By \eqref{eq:FinalError}, we therefore have
    \begin{equation}
        \label{eq:PowerExp}
        c^{-1}||\Delta_\Lambda||_m\bigl(1+||X_J||_s+\kappa^{-1}+d_*^{-1}\bigr)^P\leq C_L\Lambda^{m-(L+1)/2}\lambda_J^{sP}=C_L\Lambda^{m-(L+1)/2+sP\beta^{N_*-1}}.
    \end{equation}
    Since $c,s,P,\beta,N_*$ are independent of $L$ we can pick $L$ large enough such that the exponent in \eqref{eq:PowerExp} is negative. Thus, choosing $\Lambda$ sufficiently large makes the right-hand side at most $1$, which gives the smallness condition \eqref{eq:smallnessOfT} in Proposition \ref{prop:NashMoser}. Together with \eqref{eq:TypeCondNash} and \eqref{eq:rLambdaSupp}, this allows us to apply Proposition \ref{prop:NashMoser} and obtain $A \in C^{2m+2}(N;\mathbb R^3)$ with
    \begin{equation}
        \label{eq:SupportAlpha}
        ||A||_{2m+2,1}<\infty,\quad A=0 \quad\text{on }\{x\in N:\dist_{\sigma}(x,\p N)<d_*/2\},
    \end{equation}
    such that
     \[
    X:=X_J+\curl A,\quad |X|=1.
    \]
    The  vector field $X$ is divergence-free, and $X \in C^{2m+1}(N;\mathbb R^3)$. By \eqref{eq:SupportAlpha} and the support properties of the stage perturbations, we have
     \[
    X=X_0,\quad \text{on }\p N.
    \]
    In particular, from Proposition \ref{prop:Prep} it holds that $X\cdot\nu=0$ on $\partial N$ and the vector field has the desired indices at points of $\Sigma$. The asymptotics near the singularities follow directly from the definition of the weighted norms. This also implies easily that $X \in H^1(U;\mb R^3)$: For each $p \in \Sigma$ and $R$ sufficiently small we have by the derivative bounds for $X$ that
    \[
    \int_{U \cap B_R(p)}|DX|^2 \lesssim \int_0^R r^{-2\delta} dr <\infty,
    \]
    since $0<\delta<1/2$. Since $|X|=1$ and $DX$ is bounded away from the singularities this gives $X \in H^1(U;\mb R^3)$, which completes the proof.
\end{proof}

\section{Extension to Manifolds}
We briefly describe an extension of the main result. Suppose $(M,g)$ is a smooth compact oriented Riemannian $3$-manifold with smooth boundary. Define the singular set $\Sigma \subset \p M$ exactly as in \eqref{def:Sigma} for $M$ instead of $U$. In this more geometric setting it is convenient to work with forms instead of vector fields. Suppose that $\Omega$ is a two-form on $M\setminus \Sigma$. Note that for the vector field $X=(*\Omega)^\sharp$, one has
\[
|X|_g=|\Omega|_g, \qquad \iota_\partial^*\Omega=0 \quad\Longleftrightarrow\quad  g(X,\nu) =0,\quad \text{on }\partial M\setminus\Sigma,
\]
where $\nu$ denotes the outward unit normal along $\partial M$ and $\iota_\p :\p M\setminus\Sigma \hookrightarrow M\setminus\Sigma$ the boundary inclusion. In particular, we have the equivalence: $X$ is a magnetic vector field if and only if $\Omega:=*X^\flat$ satisfies
\begin{equation}
    \label{eq:MagneticForm}
    |\Omega|_g=1,\quad d\Omega=0,\quad \iota_\p^* \Omega=0.
\end{equation}
One can prove the following result.
\begin{theorem}
    \label{thm:mainMfd}
    There exists $\Omega \in C^1(M\setminus \Sigma;\Lambda^2 T^*M)$ such that
    \[
    d\Omega=0,\quad |\Omega|_g=1,\quad \iota_\partial^*\Omega=0.
    \]
    The associated vector field $X=(*\Omega)^\sharp$ has
    \[
    \ind_{p_i}(X|_{S_i})=\chi(S_i),
    \]
    for every boundary component $S_i\subset \p M$ with $\chi(S_i)\neq 0$.
\end{theorem}
The proof of Theorem \ref{thm:mainMfd} is a purely technical and notational modification of the proof of Theorem \ref{thm:main}. To be slightly more precise, the arguments can essentially be copied if one replaces $U$ by $M$ and the Euclidean derivatives, lengths etc. by their Riemannian counterparts as well as using the correspondence between forms and vector fields above. In fact, we had first written the proof for the general manifold case and decided afterwards that the analytic ideas are easier to grasp by working on a Euclidean domain.

	\let\oldthebibliography\thebibliography
	\let\endoldthebibliography\endthebibliography
	\renewenvironment{thebibliography}[1]{
	\begin{oldthebibliography}{#1}
	\setlength{\itemsep}{0.5pt}
	\setlength{\parskip}{0.5pt}
	}
	{
	\end{oldthebibliography}
	}

	{\small
	\IfFileExists{abbrv-dz.bst}{\bibliographystyle{abbrv-dz}}{\bibliographystyle{abbrv}}
	 \bibliography{sources}
	}

\end{document}